\documentclass[11pt,letterpaper]{article}

\usepackage[margin=1in]{geometry}
\usepackage[T1]{fontenc}
\usepackage[utf8]{inputenc}
\usepackage{lmodern}
\usepackage{amsmath,amssymb,amsthm,mathtools}
\usepackage{xcolor}
\usepackage{doi}
\usepackage{booktabs}
\usepackage{array}
\usepackage{hyperref}
\usepackage{authblk}
\usepackage{microtype}
\usepackage{enumitem}
\usepackage{needspace}
\usepackage{xspace}

\hypersetup{colorlinks=true,linkcolor=blue!50!black,citecolor=blue!50!black,urlcolor=blue!50!black,
pdftitle={Dodgson Condensation and Vaporisation: Recurrences for Structured Determinants},
pdfauthor={Max A. Alekseyev and Dmitry I. Khomovsky}}

\newtheorem{theorem}{Theorem}[section]
\newtheorem{lemma}[theorem]{Lemma}
\newtheorem{proposition}[theorem]{Proposition}
\newtheorem{corollary}[theorem]{Corollary}

\newtheorem{example}[theorem]{Example}
\newtheorem{remark}[theorem]{Remark}
\theoremstyle{definition}
\newtheorem{definition}[theorem]{Definition}

\newcommand{\Ddet}{\operatorname{Ddet}}
\newcommand{\diag}{\operatorname{diag}}
\newcommand{\F}{\mathbb F}
\newcommand{\eps}{\varepsilon}
\newcommand{\coef}[1]{\mathop{[\![#1]\!]}\xspace}

\title{Dodgson Condensation and Vaporisation:\\ Recurrences for Structured Determinants}
\author[1]{Max A. Alekseyev}
\author[2]{Dmitry I. Khomovsky}
\affil[1]{\small The George Washington University, Washington, DC, USA. Email: \href{mailto:maxal@gwu.edu}{maxal@gwu.edu}}
\affil[2]{\small Independent researcher. Email: \href{mailto:dmkhomovsky@gmail.com}{dmkhomovsky@gmail.com}}
\date{}

\begin{document}
\maketitle

\begin{abstract}
Dodgson condensation can expose banded connected-minor arrays inside a
dense structured matrix.  We construct ordinary parents of specified
late descendants directly from the original entry data, give their band
and boundary entries, and derive determinant recurrences by ordinary
row--column expansion.  Ranks outside
a fixed diagonal band give generically sharp bounds for the condensation
depth and surviving width.  For separated tails, local generator
recurrences give the entries of an ordinary parent.  An explicit boundary
normalization preserves both layers of the order-two descendant.  A Toeplitz
bulk and a fixed northwest correction then give a constant-coefficient
recurrence; arbitrary position-dependent structure need not do so.
The method is worked out for rank-one and rank-two generators,
geometric and trigonometric kernels, polynomial distance, and finite
exponential sums.  Analysis of repeated eigenvalues in the ordinary transfer gives exact
generic recurrence orders $5,41,257$ for distance powers $1,3,5$.
Functional equations classify rational and polynomial kernels with
vanishing connected minors and distinguish the sampling hypotheses needed
for continuous and analytic functions.  Agreement of the two triangular
formulas on the main diagonal can further narrow the band.
The paper includes explicit parent entries, scalar recurrence
coefficients, and exact certificates for the stated minimal orders.
\end{abstract}

\section{Introduction}\label{sec:introduction}

Dodgson condensation replaces a matrix by arrays of connected minors
\cite{Dodgson1866}.  For a dense structured matrix, vanishing minors
can make these arrays banded.  We use this information to construct an
ordinary parent with a banded bulk and a northwest correction, and
then derive the determinant recurrence by ordinary row--column
expansion.  We call the inverse construction \emph{vaporisation}.
The main constructions specify the parent's entries and verify that it
and the original matrix have a common two-layer descendant.  This
identity, with the continuation convention of
Section~\ref{sec:preliminaries}, certifies equality of their determinants.

Three questions organize the paper: which entry structure forces a
band, how many condensation steps are needed, and which ordinary
ancestor converts that band into a determinant recurrence?
Section~\ref{sec:rankone-model} answers all three for
$H_{ij}=f_i g_j$ off the main diagonal and arbitrary entries on it.
One step makes the first layer tridiagonal, and two steps make both
layers tridiagonal.  We give their entries, an ordinary tridiagonal
ancestor with a northwest correction, and the resulting second-order
recurrence with position-dependent coefficients.

The general mechanism is developed in
Section~\ref{sec:structural-route}.  If the upper and lower rectangles
outside a core of semibandwidths $b_+,b_-$ have ranks at most
$r_+,r_-$, put $R:=\max(r_-,r_+)$.  The first layer is supported on
the band of semibandwidths $b_-+r_-,b_++r_+$ after $R$ steps, and both
layers have that support after $R+1$ steps
(Theorem~\ref{thm:band-core-compression}).  The bounds are sharp for
generic parameters; a converse requires additional nonvanishing
conditions on the connected minors.  The rank hypothesis includes
independent separated formulas above and below an arbitrary diagonal
band.  Products and structured similarities retain
explicit rank bounds, whereas an arbitrary permutation need not retain
bounded-depth compression.

For separated tails, Theorem~\ref{thm:local-anchored-vaporisation}
gives the entries of an ordinary banded matrix directly from local
recurrences for the separated generators.  The boundary normalization
preserves both layers of the descendant whose first layer has order
two.  We call this the \emph{order-two descendant}; its order is
distinct from the number of condensation steps used to reach it.
The constructed parent therefore computes the original determinant,
including specializations at zero intermediate minors.  The formulas
give its band and exceptional northwest entries directly from the
structural data.

A banded Toeplitz core with generators satisfying fixed shift
recurrences gives a Toeplitz bulk and a fixed northwest
correction (Corollary~\ref{cor:band-core-stationary-parent}).  More general
Toeplitz symbols with recurrent tails have the explicit construction
of Theorem~\ref{thm:anchored-recurrent-tails}.  Row--column expansion
expresses the determinant through finitely many auxiliary minors with
linear update equations.  This determinant transfer gives the binomial bound for the recurrence order.  We include
the argument incorporating the fixed northwest correction, which
determines the initial transfer values and permits repeated roots in
the examples.  Position-dependent bulk entries generally give
position-dependent transfers, and band formation alone does not imply
a constant-coefficient recurrence.
Example~\ref{ex:banded-not-constant-recurrent} demonstrates this with
polynomial entries and factorial determinants.  For periodic transfer
data, recording the position within a period gives a constant transfer
matrix.

Section~\ref{sec:constantapps} treats separated triangles, geometric
tails around a band core, polynomial generators, cosine sums,
polynomial-distance kernels, and finite exponential kernels.  Each
family is followed through its banded descendant, ordinary parent,
and recurrence.  The cosine family has generic scalar order five.
For odd polynomial distance, the repeated unit root of the ordinary
bulk companion gives sharper bounds than the general symmetric-band
bound, with exact generic orders $5,41,257$ at powers $1,3,5$.
The signed polynomial family has exact order $(p+1)^2+1$ for nonzero
diagonal parameter.  Explicit coefficients and exactly verified
initial values give usable scalar formulas.

Section~\ref{sec:functional-classification} asks directly which
functions $f(i,j)$ give these structures.  Finite shift spaces
classify bilateral Hankel and Toeplitz kernels.  Functional equations
must distinguish one integer grid, a fixed step at every base point,
and independently variable sampling steps.  The rational and
polynomial cases admit a separated-rank classification; continuous
and analytic counterexamples identify the necessary domain
hypotheses.  Matching the upper and lower formulas on the main
diagonal can narrow the band further, as the explicit
exponential-distance examples show.

The ingredients have established antecedents.  Independent two-layer
evolution belongs to coefficient-deformed condensation and the
octahedron recurrence
\cite{RobbinsRumsey1986,Speyer2007,DiFrancesco2013}.
The rank classes are classical
\cite{EidelmanGohbergOlshevsky2005,EidelmanGohbergHaimovici2014,EidelmanGohberg2008},
as are ordinary Toeplitz determinant recurrences
\cite{Shevelev1991,Shevelev1992,Shevelev1993}.
Exact determinant formulas for rational Toeplitz symbols were obtained
by Day~\cite{Day1975}.
We use the constructive row--column results of
\cite{AlekseyevKhomovskyRecurrences2026} and the symmetry reductions
of \cite{AlekseyevKhomovskySymmetries2026}.
Our results identify the connected-minor support,
the ordinary ancestors that preserve a specified descendant, and
the boundary data entering those recurrences.  The proofs use ordinary
ancestors and determinant transfers.

The appendices give an interleaved hand calculation, the auxiliary
sharpness and transformation arguments, scalar comparisons, and the
explicit ordinary transfer for position-dependent pentadiagonal
matrices.  Inverse condensation also raises questions about minimal
depths, families of parents, and compatibility across matrix sizes.
Here each structural family is treated by
constructing one ordinary parent, verifying its descendant identity,
and deriving the recurrence from its explicit bulk and boundary data.

\section{Condensation and ordinary representatives}\label{sec:generalized}\label{sec:preliminaries}\label{sec:known-two-layer}

We work over a field $\F$.  For polynomial or rational parameter
families, a generic assertion means one valid on a nonempty
Zariski-open subset of the stated parameter space.  For ordered index
sets $I,J$, write $M[I,J]$ for the corresponding submatrix, and use
$[a:b]:=\{a,a+1,\ldots,b\}$ for an inclusive integer interval.
A \emph{connected minor} is the determinant of a square submatrix
on consecutive rows and consecutive columns.  The determinant of the
empty matrix is one.  Throughout,
$\mathsf S_-^{(n)}$ and $\mathsf S_+^{(n)}$ denote the lower and
upper one-step nilpotent shifts of order $n$, respectively.

For a bi-infinite symbol $b=(b_s)_{s\in\mathbb Z}$, put
\[
 T_n(b):=(b_{j-i})_{1\le i,j\le n}.
\]
A matrix is \emph{$(p,q)$-banded} if its entries vanish for
$j-i\notin[-p,q]$; $p$ and $q$ are its lower and upper
semibandwidth bounds.  The least such bounds are the semibandwidths;
the corresponding total bandwidth is $p+q+1$.  A \emph{bulk formula} specifies the entries away
from a stated boundary region.  A \emph{fixed northwest correction}
is an additive matrix supported in a fixed northwest $k\times k$
block, whose entries are independent of the total matrix size.
The correction may contain entries outside the bulk band.

For scalar sequences, $\mathcal E$ denotes the forward shift:
$(\mathcal E u)_n:=u_{n+1}$.  A polynomial
$P(z):=\sum_{j=0}^d p_jz^j$, with coefficients independent of $n$
and $p_d\ne0$, is an
\emph{annihilator} of $(u_n)$ on a stated range when
$P(\mathcal E)u=0$ there, that is,
$\sum_{j=0}^d p_ju_{n+j}=0$.  The recurrence is \emph{eventual}
if it holds for all sufficiently large $n$.  Its minimal order is
the least degree of a nonzero annihilator on the stated range.

For a position-dependent array $X=(X_{ij})$, index its diagonals by their
signed offsets and their entries from the northwest end:
\begin{equation}\label{eq:signed-diagonal-notation}
 x_s[i]:=X_{\,i+\max\{-s,0\},\,i+\max\{s,0\}}.
\end{equation}
Thus $x_{-1}[i]=X_{i+1,i}$, $x_0[i]=X_{ii}$, and
$x_1[i]=X_{i,i+1}$.  Negative offsets denote subdiagonals, zero denotes
the main diagonal, and positive offsets denote superdiagonals.
For two layers $A,B$ we use $a_s[i],b_s[i]$, respectively, each counted
from its own northwest end.  A monomial in these entries is balanced
when the sum of its signed offsets, weighted by their exponents, is zero.
An unindexed $x_s$ denotes a constant value on diagonal $s$ or in its
stated bulk; exceptional entries retain their position index.

Let $A=(a_{ij})$ be an $n\times n$ matrix and let $B=(b_{ij})$ be an
$(n-1)\times(n-1)$ matrix.  We call $A$ the \emph{first layer}
and $B$ the \emph{denominator layer}.  We assume first that every
entry of $B$ is nonzero.

\begin{definition}\label{def:unit}\label{def:block}
The pair $(A,B)$ is a \emph{Dodgson two-layer block}, or simply a
\emph{Dodgson block} or \emph{two-layer block}, if every internal zero of $A$ is
compatible with the corresponding $2\times2$ minor of $B$: whenever
$1<i<n$, $1<j<n$, and $a_{ij}=0$, one has
\begin{equation}\label{eq:compat}
 \det\begin{bmatrix}
 b_{i-1,j-1}&b_{i-1,j}\\
 b_{i,j-1}&b_{i,j}
 \end{bmatrix}=0.
\end{equation}
\end{definition}

The order of the block is the order $n$ of its first layer.
Starting from $A^{(0)}:=A$ and $B^{(0)}:=B$, define recursively
\begin{align}
 a_{ij}^{(k+1)}
 &:=\frac{a_{ij}^{(k)}a_{i+1,j+1}^{(k)}
       -a_{i,j+1}^{(k)}a_{i+1,j}^{(k)}}{b_{ij}^{(k)}},
 \label{eq:condA}\\
 b_{ij}^{(k+1)}&:=a_{i+1,j+1}^{(k)}.
 \label{eq:condB}
\end{align}
Here $A^{(k)}$ has order $n-k$ and $B^{(k)}$ has order $n-k-1$.
Write $\operatorname{Cond}(A,B):=(A^{(1)},B^{(1)})$ and use
$\operatorname{Cond}^k$ for its $k$th iterate.  The all-ones matrix
of order $m$ is denoted by $J_m$, or by $J$ when its order is clear.

\begin{definition}
The scalar remaining after $n-1$ iterations is
\[
 \Ddet(A,B):=a_{11}^{(n-1)}.
\]
The pair $(A^{(k)},B^{(k)})$ is the $k$th \emph{descendant} of $(A,B)$.
\end{definition}

\paragraph{Established two-layer evolution.}
Replacing each numerator $ab-cd$ by $ab+\lambda cd$ gives the
coefficient-deformed condensation of Robbins and Rumsey; the update
above is its $\lambda=-1$ specialization.  The two
initial layers may be independent
\cite{RobbinsRumsey1986,Speyer2007}.  Its octahedron-recurrence
interpretation gives Laurent and perfect-matching formulas.  Di
Francesco's Theorem~20 and Section~6.3 give a selected T-system value
as a network minor with a boundary factor
\cite{DiFrancesco2013}.  Here the structural constructions specify
which connected-minor descendant is preserved and control the ordinary
representative's band and boundary entries.  The procedure itself is
ordinary condensation.

The specialization $\lambda=1$ is not a general permanent formula:
on a $3\times3$ all-ones input it gives $8$, whereas the permanent
is $6$.  We use ordinary representatives and their determinant transfers.

\subsection{Zeros and regularization}

For ordinary input $(H,J)$, the two layers are connected minors of
consecutive orders.  These give canonical polynomial values even when
the displayed quotients have zero denominators.  Computation with an
auxiliary indeterminate $\eps$ means perturbing the underlying entries,
forming all descendants consistently, cancelling, and only then
specializing.  For general two-layer data, we first restrict the
parameters so that every required denominator is nonzero; we call this a
\emph{nonzero-denominator chart}.  A specialization is admitted only
where the resulting rational expressions, simplified on the stated
parameter family, have well-defined finite values.  The local condition
\eqref{eq:compat} alone does not guarantee this regularity.

Structural zeros require a specified continuation.  If both layers
have support $-p\le j-i\le q$, retain the zeros outside this band
and apply \eqref{eq:condA} inside it.  Outside the band the undivided
identity is $0=0$, so the support prescription is part of the rule.
For the ordinary inputs below, vanishing connected minors justify
these zeros at every later order; in particular this follows from
Theorem~\ref{thm:band-core-compression} for rank-compressed inputs.
We call such a rule on the underlying parameter family a
\emph{specified forward continuation}.  Identities of a common
continuation extend to finite specializations.  Equality of two
specialized layers without that common rule need not determine the
terminal determinant.

\begin{remark}\label{rem:zero-quotient-is-not-zero}
Replacing a quotient by zero whenever its denominator vanishes does not
implement this convention.  For example, ordinary condensation of
\[
 H_\eps:=\begin{bmatrix}1&1&0\\1&\eps&1\\0&1&1\end{bmatrix}
\]
ends with
\[
 \frac{(\eps-1)^2-1}{\eps}=\eps-2.
\]
At $\eps=0$ the determinant is $-2$, although the unsimplified quotient
is $0/0$.  For ordinary input one may always define a zero-denominator
specialization by its connected minors, or perturb the original matrix
and cancel before specializing.  Independently assigning values to
undefined intermediate quotients can destroy the minor identities.
\end{remark}

\needspace{6\baselineskip}
Using the same continuation is necessary even for ordinary parents.
Let $J_2$ denote the all-ones matrix of order two.  For $u\in\mathbb Q$,
consider
\begin{equation}\label{eq:singular-common-block}
\begin{gathered}
 H_\eps(u):=\begin{pmatrix}u&1&0\\1&\eps&1\\0&1&1\end{pmatrix},\\
 \operatorname{Cond}(H_\eps(u),J_2)
 =\left(\begin{pmatrix}u\eps-1&1\\1&\eps-1\end{pmatrix},[\eps]\right).
\end{gathered}
\end{equation}
At $\eps=0$ the displayed block is independent of $u$, but
$\det H_0(u)=-u-1$.  Thus two equal singular blocks may come from
ordinary parents with different terminal determinants.  The identities
transported below use one specified continuation on a common rational
family before specialization.

For hand calculations, Appendix~\ref{app:interleaved-layers} displays
both layers in one interleaved array, with each denominator placed in
bold at the center of its four numerator entries.  A complete
$5\times5$ example follows one perturbation through all condensation
steps and illustrates cancellation before specialization.

\subsection{Desnanot--Jacobi, transpose, and homogeneity}

For a square matrix $M$, let $M_i^j$ denote deletion of row $i$ and column
$j$, with analogous notation for two deleted rows and columns.

\begin{lemma}\label{lem:DJD}
If $(A,B)$ is a Dodgson block of order $n\ge3$, then
\begin{equation}\label{eq:DJD}
\begin{aligned}
&\Ddet(A,B)\,
 \Ddet\!\left(A_{1,n}^{1,n},B_{1,n-1}^{1,n-1}\right)\\
&\qquad=\Ddet\!\left(A_1^1,B_1^1\right)
          \Ddet\!\left(A_n^n,B_{n-1}^{n-1}\right)
 -\Ddet\!\left(A_1^n,B_1^{n-1}\right)
          \Ddet\!\left(A_n^1,B_{n-1}^{1}\right).
\end{aligned}
\end{equation}
\end{lemma}

\begin{proof}
After $k\ge1$ condensation steps, the entries of $A^{(k)}$ are the
scalar condensation values of consecutive two-layer subblocks of order
$k+1$, while the entries of $B^{(k)}$ are those of the corresponding
subblocks of order $k$.  At the final nontrivial step, equation~\eqref{eq:condA} is exactly
\eqref{eq:DJD} for the four corner descendants.  The regularized
interpretation handles vanishing intermediate denominators.
\end{proof}

\begin{lemma}\label{lem:transpose-scale}
For an order-$n$ Dodgson block and $c\ne0$,
\[
 \Ddet(A^T,B^T)=\Ddet(A,B),
 \qquad
 \Ddet(cA,B)=c^n\Ddet(A,B).
\]
Both identities extend to regularized specializations.
\end{lemma}

\begin{proof}
Transpose invariance follows at every local update.  For homogeneity, if the
$k$th $A$-layer scales as $c^{k+1}$ and the $k$th $B$-layer as $c^k$, then
\eqref{eq:condA} makes the next $A$-layer scale as $c^{k+2}$ and
\eqref{eq:condB} makes the next $B$-layer scale as $c^{k+1}$.  The unique
entry of the final $A$-layer therefore scales as $c^n$.
\end{proof}

\begin{theorem}\label{thm:classical}
Let $J_{n-1}$ be the all-ones matrix.  Then
\[
 \Ddet(A,J_{n-1})=\det A.
\]
\end{theorem}

\begin{proof}
The two starting layers represent ordinary connected minors of orders one and
zero.  Induction with the classical Desnanot--Jacobi identity shows that the
subsequent layers contain ordinary connected minors of consecutive orders.
The final entry is therefore $\det A$.
\end{proof}

\subsection{Descendants and recurrence transport}

On a nonzero-denominator chart, write $\operatorname{Cond}(U)$ for one
condensation step applied to a Dodgson block $U=(A,B)$; regularized
specializations may be included using the convention above.  If $U$
has order $N$ and $2\le r\le N$, and the required iterates are defined,
define its order-$r$ descendant by
\begin{equation}\label{eq:rdescendant}
 D_r(U):=\operatorname{Cond}^{N-r}(U).
\end{equation}
Equality of descendants always means equality of both layers.

\begin{definition}
\label{def:condequiv}
For Dodgson blocks $U,V$ of respective orders $N,M\ge r$, write
\begin{equation}\label{eq:condequiv}
 U\sim_r V
 \qquad\Longleftrightarrow\qquad
 D_r(U)=D_r(V).
\end{equation}
\end{definition}

\begin{proposition}
\label{prop:condequiv}
For every fixed $r\ge2$, the relation $\sim_r$ is an equivalence relation on
Dodgson blocks of order at least $r$ for which the order-$r$ descendant is
defined.  If both inputs respect one specified forward continuation of
their common order-$r$ descendant, then
\begin{equation}\label{eq:condequivmonotone}
 U\sim_rV\quad\Longrightarrow\quad U\sim_sV
 \qquad(2\le s\le r),
\end{equation}
and
\begin{equation}\label{eq:condequivddet}
 U\sim_rV\quad\Longrightarrow\quad \Ddet(U)=\Ddet(V).
\end{equation}
The continuation may be given on a nonzero-denominator chart or by a
specified rational rule with structural zeros justified by minor
identities.  These implications extend to finite specializations of
that common continuation.  Equality of the specialized order-$r$ arrays
alone is insufficient, as \eqref{eq:singular-common-block} shows.
\end{proposition}

\begin{proof}
Reflexivity and symmetry are immediate.  If $U\sim_rV$ and $V\sim_rW$,
then the three order-$r$ descendants are equal, proving transitivity.
Under the stated common-continuation hypothesis, applying further
condensation gives~\eqref{eq:condequivmonotone}; applying the
remaining $r-1$ steps to the common order-$r$ descendant gives
\eqref{eq:condequivddet}.
\end{proof}

An \emph{ordinary ancestor} of an order-$r$ block $C$ at depth $s$
is a matrix $M$ of order $r+s$ for which
$\operatorname{Cond}^s(M,J)=C$, with the stated continuation.
We also call the reconstructed matrix an \emph{ordinary parent};
the specified depth records how many reverse steps are involved.
For the recurrence applications, a later common descendant is enough,
with continuation understood as in the preceding proposition.
If $C:=\operatorname{Cond}^{t}(U)$ and an ordinary matrix $M$ satisfies
$\operatorname{Cond}^{s}(M,J)=C$, then
\begin{equation}\label{eq:late-descendant-suffices}
 \Ddet(U)=\Ddet(C)=\det M.
\end{equation}
Thus $M$ need not reconstruct the entire earlier block $U$.
For a family of targets, the entries and sizes of the ordinary
representatives must also be controlled so that their determinant
recurrence can be transported.  The constructions below give those
entries directly from the structural parameters; the common-descendant
identity certifies the result.

\begin{corollary}\label{cor:recurrencetransport}
Let $a,b\ge0$ be fixed integers.  If families of orders $n+a,n+b$
satisfy
\begin{equation}\label{eq:condshiftequiv}
 \operatorname{Cond}^{a}(U_{n+a})
 =\operatorname{Cond}^{b}(V_{n+b})
\end{equation}
for all sufficiently large $n$, with one specified forward continuation,
then
\begin{equation}\label{eq:ddetshiftidentity}
 \Ddet(U_{n+a})=\Ddet(V_{n+b}).
\end{equation}
Every finite scalar recurrence transfers between these sequences with
the corresponding index shift, whether its coefficients are constant
or position dependent.
\end{corollary}

\begin{proof}
Condense the common order-$n$ descendant to a scalar, obtaining
\eqref{eq:ddetshiftidentity}, and substitute this equality in the
recurrence.
\end{proof}

\subsection{Boundary anchoring of ordinary ancestors}

A determinant-preserving construction gives an ordinary ancestor of a
specified descendant only when it preserves every connected minor in
that descendant.  The following boundary normalization ensures this
stronger property.  We call it \emph{anchoring} because the first
$h$ coordinates of the row and column transformations remain fixed.

\begin{proposition}\label{prop:anchored-minor-preservation}
Let $H$ have order $n$, fix $0\le h<n$, and let
\[
 L:=\diag(I_h,L_*),\qquad R:=\diag(I_h,R_*),
\]
where $L_*$ is unit lower triangular and $R_*$ is unit upper triangular.
Put $M:=LHR$.  If $I_0=[i:i+t-1]$ and $J_0=[j:j+t-1]$ are valid intervals
with $i,j\le h+1$, then
\begin{equation}\label{eq:anchored-minor-preservation}
 \det M[I_0,J_0]=\det H[I_0,J_0].
\end{equation}
In particular,
\begin{equation}\label{eq:anchored-full-descendant}
 \operatorname{Cond}^{n-h-1}(M,J)
 =\operatorname{Cond}^{n-h-1}(H,J).
\end{equation}
For $h=0$ this preserves the final scalar; for $h=1$ it preserves both
layers of the order-two descendant.
\end{proposition}

\begin{proof}
Every row of $L$ indexed by $I_0$ uses only columns in $I_0$.  Indeed, rows
at most $h$ are unchanged, while all other nonzero columns are between
$h+1$ and the row index.  Because $i\le h+1$, none lies before $i$.
The corresponding assertion for columns of $R$ and the interval $J_0$
follows from upper triangularity.  Hence
\[
 M[I_0,J_0]=L[I_0,I_0]\,H[I_0,J_0]\,R[J_0,J_0].
\]
The two outside factors have determinant one, proving
\eqref{eq:anchored-minor-preservation}.  The first layer in
\eqref{eq:anchored-full-descendant} consists of connected minors of
order $n-h$ whose starting indices lie in $[1:h+1]$.  Its second layer
consists of connected minors of order $n-h-1$ whose starting indices
lie in $[2:h+1]$.  All of them are therefore preserved.  Ordinary
minor identities also cover zero intermediate denominators.
\end{proof}

\paragraph{Other minor choices.}\label{sec:pivot-connected-comparison}
Pivoted reductions use different index sets.  For
$H=\begin{pmatrix}C&E\\F&G\end{pmatrix}$ with an invertible
order-$k$ core $C$, put $\delta:=\det C$ and $m:=n-k$.
The bordered-minor matrix satisfies the Sylvester--Schur identities
\cite{BrualdiSchneider1983}
\begin{equation}\label{eq:pivot-sylvester-normalization}
 Z_{ij}:=\det\begin{pmatrix}C&E_{:j}\\F_{i:}&G_{ij}\end{pmatrix},
 \qquad Z=\delta(G-FC^{-1}E),\qquad
 \det Z=\delta^{m-1}\det H.
\end{equation}
Bareiss elimination uses nested pivot cores and exact division by
preceding pivot minors \cite{Bareiss1968}.  Dodgson condensation uses
consecutive intervals instead.  Reordering an original matrix and
recomputing its minors is legitimate; reordering one intermediate
layer alone need not preserve their overlaps or the terminal value.
The bandwidth claims below concern the given row and column order.

\section{From a dense kernel to a tridiagonal determinant}
\label{sec:rankone-model}

We begin with a dense matrix whose off-diagonal entries are products
of a row factor and a column factor, while the main diagonal is
arbitrary.  We compute both banded layers, construct an ordinary
tridiagonal parent of their order-two descendant, and derive the
determinant recurrence by expanding the final row of that parent.
For each $n$, let
\[
 H_n:=(H_{ij})_{1\le i,j\le n},\qquad
 H_{ij}:=\begin{cases}f_i g_j,&i\ne j,\\ h_i,&i=j,\end{cases}
 \qquad
 q_i:=f_i g_i,\qquad \delta_i:=h_i-q_i,
\]
and put $D_n:=\det H_n$, with $D_0:=1$.
We first work with generic parameters and $q_i\ne0$.
The diagonal perturbation matters: without it, $H_n$ has rank at most one
and $D_n=0$ for every $n\ge2$.

\subsection{An explicit ordinary parent and its recurrence}
\label{sec:rankone-corner-route}

\paragraph{Condensation exposes the band.}
Every connected minor of order at least two whose starting row and
column differ by at least two has two proportional rows in one of the separated triangles.
Thus one condensation step gives a tridiagonal first layer, and two
steps give two tridiagonal layers.  Fix $N\ge5$ and put $m:=N-2$.
Write $U_m:=\operatorname{Cond}^2(H_N,J)=(A_m,B_{m-1})$.
The main diagonals of the two layers are
\begin{align}\label{eq:rankone-two-layer-main}
 a_0[i]&:=\delta_i\delta_{i+1}\delta_{i+2}
       +\sum_{v=i}^{i+2}q_v
          \prod_{\substack{w=i\\w\ne v}}^{i+2}\delta_w,
       &&1\le i\le m,\notag\\
 b_0[i]&:=h_{i+1}h_{i+2}-q_{i+1}q_{i+2},
       &&1\le i\le m-1.
\end{align}
Their superdiagonals and subdiagonals are
\begin{align}\label{eq:rankone-two-layer-values}
 a_1[i]&:=f_i g_{i+3}\delta_{i+1}\delta_{i+2},&
 a_{-1}[i]&:=f_{i+3}g_i\delta_{i+1}\delta_{i+2},\notag\\
 b_1[i]&:=-f_{i+1}g_{i+3}\delta_{i+2},&
 b_{-1}[i]&:=-f_{i+3}g_{i+1}\delta_{i+2}.
\end{align}
Here the entries of $a_{\pm1}$ have $1\le i\le m-1$, and those of
$b_{\pm1}$ have $1\le i\le m-2$.  Expanding the connected minors of
orders three and two proves these formulas.  The indices $i+1,i+2$
in $b_0[i]$ record the central placement of the denominator layer.

\paragraph{Constructing the ordinary parent.}
For the determinant recurrence it is enough to retain a suitable later
descendant.  The entries of the following ordinary matrix are given
directly by the original parameters.  Let $M_n$ be tridiagonal, with entries
\begin{equation}\label{eq:rankone-corner-data}
 \begin{aligned}
 m_0[1]&:=h_1,\qquad m_0[2]:=h_2,\\
 m_0[i]&:=\delta_i+\delta_{i-1}\frac{q_i}{q_{i-1}}
                       &&(i\ge3),\\
 m_1[1]&:=f_1g_2,\qquad m_{-1}[1]:=f_2g_1,\\
 m_1[i]&:=-\delta_i\frac{g_{i+1}}{g_i},\qquad
 m_{-1}[i]:=-\delta_i\frac{f_{i+1}}{f_i}
                       &&(i\ge2).
 \end{aligned}
\end{equation}
For each size, use only entries with valid indices; in particular
$M_1=[h_1]$.

\begin{proposition}\label{prop:rankone-corner-vapor}
Assume $q_i\ne0$.  For $n\ge2$,
\begin{equation}\label{eq:rankone-corner-descendant}
 \operatorname{Cond}^{n-2}(M_n,J)
 =\operatorname{Cond}^{n-2}(H_n,J).
\end{equation}
The matrices $M_n$ are leading principal sections of one fixed
tridiagonal array, and $\det M_n=D_n$.
Extending the bulk formulas in \eqref{eq:rankone-corner-data} to
$i=2$ on the main diagonal and $i=1$ on the adjacent diagonals changes
only the northwest $2\times2$ block, independently of the total size.
\end{proposition}

\begin{proof}
Let $L_n$ and $R_n$ be unit lower and upper bidiagonal matrices,
respectively, with
\[
 l_{-1}[i]:=-\frac{f_{i+1}}{f_i},\qquad
 r_1[i]:=-\frac{g_{i+1}}{g_i}
 \qquad(2\le i\le n-1),
\]
and all other off-diagonal entries zero.  Direct multiplication gives
$M_n=L_nH_nR_n$: entries outside the tridiagonal band cancel,
the initial $2\times2$ block is unchanged, and the remaining entries
are exactly \eqref{eq:rankone-corner-data}.
Both factors have an isolated first identity coordinate.
Proposition~\ref{prop:anchored-minor-preservation} with $h=1$
therefore proves \eqref{eq:rankone-corner-descendant}.
Further condensation gives $\det M_n=\det H_n$.
All entries depend only on their indices and the original parameter
sequences, proving the leading-section assertion and the corner support.
These are connected-minor identities, so zero intermediate condensation
denominators require no additional hypothesis on the $\delta_i$.
\end{proof}

For the earlier block $U_m$ with $m=N-2$, the equality reads
\begin{equation}\label{eq:rankone-corner-route}
 \operatorname{Cond}^{N-4}(U_m)
 =\operatorname{Cond}^{N-2}(M_N,J).
\end{equation}
Thus $M_N$ is an ordinary ancestor of the complete order-two descendant.
Its entries are already explicit in the original parameters; this
identity certifies the parent without requiring the intermediate
condensation layers to be evaluated.  The entries away from the corner
are generally position dependent,
so a fixed boundary correction alone does not make the recurrence
coefficients constant.

\paragraph{The ordinary last-row recurrence.}
For $n\ge3$, expanding the final row of $M_n$ gives
\[
 \begin{aligned}
 D_n&=m_0[n]D_{n-1}-m_{-1}[n-1]m_1[n-1]D_{n-2},\\
 m_{-1}[n-1]m_1[n-1]
 &=\delta_{n-1}^{\,2}\frac{q_n}{q_{n-1}}.
 \end{aligned}
\]
Consequently
\begin{equation}\label{eq:rankone-vapor-recurrence}
 D_n=
 \left(\delta_n+\delta_{n-1}\frac{q_n}{q_{n-1}}\right)D_{n-1}
 -\delta_{n-1}^{\,2}\frac{q_n}{q_{n-1}}D_{n-2},
 \qquad D_0:=1,\quad D_1:=h_1.
\end{equation}
The same formula holds for $n=2$, because
\[
 \left(\delta_2+\delta_1\frac{q_2}{q_1}\right)h_1
 -\delta_1^2\frac{q_2}{q_1}
 =h_1h_2-q_1q_2=D_2.
\]
The constructed ordinary parent therefore proves the recurrence for every
$n\ge2$ wherever the $q$-denominators are nonzero, directly from the stated structural parameters.

For the specialization $h_i=q_i-x$, put $\Phi_n:=\det H_n$.  Then
$H_n=(f_i g_j)_{i,j=1}^n-xI_n$ and $\delta_i=-x$, so the recurrence is
\begin{equation}\label{eq:rankone-characteristic-vapor}
 \Phi_n=-x\left(1+\frac{q_n}{q_{n-1}}\right)\Phi_{n-1}
        -x^2\frac{q_n}{q_{n-1}}\Phi_{n-2},
 \qquad \Phi_0:=1,\quad\Phi_1:=q_1-x.
\end{equation}

\begin{example}\label{ex:separable-power-vapor}
Work in characteristic zero.  Let $h_i=x$, $f_i=w i^p$, and
$g_i=i^q$, where $p,q\ge0$ are integers, and put $s:=p+q$.  Then
\[
 H_{ii}=x,\qquad H_{ij}=w i^p j^q\quad(i\ne j),
\]
and \eqref{eq:rankone-vapor-recurrence} gives
\begin{align}\label{eq:separable-power-vapor-recurrence}
 D_n={}&
 \frac{\bigl((n-1)^s+n^s\bigr)x-2w\bigl(n(n-1)\bigr)^s}
      {(n-1)^s}D_{n-1}\notag\\
 &-\frac{n^s\bigl(x-w(n-1)^s\bigr)^2}{(n-1)^s}D_{n-2},
 \qquad D_0=1,\quad D_1=x.
\end{align}
The coefficients in this recurrence have no $w$ in their denominators,
so polynomial specialization also gives the case $w=0$.
For $w=1$ and $p=q=1$, this is
\[
 D_n=\frac{(2n^2-2n+1)x-2n^2(n-1)^2}{(n-1)^2}D_{n-1}
      -\frac{n^2\bigl(x-(n-1)^2\bigr)^2}{(n-1)^2}D_{n-2}.
\]
Both formulas arise from the explicitly banded descendant and its ordinary
parent, even though the original matrices are dense.
\end{example}

The determinant-lemma comparison is recorded in
Appendix~\ref{app:rankone-comparison}.  It verifies the same scalar
answer.  The corner construction above identifies an exact ordinary
ancestor of the late order-two block and gives its boundary entries explicitly.

\section{Structural compression and ordinary ancestors}
\label{sec:structural-route}

The rank-one calculation extends in three stages.  Ranks outside a
prescribed diagonal band bound the depth and width of condensation.
Local recurrences, applied while fixing the first coordinate, then specify
an ordinary ancestor of the complete order-two descendant.  When its band
entries become constant away from the boundary, we write it as a Toeplitz
bulk plus a fixed northwest correction.  The ordinary row--column
construction gives a determinant recurrence from these entries.
Sharpness, the regular converse,
and rank bounds for products, inverses, and structured similarities
are proved in Appendix~\ref{app:some-aspects}.

\subsection{Rank outside a band and the compression depth}

Let $H=(h_{ij})_{1\le i,j\le N}$ and define its connected minors by
\begin{equation}\label{eq:application-connected-minors}
 C_t(i,j):=\det H[i:i+t-1,j:j+t-1],\qquad C_0(i,j):=1.
\end{equation}
Intervals in a submatrix are inclusive.  The two layers after $k$ ordinary
condensation steps are
\begin{equation}\label{eq:application-two-layers}
 A^{(k)}_{ij}=C_{k+1}(i,j),\qquad
 B^{(k)}_{ij}=C_k(i+1,j+1).
\end{equation}
At zero intermediate denominators, these layers are defined by the
ordinary connected minors in \eqref{eq:application-connected-minors}.

The upper and lower quasiseparable ranks of $H$ are the largest ranks of
$H[1:c,c+1:N]$ and $H[c+1:N,1:c]$, respectively, over all cuts $c$.
This standard rank notion includes independently separable upper and lower
triangles and diagonal perturbations of low-rank matrices; see
\cite{EidelmanGohbergOlshevsky2005,EidelmanGohbergHaimovici2014}.

Low-rank structure need only hold outside a fixed diagonal band.  This
out-of-band rank class is classical \cite{EidelmanGohberg2008}; the
following statement gives sharp bounds on connected-minor support.  For
nonnegative integers $b_+,b_-$, define the corresponding upper and lower
rectangle ranks by
\[
 \mathfrak r_+(b_+;H):=
 \max_c\operatorname{rank}H[1:c,c+b_++1:N],\qquad
 \mathfrak r_-(b_-;H):=
 \max_c\operatorname{rank}H[c+b_-+1:N,1:c].
\]
Empty rectangles have rank zero.  In an infinite array the maxima are
replaced by suprema of the ranks of finite such rectangles.  The case
$b_+=b_-=0$ is the usual upper and lower quasiseparable rank.

\begin{theorem}\label{thm:band-core-compression}\label{thm:structural-band-compression}
Suppose $\mathfrak r_\pm(b_\pm;H)\le r_\pm$, and put
\[
 w_\pm:=b_\pm+r_\pm,\qquad R:=\max(r_-,r_+).
\]
Then the connected minors satisfy
\begin{equation}\label{eq:structural-connected-support}
\begin{aligned}
 C_t(i,j)&=0&&\text{if }t>r_+\text{ and }j-i>w_+,\\
 C_t(i,j)&=0&&\text{if }t>r_-\text{ and }i-j>w_-.
\end{aligned}
\end{equation}
After $R$ ordinary condensation steps the first layer is
$(w_-,w_+)$-banded.  After $R+1$ steps both layers are
$(w_-,w_+)$-banded.  These assertions concern every valid entry of
the indicated finite sections; at zero intermediate denominators they
are interpreted through ordinary connected minors.
\end{theorem}

\begin{proof}
Put $s:=j-i>w_+$ and $q:=\min(t,s-b_+)$.  The first $q$ rows
of the block $H[i:i+t-1,j:j+t-1]$ lie in the rectangle with cut
$c=j-b_+-1$, whose upper rank is at most $r_+$.  The remaining
$t-q$ rows add at most $t-q$ to the rank.  Hence
\[
 \operatorname{rank}H[i:i+t-1,j:j+t-1]
 \le r_++\max(0,t-s+b_+)<t.
\]
Transposition proves the lower assertion.  The layers at depth $k$
are $C_{k+1}(i,j)$ and $C_k(i+1,j+1)$, which gives the two depths.
\end{proof}

In particular, the theorem applies to
\begin{equation}\label{eq:band-core-separated-data}
 H_{ij}:=
 \begin{cases}
  F_i^+(G_j^+)^T,&j-i>b_+,\\
  d_{ij},&-b_-\le j-i\le b_+,\\
  F_i^-(G_j^-)^T,&i-j>b_-,
 \end{cases}
\end{equation}
where the upper and lower generator vectors have lengths $r_+$ and
$r_-$, respectively.  The entries $d_{ij}$ in the diagonal band
are arbitrary.  The two separated representations are independent,
and every generator component may be an arbitrary function of its
own index.  No Toeplitz, Hankel, polynomial, or constant-coefficient
assumption is required for this support statement.

For independent generator and band entries these depth and width
bounds are generically sharp
(Proposition~\ref{prop:band-core-generic-sharpness}).
A converse from vanishing connected minors requires a nonvanishing
hypothesis on the preceding minors
(Corollary~\ref{cor:regular-band-core-converse}).
The proofs and a counterexample without that hypothesis are in
Appendix~\ref{app:some-aspects}.

Band formation alone does not imply a constant-coefficient recurrence.
Already a diagonal array can have arbitrary nonzero principal
determinants: prescribe $D_0:=1$ and use diagonal entry $D_n/D_{n-1}$.
The following dense example shows that polynomial diagonal sequences
in the condensed layers do not supply the missing condition.

\begin{example}\label{ex:banded-not-constant-recurrent}
Work in characteristic zero and put
\[
 H_n:=\bigl(i(1+\mathbf 1_{i=j})\bigr)_{1\le i,j\le n},
 \qquad D_n:=\det H_n,\qquad D_0:=1.
\]
This is a dense diagonal-plus-rank-one family.  If $\boldsymbol 1_n$
denotes the all-ones column, then
\[
 H_n=\operatorname{diag}(1,\ldots,n)
       (I_n+\boldsymbol 1_n\boldsymbol 1_n^{\mathsf T}),
 \qquad D_n=(n+1)!.
\]
Consequently $D_n=(n+1)D_{n-1}$ for $n\ge1$.

For every contained connected minor of order $k\ge2$,
\[
 \det H_n[i:i+k-1,j:j+k-1]
 =\left(\prod_{v=i}^{i+k-1}v\right)
 \begin{cases}
  k+1,&j=i,\\
  (-1)^{k-1},&|j-i|=1,\\
  0,&|j-i|\ge2.
 \end{cases}
\]
Indeed, remove the row factors.  At shift zero the remaining matrix is
$I_k+\boldsymbol 1_k\boldsymbol 1_k^{\mathsf T}$; at shifts $1$ and $-1$
it is the all-ones matrix plus a single subdiagonal or superdiagonal of
ones.  Their determinants are $k+1$ and $(-1)^{k-1}$, respectively.
At larger shifts, the restricted identity has rank at most $k-2$, and
adding the all-ones matrix raises rank by at most one.

Thus the first condensation makes the top layer tridiagonal, and the
second makes both layers tridiagonal.  For $n\ge4$, write these latter layers as
$A_{n-2}$ and $B_{n-3}$, with the denominator centrally placed.  Their
signed diagonals, at valid indices, are
\begin{align*}
 a_{-1}[i]&=(i+1)(i+2)(i+3),&
 b_{-1}[i]&=-(i+2)(i+3),\\
 a_0[i]&=4i(i+1)(i+2),&
 b_0[i]&=3(i+1)(i+2),\\
 a_1[i]&=i(i+1)(i+2),&
 b_1[i]&=-(i+1)(i+2).
\end{align*}

Nevertheless, $(D_n)$ satisfies no nontrivial eventual linear recurrence
with constant coefficients.  Such a recurrence would give
$\sum_{j=0}^{r}c_j(n+j+1)!=0$ for all sufficiently large $n$.
Dividing by $(n+1)!$ gives the polynomial identity
\[
 \sum_{j=0}^{r}c_j\prod_{t=2}^{j+1}(n+t)=0,
\]
where the product for $j=0$ is $1$.  The summands have distinct degrees
$j$, so comparison of the highest degree, followed by induction, forces
every $c_j$ to vanish.

Even a constant main diagonal does not suffice.  Let
$(K_n)_{ii}:=1$ and $(K_n)_{ij}:=i$ for $i\ne j$.  Its off-diagonal
rectangles again have rank one, so both connected-minor layers are
tridiagonal after two condensation steps.  Subtracting $i$ times the
first row from row $i$, for $2\le i\le n$, gives
\[
 \det K_n=\prod_{i=2}^{n}(1-i)=(-1)^{n-1}(n-1)!,\qquad n\ge1.
\]
The same factorial argument excludes an eventual constant-coefficient
recurrence for this sequence.
\end{example}

The two depth statements distinguish support of the first layer from
support of both layers.  For rank one, one step makes the first layer
tridiagonal; two steps make both layers tridiagonal.  A denominator
with nonzero tails is permitted in the former statement.

For diagonal-plus-low-rank matrices, every surviving entry can be given
without solving a nonlinear system.  Let
\[
 H=\diag(\delta_1,\ldots,\delta_N)+UV^T,
 \qquad U,V\in\F^{N\times r}.
\]
For $t\ge r+1$, $s\ge0$, set
$I:=\{i,\ldots,i+t-1\}$, $J:=\{i+s,\ldots,i+s+t-1\}$, and
$O:=I\cap J$; all sets in minors have increasing order.

\begin{proposition}\label{prop:explicit-lowrank-layers}
With the empty determinant equal to $1$,
\begin{equation}\label{eq:explicit-lowrank-layers}
\begin{split}
 C_t(i,i+s)={}&
 \sum_{\substack{Z\subseteq O\\|Z|\ge t-r}}
 (-1)^{s|Z|}\left(\prod_{z\in Z}\delta_z\right)\\
 &\quad\times
 \sum_{\substack{S\subseteq\{1,\ldots,r\}\\|S|=t-|Z|}}
 \det U[I\setminus Z,S]\det V[J\setminus Z,S].
\end{split}
\end{equation}
In particular, on the outer upper diagonal,
\begin{equation}\label{eq:explicit-lowrank-outer}
\begin{split}
 C_t(i,i+r)={}&(-1)^{r(t-r)}
 \left(\prod_{z=i+r}^{i+t-1}\delta_z\right)\\
 &\quad\times\det U[i:i+r-1,:]\det V[i+t:i+t+r-1,:].
\end{split}
\end{equation}
The formulas with $t=r+1$ and $t=r+2$, together with their transposes,
specify both layers of the first fully banded block.
\end{proposition}

\begin{proof}
Expand the determinant by the main diagonal entries selected from
$\diag(\delta)$.  If their index set is $Z$, deleting the corresponding
rows and columns has sign $(-1)^{s|Z|}$.  The remaining minor of $UV^T$
vanishes when its size exceeds $r$; otherwise Cauchy--Binet gives the inner
sum in \eqref{eq:explicit-lowrank-layers}.  For $s=r$, the overlap has
$t-r$ elements, so only $Z=O$ contributes, proving
\eqref{eq:explicit-lowrank-outer}.
\end{proof}

Thus the characteristic-polynomial problem for
$h_{ij}=\sum_{a=1}^r f_a(i)g_a(j)$ has an explicit $(r,r)$-banded
two-layer descendant after $r+1$ steps: take $U_{ia}=f_a(i)$, $V_{ja}=g_a(j)$, and
$\delta_i=-x$ in $H-xI$.  The next construction turns these explicit descendants into ordinary
banded determinants.

\subsection{Ordinary ancestors for separated tails and a band core}

Consider leading sections $H_n$ of the fixed array
\eqref{eq:band-core-separated-data}, with
$F_i^+,G_i^+\in\F^{1\times r_+}$ and
$F_i^-,G_i^-\in\F^{1\times r_-}$.  Empty vectors are allowed when
one of $r_-,r_+$ is zero.

For $i\ge r_-+2$, suppose that the matrix with rows
$F_{i-r_-}^-,\ldots,F_{i-1}^-$ is nonsingular, and define
$\alpha_{i,0}:=1$, $\alpha_{i,1},\ldots,\alpha_{i,r_-}$ by
\begin{equation}\label{eq:local-row-annihilator}
 \sum_{a=0}^{r_-}\alpha_{i,a}F_{i-a}^-=0.
\end{equation}
Likewise, for $j\ge r_++2$, assume nonsingularity of the matrix with
rows $G_{j-r_+}^+,\ldots,G_{j-1}^+$ and define
$\beta_{j,0}:=1$, $\beta_{j,1},\ldots,\beta_{j,r_+}$ by
\begin{equation}\label{eq:local-column-annihilator}
 \sum_{b=0}^{r_+}\beta_{j,b}G_{j-b}^+=0.
\end{equation}
The empty systems when $r_-=0$ or $r_+=0$ require no assumption.
Cramer's rule expresses these coefficients as rational functions of the
indicated consecutive generator entries.

For $i\le r_-+1$ set $\alpha_{i,0}:=1$ and all coefficients with
$a>0$ equal to zero.  Make the analogous choice for
$\beta_{j,b}$ when $j\le r_++1$.  Define
\begin{equation}\label{eq:local-anchored-factors}
 (L_n)_{i,i-a}:=\alpha_{i,a},\qquad
 (R_n)_{j-b,j}:=\beta_{j,b}
\end{equation}
at valid indices, and set all other entries to zero.  Put
\begin{equation}\label{eq:local-vapor-parent}
 M_n:=L_nH_nR_n.
\end{equation}
The identity choices in the initial blocks are convenient, not essential:
one may instead use any fixed unit lower triangular block on indices
$2,\ldots,r_-+1$ and any fixed unit upper triangular block on indices
$2,\ldots,r_++1$.  The first row and column remain anchored.

\begin{theorem}\label{thm:local-anchored-vaporisation}\label{thm:band-core-anchored-vaporisation}
Put $w_-:=b_-+r_-$ and $w_+:=b_++r_+$.  Under the stated nonsingularity
assumptions, the matrices $M_n$ are leading sections of one ordinary
$(w_-,w_+)$-banded array, and for $n\ge2$,
\begin{equation}\label{eq:local-vapor-descendant}
 \operatorname{Cond}^{n-2}(M_n,J)
 =\operatorname{Cond}^{n-2}(H_n,J).
\end{equation}
Every entry of $M_n$ is given explicitly by the finite sum
\begin{equation}\label{eq:local-vapor-entry}
 (M_n)_{ij}
 =\sum_{a=0}^{\min(r_-,i-1)}
   \sum_{b=0}^{\min(r_+,j-1)}
      \alpha_{i,a}\beta_{j,b}H_{i-a,j-b}.
\end{equation}
In particular $\det M_n=\det H_n$.

The ordinary position-dependent row--column construction gives a state
realization of dimension
\begin{equation}\label{eq:local-vapor-state-dimension}
 d:=\binom{w_-+w_+}{w_-}.
\end{equation}
Its transition coefficients are rational functions of a bounded number of
consecutive generator and diagonal-band entries.
\end{theorem}

\begin{proof}
The factors in \eqref{eq:local-anchored-factors} are unitriangular.
Their first coordinate is an isolated identity coordinate.  The allowed
initial blocks use only indices starting at two, and every subsequent
row or column operation also uses indices at least two.
Proposition~\ref{prop:anchored-minor-preservation}
with $h=1$ proves \eqref{eq:local-vapor-descendant}.

Suppose $j-i>b_++r_+$.  Then $j\ge r_++2$, so the complete column
recurrence is applied.  For every participating row index $i-a$ and every
$0\le b\le r_+$, one has $j-b>i+b_+\ge i-a+b_+$.  Consequently
\[
 H_{i-a,j-b}=F_{i-a}^+(G_{j-b}^+)^T,
\]
and \eqref{eq:local-column-annihilator} makes the inner sum in
\eqref{eq:local-vapor-entry} zero.  If $i-j>b_-+r_-$, every participating row
lies below the core relative to every participating column index; summing first in
$a$ and applying \eqref{eq:local-row-annihilator} gives zero.
Thus $M_n$ is $(b_-+r_-,b_++r_+)$-banded.

All formulas use only preceding indices, so the entry formula is
independent of the chosen section size.  Each
nonzero band entry depends on a bounded range of indices, because $j-i$ is
bounded and the two local linear systems have fixed sizes.  Applying
the ordinary variable-band row--column transfer recalled in
Section~\ref{sec:background} gives \eqref{eq:local-vapor-state-dimension}
and the stated locality of its coefficients.
\end{proof}

As an identity of descendants, this construction passes through
order two.  If $R:=\max(r_-,r_+)$ and $n\ge R+3$, then $R+1$ steps
give the banded two-layer block, a further $n-R-3$ steps give its
order-two descendant, and $M_n$ is an ordinary ancestor of that
descendant at depth $n-2$.
These step counts locate the common block; they are not instructions
to compute all intervening minors before constructing the parent.
Formula~\eqref{eq:local-vapor-entry} gives every parent entry directly
from a bounded number of consecutive entries of the original data.

For $H=\diag(\delta)+UV^T$ with $U,V$ having $r$ columns, set
$F_i^+=F_i^-=U[i,:]$, $G_i^+=G_i^-=V[i,:]$, and
$b_-=b_+=0$, and $d_{ii}=\delta_i+U[i,:]V[i,:]^T$.  The theorem gives an ordinary
$(r,r)$-banded ancestor and a state vector of dimension $\binom{2r}{r}$.
Indeed, for the identity initial-block choice,
\[
 M_n=L_n\diag(\delta)R_n+(L_nU)(V^TR_n),
\]
and the second summand is supported in the first $(r+1)\times(r+1)$
corner.  This is also a direct description of the boundary data retained
by the ordinary parent.  More generally, for $H=B+UV^T$ with an
arbitrary banded $B$, the same identity is $M=LBR+(LU)(V^TR)$;
the low-rank summand is again confined to that fixed corner.

\begin{corollary}\label{cor:band-core-stationary-parent}
Let $B_n:=T_n(b)$, where $b_s=0$ for $s\notin[-b_-,b_+]$, and let
\[
 H_n:=B_n+U_nV_n^T
\]
be leading sections of fixed arrays with $\ell$ generator columns.
Suppose that the columns of $U$ and $V$ satisfy fixed backward recurrences
\[
 \sum_{a=0}^{a_-}\alpha_a U_{i-a,:}=0,\qquad
 \sum_{b=0}^{a_+}\beta_b V_{j-b,:}=0,
 \qquad \alpha_0=\beta_0=1,
\]
whenever all participating indices are positive.  Define
\[
 P(z):=\sum_{a=0}^{a_-}\alpha_a z^a,\qquad
 Q(z):=\sum_{b=0}^{a_+}\beta_b z^b,
 \qquad p:=b_-+a_-,\quad q:=b_++a_+.
\]
For $n=1$ set $L_1=R_1=(1)$.  For $n\ge2$, the ordinary matrices
\begin{equation}\label{eq:band-core-stationary-parent}
\begin{split}
 L_n&:=\diag\bigl(1,P(\mathsf S_-^{(n-1)})\bigr),\qquad
 R_n:=\diag\bigl(1,Q(\mathsf S_+^{(n-1)})\bigr),\\
 M_n&:=L_nH_nR_n=T_n(c)+E_n,\qquad
 c_s:=\sum_{a=0}^{a_-}\sum_{b=0}^{a_+}
                 \alpha_a\beta_b b_{s+a-b}
\end{split}
\end{equation}
are leading sections of one $(p,q)$-banded array.  Their northwest
correction $E_n$ is the leading section of a fixed matrix supported in a
northwest $(p+1)\times(q+1)$ rectangle.  For $n\ge2$,
\[
 \operatorname{Cond}^{n-2}(M_n,J)
 =\operatorname{Cond}^{n-2}(H_n,J),\qquad \det M_n=\det H_n.
\]
Thus this construction gives an explicit ordinary ancestor of the complete
order-two descendant at depth $n-2$.  This is the constructed depth,
not a minimal-depth assertion for that target.  The determinant
sequence has an eventual constant-coefficient recurrence of order at most
\begin{equation}\label{eq:band-core-stationary-order}
 \binom{p+q}{p}.
\end{equation}
The recurrence is the characteristic-polynomial recurrence of the ordinary
row--column transfer for the bulk $c$, with initial boundary state specified
by $E$.
\end{corollary}

\begin{proof}
The constant recurrence makes $L_nU_n$ zero below row $a_-+1$ and
$V_n^TR_n$ zero after column $a_++1$.  Hence their product contributes only
to the northwest $(a_-+1)\times(a_++1)$ rectangle.  This rectangle is
contained in the $(p,q)$ band.

Away from the anchored boundary, multiplication of the three fixed-band
matrices $L_n,B_n,R_n$ gives precisely $c_{j-i}$.  To locate the boundary
correction, compare the truncated sums with the full double sum in
\eqref{eq:band-core-stationary-parent}.  A missing row term has
$i\le a_-+1$ and $i-a\le1$; for it to contribute, the band of $b$ requires
\[
 j-i+a-b\le b_+,\qquad j\le b_++a_++1=q+1.
\]
A missing column term analogously requires $j\le a_++1$ and
$i\le b_-+a_-+1=p+1$.  Thus the correction lies in the stated rectangle.
The product $L_nB_nR_n$ is $(p,q)$-banded, so the same is true of $M_n$.

The first coordinate is isolated in both unitriangular factors.
Proposition~\ref{prop:anchored-minor-preservation} proves equality of the
two complete descendant layers.  Once the growing ordinary determinant
passes the fixed corner, its row--column state evolves by the constant
bulk transfer of dimension $\binom{p+q}{p}$.  Cayley--Hamilton gives
\eqref{eq:band-core-stationary-order}.
\end{proof}

Every entry of the northwest correction is explicit.  In the local expression
\eqref{eq:local-vapor-entry}, replace $r_-,r_+$ by $a_-,a_+$ and use the constant coefficient $\alpha_a$ when
$a=0$ or $1\le a\le i-2$, and zero for a term that would cross the
isolated first coordinate; apply the analogous rule to $\beta_b$.
Subtracting $c_{j-i}$ gives $E_{ij}$.  Thus the exceptional entries of
this parent are determined by a fixed finite portion of the original
data.  The recurrence orders $a_\pm$ of the generators need not
equal their separated ranks: local rank and constant recurrence order are
different quantities.

\subsection{Ordinary parents for recurrent Toeplitz tails}
\label{sec:tailannihilator}

For a Toeplitz symbol $b=(b_s)_{s\in\mathbb Z}$, define
\begin{equation}\label{eq:toeplitzconnectedminor}
 D_t(s;b):=\det(b_{s+j-i})_{0\le i,j<t},\qquad D_0(s;b):=1.
\end{equation}
Fixed recurrences for the two symbol tails directly specify both the
compression depth and an anchored ordinary parent.  The two recurrence
orders and the two surviving widths may be different.

Suppose there are integers $m_+,m_-,r_+,r_-\ge0$ and coefficients
\[
 \beta_0=1,\beta_1,\ldots,\beta_{r_+},\qquad
 \alpha_0=1,\alpha_1,\ldots,\alpha_{r_-}
\]
such that
\begin{align}
 \sum_{q=0}^{r_+}\beta_q b_{s-q}&=0
 &&(s\ge m_++1),\label{eq:tailannpos}\\
 \sum_{a=0}^{r_-}\alpha_a b_{s+a}&=0
 &&(s\le-m_--1).\label{eq:tailannneg}
\end{align}

Apply the two recurrences while fixing the first coordinate.  For $n\ge2$
put
\begin{equation}\label{eq:anchored-tail-factors}
 \widehat L_n:=\diag\left(1,
   \sum_{a=0}^{r_-}\alpha_a(\mathsf S_-^{(n-1)})^a\right),\qquad
 \widehat R_n:=\diag\left(1,
   \sum_{q=0}^{r_+}\beta_q(\mathsf S_+^{(n-1)})^q\right),
\end{equation}
and set $\widehat L_1=\widehat R_1=(1)$.  Define
\begin{equation}\label{eq:anchored-tail-parent}
 M_n:=\widehat L_nT_n(b)\widehat R_n,
 \qquad
 c_s:=\sum_{a=0}^{r_-}\sum_{q=0}^{r_+}
       \alpha_a\beta_q b_{s+a-q}.
\end{equation}
The anchoring retains the connected boundary minors needed for a complete
order-two Dodgson block.  It therefore makes the resulting ordinary
matrix an exact ancestor of that block.

\begin{theorem}\label{thm:anchored-recurrent-tails}
Assume \eqref{eq:tailannpos}--\eqref{eq:tailannneg}, and put
$R:=\max(r_-,r_+)$.  For the connected minors in
\eqref{eq:toeplitzconnectedminor},
\begin{equation}\label{eq:anchored-tail-minor-support}
 D_t(s;b)=0\quad\text{if}\quad
 \bigl(t>r_+\text{ and }s>m_+\bigr)
 \quad\text{or}\quad
 \bigl(t>r_-\text{ and }s<-m_-\bigr).
\end{equation}
Consequently, $R$ condensation steps give a first layer supported on
$[-m_-,m_+]$, and $R+1$ steps give two layers with that support, whenever
the indicated layers are present.

The ordinary matrices $M_n$ are leading sections of one fixed array and satisfy
\begin{equation}\label{eq:anchored-tail-late-block}
 \operatorname{Cond}^{n-2}(M_n,J)
 =\operatorname{Cond}^{n-2}(T_n(b),J)\qquad(n\ge2).
\end{equation}
In particular, $\det M_n=\det T_n(b)$.  Its Toeplitz bulk is supported on
$[-m_-,m_+]$.  More precisely, if $\tau(i):=\max(i-2,0)$, then
\begin{equation}\label{eq:anchored-tail-entry}
 M_{ij}=
 \sum_{a=0}^{\min(r_-,\tau(i))}
 \sum_{q=0}^{\min(r_+,\tau(j))}
 \alpha_a\beta_q b_{j-i+a-q},
 \qquad c_s=0\quad(s>m_+\text{ or }s<-m_-).
\end{equation}
Thus $M_n=T_n(c)+E_n$, where $E_n$ is the leading section of a fixed
matrix and
\begin{equation}\label{eq:anchored-tail-corner}
 E_{ij}=M_{ij}-c_{j-i}=0
 \quad\text{if }i>K_-\text{ or }j>K_+,
 \qquad
 \begin{cases}
 K_-:=1+\max(r_-,r_++m_-),\\
 K_+:=1+\max(r_+,r_-+m_+).
 \end{cases}
\end{equation}
If $r_-\ge m_-$ and $r_+\ge m_+$, the corner is contained in the smaller
$(r_-+1)\times(r_++1)$ rectangle.
\end{theorem}

\begin{proof}
Suppose $t>r_+$ and $s>m_+$.  In the size-$t$ Toeplitz minor, take the
linear combination of its last $r_++1$ columns with coefficients
$\beta_0,\ldots,\beta_{r_+}$ in reverse column order.  In row $i$,
$0\le i<t$, it is
\[
 \sum_{q=0}^{r_+}\beta_q b_{s+t-1-i-q}=0,
\]
because $s+t-1-i\ge s\ge m_++1$.  The last column is therefore a linear
combination of preceding columns.  The negative-tail assertion follows
by combining the last $r_-+1$ rows: the recurrence index
$s+j-t+1$ is at most $s\le-m_--1$.  This proves
\eqref{eq:anchored-tail-minor-support} and the stated condensation depths.

The factors in \eqref{eq:anchored-tail-factors} are unitriangular and have
an isolated first identity coordinate.  Proposition~\ref{prop:anchored-minor-preservation}
with $h=1$ proves \eqref{eq:anchored-tail-late-block}, including the central
denominator minor.  Direct multiplication gives
\eqref{eq:anchored-tail-entry}.  All indices in this formula precede the
entry being computed, so it is independent of the chosen leading section.
Thus the same formula defines every leading section.

For $s>m_+$, apply \eqref{eq:tailannpos} at $s+a$ for each $a$ in the full
bulk sum.  It gives $c_s=0$.  For $s<-m_-$, apply
\eqref{eq:tailannneg} at $s-q$ for each $q$.  This proves the bulk support.
Both sums have their full ranges when $i>r_-+1$ and $j>r_++1$, so the
correction initially lies in the union of the first $r_-+1$ rows and the first
$r_++1$ columns.

If $i\le r_-+1$ and $j>K_+$, then the complete positive-tail recurrence is
applied and $j-i>m_+$.  For every retained row index $a$, the positive
recurrence annihilates its whole column-filter sum.  Both $M_{ij}$ and
$c_{j-i}$ are zero.  Transposition gives the corresponding assertion
when $j\le r_++1$ and $i>K_-$.  Together with the region where both sums are complete,
these facts prove \eqref{eq:anchored-tail-corner}.

Finally assume $r_+\ge m_+$ and $r_-\ge m_-$.  When $j>r_++1$, every
row-filter term missing at row $i\le r_-+1$ has $a\ge i-1$, and hence
$j-i+a\ge j-1>m_+$.  The complete positive-tail recurrence annihilates each such term.
Thus no correction remains above the bulk outside the first $r_++1$ columns.
The reflected argument proves the smaller rectangle claim.
\end{proof}

The bound~\eqref{eq:anchored-tail-corner} concerns the support of $E_n$;
it need not place that support inside the band of $T_n(c)$.  The
stronger conclusion that the whole ordinary parent is banded with the
stated widths is proved under the hypotheses of
Corollary~\ref{cor:band-core-stationary-parent}.

\subsection{Ordinary row--column recurrence input}\label{sec:background}

For fixed nonnegative integers $m_1,m_2$, let
\[
 T_n:=(a_{j-i})_{i,j=1}^n,
 \qquad a_s=0\quad\text{for }s\notin[-m_1,m_2].
\]
We use the ordinary row--column construction of
\cite{AlekseyevKhomovskyRecurrences2026} and the reflection refinements
of \cite{AlekseyevKhomovskySymmetries2026}.  These recurrence bounds have
classical antecedents: Shevelev's determinant algomatrix is a compound
of a companion-type operator, with fixed-cardinality root products and
the binomial state dimension \cite{Shevelev1991}.  The broader
restricted-position methods are surveyed in
\cite{Shevelev1992,Shevelev1993}.

\begin{theorem}
For fixed $m_1,m_2$, the sequence $(\det T_n)_{n\ge0}$ is annihilated by a
homogeneous constant-coefficient recurrence of order at most
\[
 N(m_1,m_2):=\binom{m_1+m_2}{m_1}.
\]
The recurrence can be constructed from a finite family of boundary minors.
The same boundary-state construction extends to position-dependent banded
matrices, in which case the transfer matrix depends on the position.
\end{theorem}

\paragraph{Scalar recurrences and nonvanishing conditions.}
Write the ordinary realization, after its fixed initial segment, as
\[
 z_{n+1}=\mathcal T_n z_n,\qquad D_n=\lambda_nz_n,
 \qquad z_n\in\F^d.
\]
For $j\ge0$ define row vectors
\begin{equation}\label{eq:local-vapor-observation-rows}
 \nu_{n,0}:=\lambda_n,\qquad
 \nu_{n,j}:=\lambda_{n+j}\mathcal T_{n+j-1}\cdots \mathcal T_n\quad(j\ge1).
\end{equation}
Since $D_{n+j}=\nu_{n,j}z_n$, these vectors recover successive
determinants from one common state; we call them \emph{observation rows}.
They involve products of successive transfer matrices.  For variable
data these factors must retain their respective position indices.

If the first $d$ rows form a nonsingular matrix
\[
 \mathcal O_n:=\begin{bmatrix}
 \nu_{n,0}\\\vdots\\\nu_{n,d-1}
 \end{bmatrix},
\]
then the scalar coefficients are explicitly
\begin{equation}\label{eq:local-vapor-scalar-coefficients}
 (c_{n,0},\ldots,c_{n,d-1}):=\nu_{n,d}\mathcal O_n^{-1},
 \qquad
 D_{n+d}=\sum_{j=0}^{d-1}c_{n,j}D_{n+j}.
\end{equation}
They are rational functions of a bounded number of consecutive entries of
the original data.

More generally, in the rational function field of independent formal
generator sequences, choose the first dependent row among
$\nu_{n,0},\ldots,\nu_{n,d}$.  Its index $\nu$ is at most $d$, and
its dependence on the preceding independent rows gives a monic
recurrence of order $\nu$, with coefficients found from a nonzero
minor of those rows.  The same formula shifts with $n$, because all
transfer formulas beyond the initial segment use the same relative
indices.  Hence the generic formal determinant sequence has a scalar
recurrence of order at most $d$ with rational local coefficients.

At a numerical or functional specialization this monic formula applies
where its chosen generator and observation denominators are nonzero.
When an observation minor vanishes, the state update remains
valid; clearing the denominators gives a scalar identity,
but its leading coefficient can vanish.  If the last observation row
still lies in the span of the preceding rows, a nonsingular smaller
minor gives a lower-order local formula.  Dimension alone does not
justify division by a vanished observation determinant or a monic
recurrence at every specialized index.

\paragraph{When constant coefficients are available.}
A scalar sequence has an eventual homogeneous linear recurrence with
constant coefficients precisely when its shifts span a finite-dimensional
vector space after a fixed initial segment is discarded.  Indeed, a monic
recurrence expresses all later shifts through finitely many preceding
ones; conversely, a dependence among finitely many shifts is such a
recurrence.  Equivalently, the sequence has a realization
$D_n=\lambda \mathcal T^{n-n_0}v$ for $n\ge n_0$, with $\mathcal T,\lambda,v$ independent
of $n$.  The forward implication uses the companion matrix, and the
reverse implication follows from Cayley--Hamilton.  This is the
constant-coefficient setting of \cite{KauersZeilberger2018}.
The position-dependent realization supplied by band compression need
not have this stronger property, even when each of its local coefficient
sequences is polynomial in the position; see
Example~\ref{ex:banded-not-constant-recurrent}.

Periodicity of the complete transfer and observation data is another
sufficient condition and gives an explicit cost in the order.

\begin{proposition}\label{prop:periodic-ordinary-transfer}
Suppose that, for $n\ge n_0$,
\[
 z_{n+1}=\mathcal T_n z_n,\qquad D_n=\lambda_nz_n,\qquad z_n\in\F^d,
\]
and that $\mathcal T_{n+p}=\mathcal T_n$ and $\lambda_{n+p}=\lambda_n$ for a fixed
integer $p\ge1$.  Then $(D_n)$ has an eventual homogeneous
constant-coefficient recurrence of order at most $pd$.
\end{proposition}

\begin{proof}
Use $p$ copies of $\F^d$, labelled $0,\ldots,p-1$.  Embed $z_n$ in
copy $s$, where $s\equiv n-n_0\pmod p$, and put zero in the other
copies.  One fixed linear map sends copy $s$ to copy $s+1$ modulo $p$
by $\mathcal T_{n_0+s}$.  One fixed row vector restricts to $\lambda_{n_0+s}$
on copy $s$.  These give a constant realization of dimension $pd$;
its characteristic polynomial annihilates the scalar sequence.
\end{proof}

For leading sections of one ordinary banded array, with a
size-independent northwest correction and diagonal sequences eventually
periodic with a common period $p$, the row--column transfer and output
maps are eventually periodic with period $p$.
If its lower and upper semibandwidths are $w_-,w_+$, one may take
$d=\binom{w_-+w_+}{w_-}$.  A fixed Toeplitz bulk is the case $p=1$.
Corollary~\ref{cor:band-core-stationary-parent} constructs that case
from a Toeplitz band core plus separated generators satisfying fixed
linear recurrences.  An arbitrary diagonal band is not covered by the
corollary, even when its entries and all generators are polynomial or
exponential-polynomial sequences.  A separate constant realization,
rather than the regularity of the individual entries alone, is what
justifies replacing variable coefficients by constant ones.

In Examples~\ref{ex:ranktwo-power-vapor}
and~\ref{ex:ranktwo-vapor-corrected}, the constant Toeplitz core and
fixed shift recurrences of the generators supply the extra hypothesis.
The smaller position-dependent state and a constant realization may
have different dimensions, which accounts for their different scalar
recurrence orders.

For a stationary banded bulk, the fixed northwest correction supplies
the initial boundary state.  Changing this size-independent correction
preserves the eventual recurrence supplied by the bulk transfer, although
special boundary data can lower the minimal scalar order.  A bounded
correction support alone is insufficient when its entries vary with the
total matrix size.

\paragraph{The recurrence of the constructed parent.}
The explicit parent in Theorem~\ref{thm:anchored-recurrent-tails} has
a fixed northwest correction that can extend outside its bulk band.
The following ordinary determinant lemma supplies the recurrence for
this parent and remains valid at repeated roots.  Its exact boundary
formula justifies the specializations used in the examples.

\begin{lemma}\label{lem:ordinary-fixed-corner-state}\label{lem:exact-boundary-companion}
Let $p,q\ge0$, let $c_s=0$ for $s\notin[-p,q]$, and let $E$ be a
fixed infinite matrix supported in a finite northwest rectangle.
Write $E_n:=E[1:n,1:n]$ and
\[
 M_n:=T_n(c)+E_n,\qquad D_n:=\det M_n.
\]
Then $(D_n)$ has an eventual homogeneous constant-coefficient recurrence
of order at most $\binom{p+q}{p}$.

If $p,q>0$, put $\ell:=p+q$ and suppose $c_q\ne0$.  When
one fixed leading block containing the correction is invertible, a constant
transfer is
\begin{equation}\label{eq:exact-boundary-companion-transfer}
 \gamma\Lambda^q C,\qquad \gamma:=(-1)^q c_q,
\end{equation}
where $C$ is the $\ell\times\ell$ forward companion defined by
\[
 C_{i,i+1}:=1\quad(1\le i<\ell),\qquad
 C_{\ell,j}:=-\frac{c_{j-p-1}}{c_q}\quad(1\le j\le\ell),
\]
and all remaining entries are zero.  This is the transpose of the
column-companion convention used in
\cite{AlekseyevKhomovskyRecurrences2026}.  Here $\Lambda^q C$ denotes the
$q$th exterior power, whose matrix entries in the standard basis are
the $q\times q$ minors of $C$.  No distinct-root hypothesis is required.
\end{lemma}

\begin{proof}
Choose a fixed $K$ so that the entire rectangular support of $E$ lies
in $[1:K]^2$.  If $q=0$, the top-right block of $M_{K+m}$ is zero;
if $p=0$, its bottom-left block is zero.  In either case
\[
 D_{K+m}=D_Kc_0^m\qquad(m\ge0),
\]
which proves the order-one bound, including $p=q=0$.

Now assume $p,q>0$ and work first where $c_q\det M_K\ne0$.
Eliminating $M_K$ leaves a Schur complement $S_m$ of order $m$ and
\[
 D_{K+m}=\det M_K\det S_m.
\]
The lower-left coupling to $M_K$ is confined to the first $p$
remaining rows, and the upper-right coupling is confined to the first
$q$ remaining columns.  Consequently $S_m$ has Toeplitz bulk $c$ and
a fixed correction supported in its northwest $p\times q$ rectangle.
This conclusion allows the original corner to extend outside the bulk
band.

Take $m\ge\max\{p,q\}$.  Extend the bulk row equations to variables
$z_1,\ldots,z_{m+q}$, with $z_i=0$ for $i\le0$.  The first $p$
row equations determine $z_{q+1},\ldots,z_{q+p}$ uniquely from
$z_1,\ldots,z_q$: their coefficient matrix on the variables being
solved for is triangular with diagonal $c_q$.  They therefore give a
fixed solution frame $V\in\F^{\ell\times q}$, whose first $q$ rows
form $I_q$.  Every later row advances the $\ell$ consecutive
coordinates by $C$.  With $P:=[0_{q\times p}\ I_q]$, the right
boundary map is
\[
 PC^{m-p}V.
\]

For the normalization, let $F$ map
$(z_1,\ldots,z_{m+q})$ to its first $q$ coordinates followed by the
$m$ row-equation residuals.  Its determinant is $c_q^m$, since each
successive residual has last coefficient $c_q$.  The block of
$F^{-1}$ taking the initial $q$ coordinates to the final $q$
coordinates is $PC^{m-p}V$.  The complementary-minor identity for
$F^{-1}$ therefore gives
\begin{equation}\label{eq:exact-boundary-shooting}
 \det S_m=\gamma^m\det(PC^{m-p}V).
\end{equation}
Indeed, its complementary minor in $F$ is $S_m$ and the sign is
$(-1)^{qm}$.  Taking $q$th exterior powers in
\eqref{eq:exact-boundary-shooting} expresses $D_{K+m}$ as a
fixed matrix coefficient of
$(\gamma\Lambda^q C)^{m-p}$, with the constant factor
$\det M_K\gamma^p$ absorbed into the boundary covector.  The state
space has dimension $\binom{p+q}{q}=\binom{p+q}{p}$, so
Cayley--Hamilton proves the asserted bound under these nonvanishing assumptions.

The starting size $K+\max\{p,q\}$ is independent of the entries.  Every minor
of order greater than $\binom{p+q}{p}$ in the corresponding eventual
Hankel matrix of $(D_n)$ is a polynomial in the finitely many bulk and
corner entries, and vanishes whenever $c_q\det M_K\ne0$.  It therefore
vanishes identically, including when $c_q=0$ or $M_K$ is singular.
Finite rank of that Hankel matrix gives a constant-coefficient
recurrence of at most the same order: its columns are the successive
shifts of the tail sequence and span a shift-invariant vector space of
that dimension.  This proves the specialization assertion without
requiring a nonsingular companion or distinct roots.  More generally,
polynomial determinant identities extend any resulting transfer
annihilator to singular leading blocks whenever
its coefficients are regular.
\end{proof}

For $p=q=r$, the matrix $C$ is the normalized forward companion of
$\sum_{s=-r}^r c_s z^{s+r}$.  In this case the boundary identity
\eqref{eq:exact-boundary-shooting} already holds for $m\ge r$, and the
ordinary determinant transfer is $(-1)^r c_r\Lambda^r C$.

\begin{theorem}
For $m\ge1$, an $(m,m)$-banded Toeplitz determinant has an annihilator of order at
most $\binom{2m}{m}$.  If the Toeplitz matrix is symmetric, the order
may be reduced to
\[
 \frac{3^m+1}{2};
\]
in the skew-symmetric case it may be reduced to
\[
 2\cdot 3^{m-1}.
\]
In characteristic zero, the earlier symmetry paper
\cite{AlekseyevKhomovskySymmetries2026} proves that, outside a proper
algebraic exceptional set, these reduced bounds are the
generic minimal recurrence degrees, using a Hankel--Vandermonde nonvanishing
argument.
\end{theorem}

We now apply the ordinary boundary transfer to the parent of
Theorem~\ref{thm:anchored-recurrent-tails}, including its fixed northwest
correction.

\begin{corollary}\label{cor:anchored-tail-recurrence}
Under the hypotheses of Theorem~\ref{thm:anchored-recurrent-tails},
$D_n:=\det T_n(b)$ is eventually annihilated by a homogeneous
constant-coefficient recurrence of order at most
\begin{equation}\label{eq:anchored-tail-order}
 \binom{m_-+m_+}{m_-}.
\end{equation}
This applies, in particular, when both symbol tails have rational
generating functions after deletion of finite prefixes.  If the original
Toeplitz family is symmetric, $m_-=m_+=m$, and
$\alpha_a=\beta_a$, with equal filter lengths, the bound improves to
\begin{equation}\label{eq:anchored-tail-symmetric-order}
 \frac{3^m+1}{2}.
\end{equation}
\end{corollary}

\begin{proof}
For $n\ge R+3$, first condense $R+1$ times to the two-layer banded
block, and then $n-R-3$ times to its order-two descendant.
Equation~\eqref{eq:anchored-tail-late-block} vaporizes that complete
descendant into $M_n$ at depth $n-2$.  The explicit bulk and corner in
\eqref{eq:anchored-tail-parent}--\eqref{eq:anchored-tail-corner} then give
the constant ordinary state of
Lemma~\ref{lem:exact-boundary-companion}, of dimension
\eqref{eq:anchored-tail-order}.  This proves the first bound.
The recurrence conditions \eqref{eq:tailannpos}--\eqref{eq:tailannneg}
are precisely the eventual coefficient recurrences supplied by rational
tail generating functions.

For the symmetric bound, consider first a freely parameterized symmetric
bulk and northwest corner.  On the open set where the outer bulk
coefficient is nonzero and the bulk companion has distinct roots, these
roots occur in reciprocal pairs
$\rho_1,\rho_1^{-1},\ldots,\rho_m,\rho_m^{-1}$.  The direct ordinary
boundary transfer of Lemma~\ref{lem:exact-boundary-companion} is a scalar
multiple of the middle compound.  Its eigenvalues, with their common
scalar suppressed, are among
\[
 \prod_{j=1}^m\rho_j^{e_j},\qquad
 e_j\in\{-1,0,1\},\qquad
 \#\{j:e_j\ne0\}\equiv m\pmod2.
\]
There are at most $(3^m+1)/2$ such products.  The companion and its compound
are diagonalizable on this open set, so this is a matrix annihilator
bound for every boundary state.  The finite Hankel identities are
polynomial in the bulk and corner entries and therefore extend to all
their values.  Specializing to the anchored family proves
\eqref{eq:anchored-tail-symmetric-order}, including repeated roots; $m=0$
is immediate.
\end{proof}

The fixed first coordinate preserves both layers of the order-two
block.  Appendix~\ref{app:scalar-comparisons} gives the unanchored,
determinant-preserving variant with its smaller boundary correction.

For the polynomial kernels below, the repeated-root part can be evaluated
without diagonalizing the companion.

\begin{lemma}\label{lem:exterior-jordan-length}
In characteristic zero, if $J_s(1)$ is a single unipotent Jordan block,
then the largest Jordan block of $\Lambda^k J_s(1)$ has size
\begin{equation}\label{eq:exterior-jordan-length}
 k(s-k)+1\qquad(0\le k\le s).
\end{equation}
\end{lemma}

\begin{proof}
The entries of $J_s(1)^n$ are $\binom n{j-i}$, with a negative lower
index interpreted as zero.  Every $k\times k$ minor has degree in $n$
at most $k(s-k)$.  For the first $k$ rows and last $k$ columns, put
$h:=s-k$.  The coefficient of the maximal power is
\[
 \det\left(\frac1{(h+j-i)!}\right)_{0\le i,j<k}
 =\frac{\prod_{j=0}^{k-1}j!}{\prod_{j=0}^{k-1}(h+j)!}\ne0.
\]
Here $1/t!$ is zero when $t<0$.  Factoring $1/(h+j)!$ from column $j$
leaves the falling-factorial polynomials in $h+j$, whose determinant is
a Vandermonde determinant.  Hence the matrix powers of the compound
have exact maximal polynomial degree $k(s-k)$.  Finite differences
prove \eqref{eq:exterior-jordan-length}.
\end{proof}

For a tensor product of $z$ two-dimensional unipotent blocks, the same
argument gives largest Jordan length $z+1$: its matrix powers have degree
at most $z$, and the product of the $z$ upper-right entries has degree $z$.

\subsubsection*{Products of established determinant transfers}

The rank estimates in Section~\ref{sec:structured-closure} control
condensation of a matrix product itself.  A second route becomes available
once the factors have ordinary stationary determinant transfers.
The following standard closure argument for constant-coefficient
recurrences \cite{KauersZeilberger2018} constructs the product transfer
without multiplying the matrix entries.

\begin{proposition}\label{prop:product-determinant-transfer}
Fix $k$ families of $n\times n$ matrices $H_n^{(j)}$, $1\le j\le k$.
Suppose that, for $n\ge n_0$, each determinant has a realization
\[
 z_{n+1}^{(j)}=\mathcal T_j z_n^{(j)},\qquad
 \det H_n^{(j)}=\lambda_jz_n^{(j)},
\]
where $z_n^{(j)}$ has dimension $d_j$, and the matrices $\mathcal T_j$ and
row vectors $\lambda_j$ are independent of $n$.  Then the sequence
\[
 h_n:=\det\bigl(H_n^{(1)}\cdots H_n^{(k)}\bigr)
\]
has a stationary realization with transfer
\begin{equation}\label{eq:product-determinant-transfer}
 \mathcal T:=\mathcal T_1\otimes\cdots\otimes \mathcal T_k
\end{equation}
and an eventual constant-coefficient annihilator of degree at most
$d_1\cdots d_k$.  In particular, constant-coefficient scalar factor recurrences of orders
$d_j$ give the same bound.
\end{proposition}

\begin{proof}
Put
\[
 z_n:=z_n^{(1)}\otimes\cdots\otimes z_n^{(k)},\qquad
 \lambda:=\lambda_1\otimes\cdots\otimes\lambda_k.
\]
The tensor product identity gives $z_{n+1}=\mathcal T z_n$, while determinant
multiplicativity gives
\[
 h_n=\prod_{j=1}^k\det H_n^{(j)}
     =\prod_{j=1}^k(\lambda_jz_n^{(j)})=\lambda z_n.
\]
Cayley--Hamilton supplies the characteristic-polynomial annihilator
of $\mathcal T$, whose dimension is $d_1\cdots d_k$.  A constant-coefficient scalar
recurrence has a stationary companion realization of its stated order,
proving the last assertion.
\end{proof}

No commutativity of the factors is required.  The degree bound is not
an assertion of minimality: repeated spectral products and unobserved
state components can reduce it.  Position-dependent factor transfers
instead give the exact state update
$\mathcal T_1(n)\otimes\cdots\otimes \mathcal T_k(n)$, with no automatic conclusion of
constant scalar coefficients.  Thus compression of a product and
stationarity of its determinant transfer remain separate assertions.

The product here is taken before changing the matrix size.  For
well-defined infinite products, let $A_n,B_n$ denote leading sections.
Then
\begin{equation}\label{eq:finite-product-tail-correction}
 (AB)_n=A_nB_n+
 A[1:n,n+1:\infty]\,B[n+1:\infty,1:n].
\end{equation}
The second term need not vanish.  For banded factors it is confined
to a bounded southeast corner and enters through terminal boundary
data, as in \cite{AlekseyevKhomovskyAlgebra2026}.

\section{Structured families and determinant recurrences}
\label{sec:constantapps}

We now apply the structural route to separated, trigonometric,
polynomial, and exponential kernels.  In each family the connected
minors determine when the two layers become banded, and an explicit
ordinary ancestor preserves their common order-two descendant.
Its band entries give the ordinary determinant transfer.  When these
entries are constant outside a fixed northwest corner, that corner
specifies the initial boundary state for a constant transfer.
The examples also distinguish position-dependent coefficients from a
constant realization and show how repeated bulk roots reduce the
scalar order below the general state bound.

\subsection{Separated tails and Toeplitz specializations}

\subsubsection{Independent triangular factors}

Three-term characteristic-polynomial recurrences for order-one
quasiseparable matrices are classical
\cite[Theorem~3.3 and Corollary~3.4]{EidelmanGohbergOlshevsky2005}.
The scalar recurrence below is the specialization with unit transition
generators.  We prove it here by constructing an ordinary parent and
establishing equality of both descendant layers.  This also records
precisely which structural information is retained by vaporisation.

Let $H_n=(a_{ij})_{1\le i,j\le n}$ be given by
\begin{equation}\label{eq:separable-matrix}
 a_{ij}:=\begin{cases}
   f_i g_j,&i<j,\\
   h_i,&i=j,\\
   \widehat f_i \widehat g_j,&i>j.
 \end{cases}
\end{equation}
Put
\begin{equation}\label{eq:separable-local-data}
 \tau_i^+:=f_i g_i,\qquad \tau_i^-:=\widehat f_i\widehat g_i,\qquad
 w_i:=g_i\widehat f_i,\qquad v_i:=f_i\widehat g_i,\qquad
 \varepsilon_i^\pm:=h_i-\tau_i^\pm.
\end{equation}
Thus $\tau_i^+\tau_i^-=v_iw_i$.  Assume first that $w_i\ne0$.
Write $D_{i,j}:=\det H_n[i:j,i:j]$, with $D_{i,i-1}:=1$, and let
$C_t(i,j)$ denote the connected minor of order $t$ beginning in row $i$
and column $j$.  Write
$(A^{(r)},B^{(r)}):=\operatorname{Cond}^r(H_n,J)$ and use
$a_s^{(r)}[i],b_s^{(r)}[i]$ for the signed diagonal entries of its layers.

\begin{proposition}\label{prop:separable-band}
For $t\ge2$, $C_t(i,j)=0$ whenever $|i-j|>1$.  For $t\ge1$,
\begin{equation}\label{eq:separable-side-minors}
\begin{aligned}
 C_t(i,i+1)&=(-1)^{t-1}f_i g_{i+t}
            \prod_{k=i+1}^{i+t-1}\varepsilon_k^+,\\
 C_t(i+1,i)&=(-1)^{t-1}\widehat f_{i+t}\widehat g_i
            \prod_{k=i+1}^{i+t-1}\varepsilon_k^-.
\end{aligned}
\end{equation}
Consequently one condensation step gives a tridiagonal first layer, and
two steps give two tridiagonal layers.  The first-step entries are
\[
\begin{aligned}
 a_0^{(1)}[i]&=h_i h_{i+1}-v_iw_{i+1},\\
 a_1^{(1)}[i]&=-f_i g_{i+2}\varepsilon_{i+1}^+,\qquad
 a_{-1}^{(1)}[i]=-\widehat f_{i+2}\widehat g_i\varepsilon_{i+1}^-.
\end{aligned}
\]
\end{proposition}

\begin{proof}
For a shift of at least two above the diagonal, the first two rows of a
connected block are proportional, since both lie wholly in the upper
separable part.  Transposition proves the lower claim.  The side formulas
are immediate for $t=1,2$.  For larger $t$, Desnanot--Jacobi at shift one
has a zero shift-two term and gives
\[
 C_t(i,i+1)C_{t-2}(i+1,i+2)
 =C_{t-1}(i,i+1)C_{t-1}(i+1,i+2).
\]
Induction proves the upper formula; the lower one follows in the same way.
The identities hold as polynomial identities, including at vanishing
intermediate minors.  The first and second layers after two steps contain
minors of orders three and two, so both are tridiagonal.
\end{proof}

The second-step first-layer main diagonal is also explicit:
\begin{equation}\label{eq:separable-second-diagonal}
\begin{split}
 a_0^{(2)}[i]={}&h_i h_{i+1}h_{i+2}
       -h_i v_{i+1}w_{i+2}-h_{i+2}v_iw_{i+1}\\
       &+v_iw_{i+2}
          (\tau_{i+1}^++\tau_{i+1}^--h_{i+1}).
\end{split}
\end{equation}
Its adjacent diagonals and the denominator layer are
\[
\begin{aligned}
 a_1^{(2)}[i]&=f_i g_{i+3}\varepsilon_{i+1}^+\varepsilon_{i+2}^+,\\
 a_{-1}^{(2)}[i]&=\widehat f_{i+3}\widehat g_i
                      \varepsilon_{i+1}^-\varepsilon_{i+2}^-,\\
 b_s^{(2)}[i]&=a_s^{(1)}[i+1]\qquad(s=-1,0,1).
\end{aligned}
\]
The shift in the last formula records the central placement of the
denominator layer.  These formulas specify the first fully banded block.

Define the regular tridiagonal filling by $t_0[1]:=h_1$ and
\begin{equation}\label{eq:separable-parent-data}
\begin{aligned}
 t_0[i]&:=h_i+
   (h_{i-1}-\tau_{i-1}^+-\tau_{i-1}^-)
       \frac{w_i}{w_{i-1}},\qquad i\ge2,\\
 t_1[i]&:=-\varepsilon_i^+\frac{g_{i+1}}{g_i},\qquad
 t_{-1}[i]:=-\varepsilon_i^-\frac{\widehat f_{i+1}}{\widehat f_i},\\
 p_i&:=t_1[i]t_{-1}[i]
      =\varepsilon_i^+\varepsilon_i^-\frac{w_{i+1}}{w_i}.
\end{aligned}
\end{equation}
Let $T_n$ be the ordinary tridiagonal matrix with signed diagonal
entries $t_{-1}[i],t_0[i],t_1[i]$.  Set $M_1:=[h_1]$; for $n\ge2$,
obtain $M_n$ from $T_n$ by replacing its northwest $2\times2$ block by
\[
 \begin{pmatrix}h_1&f_1g_2\\
 \widehat f_2\widehat g_1&h_2\end{pmatrix}.
\]
Thus $m_s[i]=t_s[i]$ outside this corner.  The product $p_i$ is retained
only as a balanced product of the adjacent diagonals.

\begin{theorem}\label{thm:separable-corner-vapor}
For $n\ge3$,
\begin{equation}\label{eq:separable-corner-equality}
 \operatorname{Cond}^{n-2}(H_n,J)
 =\operatorname{Cond}^{n-2}(M_n,J).
\end{equation}
The matrices $M_n$ are leading sections of one fixed array.  Its
determinants equal $D_n:=\det H_n$ and satisfy the ordinary row--column
recurrence
\begin{equation}\label{eq:separable-recurrence}
 D_n=t_0[n]D_{n-1}-p_{n-1}D_{n-2}\quad(n\ge2),
 \qquad D_0=1,\qquad D_1=h_1.
\end{equation}
The discrepancy $M_n-T_n$ is supported in the northwest
$2\times2$ block.
\end{theorem}

\begin{proof}
For each interval $[i,j]$, let $K_{i,j}$ denote the ordinary continuant with
first main diagonal entry $h_i$, subsequent entries
$t_0[i+1],\ldots,t_0[j]$, and edge products $p_i,\ldots,p_{j-1}$.
For $j\ge i+1$, define $\widetilde M_{i,j}$ by replacing its first two main diagonal entries by
$h_i,h_{i+1}$ and its first upper and lower entries by
$f_i g_{i+1},\widehat f_{i+1}\widehat g_i$.  Its leading determinants equal the corresponding
$K$'s: the initial one-by-one determinants agree, and
\begin{align*}
 t_0[i+1]h_i-p_i
 &=h_i h_{i+1}
 +\bigl[h_i(h_i-\tau_i^+-\tau_i^-)
        -(h_i-\tau_i^+)(h_i-\tau_i^-)\bigr]
          \frac{w_{i+1}}{w_i}\\
 &=h_i h_{i+1}-\tau_i^+\tau_i^-\frac{w_{i+1}}{w_i}
 =h_i h_{i+1}-v_iw_{i+1}.
\end{align*}
All subsequent last-row expansions use the same main-diagonal entries and edge products.
The trailing principal interval of $\widetilde M_{i,j}$ beginning at its
second vertex is the continuant $K_{i+1,j}$ itself.

Desnanot--Jacobi for $\widetilde M_{i,j}$ therefore becomes
\begin{equation}\label{eq:separable-parent-dj}
 K_{i,j}K_{i+1,j-1}
 =K_{i,j-1}K_{i+1,j}
  -v_iw_j\prod_{k=i+1}^{j-1}
               \varepsilon_k^+\varepsilon_k^-.
\end{equation}
Indeed, its side-minor product is the product of all its edge products,
\[
 v_iw_{i+1}\prod_{k=i+1}^{j-1}
      \left(\varepsilon_k^+\varepsilon_k^-
            \frac{w_{k+1}}{w_k}\right)
 =v_iw_j\prod_{k=i+1}^{j-1}
               \varepsilon_k^+\varepsilon_k^-.
\]
Proposition~\ref{prop:separable-band} shows that the original interval
determinants $D_{i,j}$ satisfy the same identity.  Their empty and
one-by-one values agree.  Induction on interval length thus gives
$D_{i,j}=K_{i,j}$ when the intermediate minors are nonzero, and then as a rational identity
in the original parameters.  The result consequently includes vanishing
intermediate minors whenever the $w$-denominators remain nonzero.

The first layer of the order-two descendant contains principal entries
$D_{1,n-1},D_{2,n}$, while the second layer is $D_{2,n-1}$.  These agree
for $H_n$ and $M_n$ by the preceding interval identity.  The upper side
minor of $M_n$ is
\[
 f_1g_2\prod_{k=2}^{n-1}
       \left(-\varepsilon_k^+\frac{g_{k+1}}{g_k}\right)
 =(-1)^{n-2}f_1g_n\prod_{k=2}^{n-1}\varepsilon_k^+,
\]
and its lower side minor is the reflected expression with $\widehat f,\widehat g$ and
$\varepsilon^-$.  These are precisely the two side entries from
\eqref{eq:separable-side-minors}.  This proves equality of both layers in
\eqref{eq:separable-corner-equality}.  The ordinary last-row expansion
of $M_n$ now gives \eqref{eq:separable-recurrence} for $n\ge3$; the
initial two-by-two calculation proves it also for $n=2$.
\end{proof}

For $n\ge4$, the route takes two steps to the explicit tridiagonal
block, another $n-4$ steps to order two, and vaporisation of depth $n-2$
into $M_n$.  The coefficients in
\eqref{eq:separable-recurrence} are generally position dependent;
constant coefficients require additional assumptions on the local data.

\subsubsection{Geometric Toeplitz tails}\label{sec:geometrictails}

Let
\begin{equation}\label{eq:Btoeplitz}
 B_n:=\begin{bmatrix}
 b_0&b_1&b_2&\cdots&b_{n-1}\\
 b_{-1}&\ddots&\ddots&\ddots&\vdots\\
 b_{-2}&\ddots&\ddots&\ddots&b_2\\
 \vdots&\ddots&\ddots&\ddots&b_1\\
 b_{-n+1}&\cdots&b_{-2}&b_{-1}&b_0
 \end{bmatrix},
 \qquad
 b_{\pm s}:=b_{\pm1}\rho_\pm^{s-1}\quad(s\ge1).
\end{equation}
When $b_1b_{-1}\ne0$, the same tails can be written as
\[
 \rho_+:=b_2/b_1,\qquad \rho_-:=b_{-2}/b_{-1},\qquad
 b_{\pm s}=\frac{b_{\pm2}^{s-1}}{b_{\pm1}^{s-2}}.
\]
Define the ordinary parent bulk by
\begin{align}\label{eq:geometric-vapor-parameters}
 m_1&:=b_1-b_0\rho_+,&
 m_{-1}&:=b_{-1}-b_0\rho_-,\notag\\
 m_0&:=b_0-b_1\rho_--b_{-1}\rho_+
                   +b_0\rho_+\rho_-.
\end{align}
Write $D_n:=\det B_n$.  One condensation step has first-layer diagonals
\[
 a_0=b_0^2-b_{-1}b_1,
 \qquad a_1=b_1^2-b_0b_2=b_1m_1,
 \qquad a_{-1}=b_{-1}^2-b_0b_{-2}=b_{-1}m_{-1}.
\]
Its denominator layer is the dense central matrix $B_{n-2}$.  After two
steps, both layers are tridiagonal, with entries
\begin{equation}\label{eq:geometric-two-layer-values}
 \begin{array}{c|ccc}
 &s=-1&s=0&s=1\\ \hline
 A^{(2)}&b_{-1}m_{-1}^2&D_3&b_1m_1^2\\
 B^{(2)}&b_{-1}m_{-1}&D_2&b_1m_1
 \end{array}
\end{equation}
where
\begin{align}\label{eq:geometric-small-determinants}
 D_2&=b_0^2-b_1b_{-1},\notag\\
 D_3&=b_0^3-2b_0b_1b_{-1}
       +b_1^2b_{-2}+b_{-1}^2b_2-b_0b_2b_{-2}.
\end{align}
Write $C_t(s):=D_t(s;b)$ for the connected Toeplitz minor defined
in \eqref{eq:toeplitzconnectedminor}.  For every $t\ge1$ the two
side minors are
\begin{equation}\label{eq:geometric-shifted-minors}
 C_t(1)=b_1m_1^{t-1},\qquad C_t(-1)=b_{-1}m_{-1}^{t-1},
\end{equation}
and for $t\ge2$ all shifts of absolute value greater than one vanish.
These are the Toeplitz specialization of
Proposition~\ref{prop:separable-band}.

\begin{proposition}\label{prop:Brec}
Let $M_n$ be the ordinary tridiagonal matrix with main diagonal
$b_0,b_0,m_0,\ldots,m_0$, superdiagonal $b_1,m_1,\ldots,m_1$,
and subdiagonal $b_{-1},m_{-1},\ldots,m_{-1}$.  For $n\ge3$,
\begin{equation}\label{eq:geometric-corner-cond}
 \operatorname{Cond}^{n-2}(B_n,J)
 =\operatorname{Cond}^{n-2}(M_n,J).
\end{equation}
Relative to the ordinary tridiagonal Toeplitz bulk, its fixed northwest
correction is
\begin{equation}\label{eq:geometric-northwest-defect}
 E[1:2,1:2]=
 \begin{bmatrix}b_0-m_0&b_1-m_1\\b_{-1}-m_{-1}&b_0-m_0\end{bmatrix},
 \qquad E_{ij}=0\quad\text{outside this corner}.
\end{equation}
Consequently, with $D_0=1$ and $D_1=b_0$,
\begin{equation}\label{eq:geometric-continuant-recurrence}
 D_n=m_0 D_{n-1}-m_1m_{-1} D_{n-2}\qquad(n\ge2).
\end{equation}
In terms of $b_0,b_{\pm1},b_{\pm2}$, when $b_1b_{-1}\ne0$, this is
\begin{align}\label{eq:recB}
 D_n={}&
 \frac{b_{-1}b_0b_1-b_{-1}^2b_2-b_{-2}b_1^2+b_{-2}b_0b_2}
      {b_{-1}b_1}D_{n-1}\notag\\
 &-\frac{(b_1^2-b_0b_2)(b_{-1}^2-b_0b_{-2})}
         {b_{-1}b_1}D_{n-2}.
\end{align}
\end{proposition}

\begin{proof}
First suppose $\rho_+\rho_-\ne0$ and represent the upper and lower triangular
parts as in \eqref{eq:separable-matrix} by
\[
 f_i:=b_1\rho_+^{-i},\quad g_i:=\rho_+^{i-1},\qquad
 \widehat f_i:=\rho_-^{i-1},\quad
 \widehat g_i:=b_{-1}\rho_-^{-i},\quad h_i:=b_0.
\]
The local ratios $w_i/w_{i-1}$ equal $\rho_+\rho_-$, and substitution
in \eqref{eq:separable-parent-data} gives exactly
$(m_0,m_1,m_{-1})$ from \eqref{eq:geometric-vapor-parameters}.
Theorem~\ref{thm:separable-corner-vapor} therefore supplies the exact
ordinary ancestor $M_n$, including equality of both layers in
\eqref{eq:geometric-corner-cond}.  Its entries are independent of the
total size, and its last-row expansion is
\eqref{eq:geometric-continuant-recurrence} for $n\ge3$.  For $n=2$,
\[
 m_0 b_0-m_1m_{-1}=b_0^2-b_1b_{-1}=D_2.
\]
All formulas in $b_0,b_1,b_{-1},\rho_+,\rho_-$ are polynomial.
Consequently the identities also hold at zero tail ratios or zero
intermediate minors.  Substitution of $\rho_+=b_2/b_1$ and
$\rho_-=b_{-2}/b_{-1}$ gives \eqref{eq:recB}.
\end{proof}

If $b_{-2}=b_2=0$, \eqref{eq:recB} specializes to the usual
tridiagonal recurrence.  If $b_{-2}=b_2=b_{-1}=b_1$, the matrix has
constant off-diagonal entries and
\[
 D_n=2(b_0-b_1)D_{n-1}-(b_0-b_1)^2D_{n-2},\qquad
 D_n=(b_0+(n-1)b_1)(b_0-b_1)^{n-1}.
\]

\subsubsection{A tridiagonal core}
\label{sec:rankone-band-core-example}

The prescribed band may contain more than the main diagonal.  Fix $h_0,h_1,h_{-1},u,v$
and set, at positive indices,
\begin{equation}\label{eq:rankone-band-core-input}
 H_{ij}:=\begin{cases}
 h_0,&i=j,\\
 h_1,&j=i+1,\\
 h_{-1},&i=j+1,\\
 u^iv^j,&|i-j|>1.
 \end{cases}
\end{equation}
Here $b_-=b_+=1$ and $r_-=r_+=1$.  One condensation therefore gives a
pentadiagonal first layer, although its denominator is still the central
section of the dense original matrix.  Write $(A^{(1)},B^{(1)}):=\operatorname{Cond}(H_n,J)$.  Its five first-layer diagonals are
\begin{align}\label{eq:rankone-band-core-first-layer}
 a_0^{(1)}[i]&=h_0^2-h_1h_{-1},\\
 a_1^{(1)}[i]&=h_1^2-h_0u^iv^{i+2},&
 a_{-1}^{(1)}[i]&=h_{-1}^2-h_0u^{i+2}v^i,\notag\\
 a_2^{(1)}[i]&=u^iv^{i+3}(u^{i+1}v^{i+2}-h_1),&
 a_{-2}^{(1)}[i]&=u^{i+3}v^i(u^{i+2}v^{i+1}-h_{-1}).\notag
\end{align}
Every other entry is zero.  These formulas follow by expanding the
corresponding $2\times2$ minors.  At the next step both layers are
pentadiagonal, with first layer $C_3(i,j)$ and second layer
$C_2(i+1,j+1)$.

For the values $h_0=-3$, $h_1=h_{-1}=7$, $u=3$, $v=5$, the first entries are
\[
 a_0^{(1)}[1]=-40,\quad a_1^{(1)}[1]=1174,\quad a_{-1}^{(1)}[1]=454,
 \qquad a_2^{(1)}[1]=2096250.
\]
Thus a tridiagonal core enlarges the surviving semibandwidth from one
to two without delaying the first compression step.

The ordinary ancestor of the order-two descendant is also explicit.
Use
\[
 L_n:=\diag(1,I-u\mathsf S_-^{(n-1)}),\qquad
 R_n:=\diag(1,I-v\mathsf S_+^{(n-1)}).
\]
For $M_n:=L_nH_nR_n$, every entry is the four-term expression
\begin{equation}\label{eq:rankone-band-core-parent-entry}
\begin{split}
 (M_n)_{ij}={}&H_{ij}
 -u\,\mathbf1_{i\ge3}H_{i-1,j}
 -v\,\mathbf1_{j\ge3}H_{i,j-1}\\
 &+uv\,\mathbf1_{i\ge3}\mathbf1_{j\ge3}H_{i-1,j-1},
\end{split}
\end{equation}
where terms multiplied by a zero indicator are omitted.  The parent is
pentadiagonal, and Theorem~\ref{thm:local-anchored-vaporisation} gives
equality of the complete order-two descendants at depth $n-2$.  The ordinary
row--column state has dimension six.  For general $u,v$ its band entries
depend on position, so the scalar recurrence is obtained from successive transfer matrices
by \eqref{eq:local-vapor-observation-rows}--\eqref{eq:local-vapor-scalar-coefficients}.  When $uv=1$, the original family is
Toeplitz and the delayed-geometric specialization below
gives a constant transfer.  A constant diagonal by itself does not make
the general $u^iv^j$ family Toeplitz.

\subsubsection{Delayed geometric tails}

The same route applies when the recurrent tail begins one diagonal later
than in Section~\ref{sec:geometrictails}.  Let
\begin{equation}\label{eq:delayed-geometric-symbol}
 b_{\pm s}:=b_{\pm2}\rho_\pm^{s-2}\qquad(s\ge2),
 \qquad D_n:=\det T_n(b).
\end{equation}
The seven parameters $b_0,b_{\pm1},b_{\pm2},\rho_+,\rho_-$ are independent.

\begin{proposition}\label{prop:delayed-geometric-vapor}
One condensation step gives a pentadiagonal first layer, with entries
\begin{align}\label{eq:delayed-geometric-first-layer}
 a_0^{(1)}&:=b_0^2-b_1b_{-1},\notag\\
 a_{\pm1}^{(1)}&:=b_{\pm1}^2-b_0b_{\pm2},\qquad
 a_{\pm2}^{(1)}:=b_{\pm2}^2-b_{\pm1}b_{\pm3}.
\end{align}
Its denominator is the central section of $T_n(b)$; a second step makes
both layers pentadiagonal.

An exact ordinary ancestor of the complete order-two descendant is
\begin{equation}\label{eq:delayed-geometric-parent}
 M_n:=\diag(1,I-\rho_-\mathsf S_-^{(n-1)})\,
 T_n(b)\,
 \diag(1,I-\rho_+\mathsf S_+^{(n-1)}).
\end{equation}
It is pentadiagonal, with Toeplitz bulk
\begin{equation}\label{eq:delayed-geometric-bulk}
\begin{aligned}
 m_0&:=(1+\rho_-\rho_+)b_0-\rho_-b_1-\rho_+b_{-1},\\
 m_{\pm1}&:=(1+\rho_-\rho_+)b_{\pm1}
             -\rho_{\mp}b_{\pm2}-\rho_\pm b_0,\\
 m_{\pm2}&:=b_{\pm2}-\rho_\pm b_{\pm1},
 \qquad m_s:=0\quad(|s|>2),
\end{aligned}
\end{equation}
and the following northwest correction, extended by zero:
\begin{equation}\label{eq:delayed-geometric-corner}
 E=
 \begin{pmatrix}
 b_0-m_0&b_1-m_1&0\\
 b_{-1}-m_{-1}&b_0-m_0&\rho_-m_2\\
 0&\rho_+m_{-2}&0
 \end{pmatrix}.
\end{equation}
Put $\omega:=m_{-2}m_2$.  The determinant sequence has an eventual
annihilator of order at most six, explicitly
\begin{align}\label{eq:delayed-geometric-annihilator}
 P_6(z):={}&z^6-m_0z^5+(m_{-1}m_1-\omega)z^4\notag\\
 &-\bigl(m_{-2}m_1^2+m_2m_{-1}^2-2m_0\omega\bigr)z^3\notag\\
 &+\omega(m_{-1}m_1-\omega)z^2-m_0\omega^2z+\omega^3.
\end{align}
Writing $(\mathcal E D)_n:=D_{n+1}$, one may take
$P_6(\mathcal E)D_n=0$ for $n\ge5$.
\end{proposition}

\begin{proof}
For $|s|\ge3$, the geometric relation gives
$b_s^2-b_{s-1}b_{s+1}=0$; the other five entries give
\eqref{eq:delayed-geometric-first-layer}.  Its first layer retains this
support under another condensation, while its preceding first layer
becomes the denominator.  This proves the stated condensation depths.

The recurrences $b_s-\rho_+b_{s-1}=0$ for $s\ge3$ and
$b_s-\rho_-b_{s+1}=0$ for $s\le-3$ have $m_-=m_+=2$ and
$r_-=r_+=1$.  Theorem~\ref{thm:anchored-recurrent-tails} therefore gives
the exact ancestor \eqref{eq:delayed-geometric-parent} and its complete
descendant identity \eqref{eq:anchored-tail-late-block}.  Away from the boundary, the complete sum is
\[
 m_s=(1+\rho_-\rho_+)b_s-\rho_-b_{s+1}-\rho_+b_{s-1}.
\]
Substitution gives \eqref{eq:delayed-geometric-bulk}.  The truncated first
two rows and columns give \eqref{eq:delayed-geometric-corner}; all other
entries equal their bulk values.  In particular the actual corner is
smaller than the general rectangle bound.

Apply the ordinary row--column transfer to this pentadiagonal ancestor.
For an explicit annihilator, assume initially $m_2\ne0$ and
$\det M_3\ne0$.  Lemma~\ref{lem:exact-boundary-companion} uses $K=3$ and
$r=2$, so from $n=5$ onward its transfer is $m_2\Lambda^2 C$, where $C$
is the four-dimensional forward-shift companion of
\[
 m_2z^4+m_1z^3+m_0z^2+m_{-1}z+m_{-2}.
\]
Its characteristic polynomial is exactly
\eqref{eq:delayed-geometric-annihilator}, as follows by taking the six
pairwise products of the four companion roots, with the common factor
$m_2$.  This proves the stated annihilation.  Its coefficients and all
$D_n$ are polynomial in the seven original parameters, so the identity
continues to hold when either of the assumed nonzero quantities vanishes.  For $n\ge4$, the full structural
route consists of two condensation steps to a two-layer pentadiagonal
block, $n-4$ further steps to order two, and vaporisation of depth $n-2$
into \eqref{eq:delayed-geometric-parent}.
\end{proof}

When $b_{\pm2}\ne0$, the equivalent tail formula is
\[
 b_{\pm s}=\frac{b_{\pm3}^{\,s-2}}{b_{\pm2}^{\,s-3}}\qquad(s\ge3).
\]
The exponent shift is fixed by the required identity
$b_{\pm3}^2-b_{\pm2}b_{\pm4}=0$ at the first omitted diagonal.

\subsection{Polynomial generators and variable coefficients}

\begin{example}\label{ex:ranktwo-power-vapor}
Let $p_1,p_2,q_1,q_2$ be nonnegative integers and consider
\[
 H_{ij}:=w_1i^{p_1}j^{q_1}+w_2i^{p_2}j^{q_2}
             -x\,\mathbf1_{i=j}.
\]
Assume first that $p_1\ne p_2$, $q_1\ne q_2$, and $w_1w_2\ne0$,
and work in characteristic zero.  Define
\[
 \Delta_p(a,b):=a^{p_1}b^{p_2}-a^{p_2}b^{p_1},\qquad
 \Delta_q(a,b):=a^{q_1}b^{q_2}-a^{q_2}b^{q_1}.
\]
The local ordinary parent is completely specified by
\eqref{eq:local-vapor-entry}, the anchored initial blocks, and, for
$i\ge4$,
\[
 \alpha_{i,1}=-\frac{\Delta_p(i-2,i)}{\Delta_p(i-2,i-1)},\qquad
 \alpha_{i,2}=\frac{\Delta_p(i-1,i)}{\Delta_p(i-2,i-1)}.
\]
The coefficients $\beta_{j,1},\beta_{j,2}$ are the same expressions
with $p,i$ replaced by $q,j$.  These identities follow by solving the
two local generator equations.  Their denominators are nonzero at the
positive indices in question under the stated exponent assumptions.

Two condensation steps give a pentadiagonal first layer, and three give
a full pentadiagonal block.  Theorem~\ref{thm:local-anchored-vaporisation}
then supplies an ordinary pentadiagonal parent of the order-two
descendant.  Its ordinary six-state transfer and the observation rows
\eqref{eq:local-vapor-observation-rows} determine a recurrence of order
at most six where the required observation minor is nonzero.
For these fixed exponents all entries and observation rows are rational
functions of the index and the parameters.  Their first row dependence
therefore gives rational coefficient functions after specialization to
this family; one must compute that dependence for the family rather than
divide a generic observation determinant that vanishes identically there.
Equal exponents or zero weights reduce the separated rank and are treated
by the lower-rank construction.

An explicit scalar recurrence is available for these fixed exponents.
For every nonnegative integer $t$, put
\[
 S_t(n):=\sum_{i=1}^n i^t,\qquad S_t(0):=0,
\]
and define the polynomial
\begin{align}\label{eq:ranktwo-power-polynomial}
 G(n):={}&x^2-x\bigl(w_1S_{p_1+q_1}(n)
                         +w_2S_{p_2+q_2}(n)\bigr)\notag\\
 &+w_1w_2\bigl(S_{p_1+q_1}(n)S_{p_2+q_2}(n)
                   -S_{p_1+q_2}(n)S_{p_2+q_1}(n)\bigr).
\end{align}
Expansion through the rank-two generators gives, for $x\ne0$,
\[
 D_n:=\det H_n=(-x)^{n-2}G(n).
\]
This also holds at $n=0,1$ after cancellation.  Put
\[
 L:=p_1+p_2+q_1+q_2+3.
\]
Each $S_t(n)$ has degree $t+1$ and leading coefficient $1/(t+1)$.
Consequently $G$ has degree $L-1$: its leading coefficient is
\[
 \frac{w_1w_2(p_1-p_2)(q_1-q_2)}
      {\prod_{s=1}^2\prod_{t=1}^2(p_s+q_t+1)},
\]
which is nonzero under the stated hypotheses.  Taking $L$ successive
differences of the polynomial factor gives the explicit recurrence
\begin{equation}\label{eq:ranktwo-power-constant-recurrence}
 \sum_{j=0}^{L}\binom Lj x^jD_{n-j}=0,
 \qquad n\ge L.
\end{equation}
Its minimal constant-coefficient order is $L$ when $x\ne0$.
This order can exceed six: the earlier six-state bound allows variable
coefficients and is a different assertion.

The compact six-term formula for geometric generators is given in
Proposition~\ref{prop:geometric-ranktwo-six}.  It requires constant
geometric ratios; substituting the current ratios of these power
functions into that formula is not valid.
\end{example}

\begin{example}\label{ex:ranktwo-vapor-corrected}
Work in characteristic zero with $x\ne0$, and put
\[
 H_n:=\bigl(\xi+ij-x\,\mathbf1_{i=j}\bigr)_{1\le i,j\le n},
 \qquad D_n:=\det H_n.
\]
The two generators $1,i$ satisfy the second-difference relation.
Choose
\[
 S_n:=\diag\bigl(1,(I-\mathsf S_-^{(n-1)})^2\bigr),
 \qquad M_n:=S_nH_nS_n^T.
\]
This is the allowed anchored initial-block choice in
Theorem~\ref{thm:local-anchored-vaporisation}.  It gives an exact
ordinary pentadiagonal ancestor of the order-two descendant.
Its symmetric Toeplitz bulk has coefficients
\[
 (c_0,c_1,c_2)=(-6x,4x,-x),
\]
and its northwest correction is the fixed $3\times3$ matrix
\begin{equation}\label{eq:ranktwo-vapor-corner}
 E=\begin{bmatrix}
 \xi+1+5x&\xi+2-4x&x-\xi-1\\
 \xi+2-4x&\xi+4+5x&-\xi-2-2x\\
 x-\xi-1&-\xi-2-2x&\xi+1+x
 \end{bmatrix}.
\end{equation}
For smaller orders, take the corresponding leading corner section.

The bulk polynomial is $-x(z-1)^4$.  The four-dimensional
bulk companion matrix may therefore be taken as
\[
 T:=\begin{bmatrix}
 0&1&0&0\\0&0&1&0\\0&0&0&1\\-1&4&-6&4
 \end{bmatrix}.
\]
Since $T$ is the companion of $(z-1)^4$, it is similar to $J_4(1)$.
Lemma~\ref{lem:exterior-jordan-length} gives
$\bigl(\Lambda^2T-I\bigr)^5=0$.
Lemma~\ref{lem:exact-boundary-companion} identifies the determinant
transfer with this compound multiplied by $-x$, including at the
repeated-root bulk.  Hence the ordinary row--column sequence, and
hence $D_n$, satisfies
\begin{equation}\label{eq:ranktwo-vapor-recurrence}
 \sum_{j=0}^5\binom5j x^jD_{n-j}=0
\end{equation}
after the fixed initial boundary segment.

For comparison, expansion through the two generators gives
\begin{equation}\label{eq:ranktwo-vapor-comparison}
 D_n=(-x)^{n-2}\left(
 x^2-x\left(\xi n+\frac{n(n+1)(2n+1)}6\right)
       +\xi\frac{n^2(n^2-1)}{12}\right).
\end{equation}
It verifies \eqref{eq:ranktwo-vapor-recurrence} for every $n\ge5$.
For $\xi\ne0$, the polynomial factor has degree four, so the minimal
constant-coefficient order is exactly five.  At $\xi=0$ its degree is
three, and the minimal constant-coefficient order drops to four:
\begin{equation}\label{eq:ranktwo-vapor-rankone-specialization}
 \sum_{j=0}^4\binom4j x^jD_{n-j}=0,
 \qquad n\ge4.
\end{equation}
The same pentadiagonal bulk remains available; this order drop is
caused by its boundary data.  Alternatively, the rank-one construction
at $\xi=0$ gives a tridiagonal position-dependent realization.
Equation~\eqref{eq:rankone-characteristic-vapor} then reads
\begin{equation}\label{eq:ranktwo-xizero-variable-recurrence}
 D_n=-x\left(1+\frac{n^2}{(n-1)^2}\right)D_{n-1}
       -x^2\frac{n^2}{(n-1)^2}D_{n-2},
 \qquad n\ge2,
\end{equation}
with $D_0=1$ and $D_1=1-x$.  Its order two permits variable
coefficients, whereas the minimal constant-coefficient order is four.
\end{example}

\begin{remark}\label{rem:rankone-different-main-diagonal}
At $\xi=0$, the main-diagonal entries of
Example~\ref{ex:ranktwo-vapor-corrected} are $i^2-x$.
The matrix with main-diagonal entries $x$ and off-diagonal entries $ij$
is instead
\[
 \widehat H_n=\diag(x-1^2,\ldots,x-n^2)+(ij)_{i,j=1}^n.
\]
It is the case $w=1$, $p=q=1$ of
Example~\ref{ex:separable-power-vapor}, whose displayed second-order
recurrence has coefficients depending on $n$.
These are different matrix families, as well as different types of
recurrence from the constant-coefficient orders five and four above.
\end{remark}

\subsection{A cosine-sum kernel}\label{sec:cosine-vaporisation}

Put $u_t:=\cos(Ct)$, $c:=\cos C$, and $\sigma:=\sin C$.  The
unperturbed matrix $K_n:=(u_{i+j})_{1\le i,j\le n}$ has rank at most
two, since
\[
 u_{i+j}=\cos(Ci)\cos(Cj)-\sin(Ci)\sin(Cj).
\]
Consequently $\det K_n=0$ for $n\ge3$.  A nontrivial determinant
sequence is obtained by adding a scalar diagonal:
\begin{equation}\label{eq:cosine-input}
 A_n:=xI_n+K_n,\qquad D_n:=\det A_n,\qquad D_0:=1.
\end{equation}
Thus $(A_n)_{ii}=x+\cos(2Ci)$; prescribing $(A_n)_{ii}=x$ instead
defines a different family, considered below.

Write
\[
 C_k(i,j):=\det A[i:i+k-1,j:j+k-1],\qquad s:=j-i.
\]
These formulas refer to the infinite kernel; only intervals contained in
the selected finite section are retained.  At condensation depth $q$,
the first layer consists of $C_{q+1}(i,j)$ and the second of
$C_q(i+1,j+1)$.  As for the other ordinary descendants, identities at
vanishing intermediate denominators mean their connected-minor
specializations.

\begin{proposition}\label{prop:cosine-connected-minors}
The first three nontrivial connected-minor layers are
\begin{equation}\label{eq:cosine-C2}
 C_2(i,j)=
 \begin{cases}
 x^2+2cxu_{i+j+1}-\sigma^2,&s=0,\\
 -xu_{i+j+1}-\sigma^2,&|s|=1,\\
 -\sigma^2,&|s|\ge2,
 \end{cases}
\end{equation}
\begin{equation}\label{eq:cosine-C3}
 C_3(i,j)=
 \begin{cases}
 x^3+(4c^2-1)x^2u_{i+j+2}
       -2(1+2c^2)\sigma^2x,&s=0,\\
 x^2u_{i+j+2}+4c\sigma^2x,&|s|=1,\\
 -\sigma^2x,&|s|=2,\\
 0,&|s|>2,
 \end{cases}
\end{equation}
and
\begin{equation}\label{eq:cosine-C4}
 C_4(i,j)=
 \begin{cases}
 x^4+4c(2c^2-1)x^3u_{i+j+3}
       -4(4c^4+1)\sigma^2x^2,&s=0,\\
 -x^3u_{i+j+3}-2(6c^2-1)\sigma^2x^2,&|s|=1,\\
 -\sigma^2x^2,&|s|=2,\\
 0,&|s|>2.
 \end{cases}
\end{equation}
In particular two condensation steps produce a pentadiagonal first
layer, and three produce a pentadiagonal two-layer block.  For
$x\sigma\ne0$ the displayed outer diagonals survive.  The denominator
at depth two still has the nonzero far tail $-\sigma^2$, so it has not
yet acquired a fixed pentadiagonal support.
\end{proposition}

\begin{proof}
For increasing row indices $a<b$ and column indices $d<e$, the kernel
has the minor identity
\begin{equation}\label{eq:cosine-kernel-minor}
 \det\begin{bmatrix}u_{a+d}&u_{a+e}\\u_{b+d}&u_{b+e}\end{bmatrix}
 =-\sin(C(b-a))\sin(C(e-d)).
\end{equation}
All its minors of order at least three vanish.  Expand the indicated
connected minors in their scalar main diagonal entries $x$.  At displacement
two, the only surviving term selects every available diagonal entry,
leaving a contiguous kernel minor of order two; this gives
$-\sigma^2x$ and $-\sigma^2x^2$.  At displacement greater than two,
too few main diagonal entries are available to leave a kernel minor of
order at most two.  The remaining terms in
\eqref{eq:cosine-C2}--\eqref{eq:cosine-C4} follow from
\eqref{eq:cosine-kernel-minor} and
$\sin(2C)=2c\sigma$, $\sin(3C)=(4c^2-1)\sigma$.
For example the coefficient of $x^2$ in the principal minor $C_4(i,i)$ is
\[
 -3\sin^2 C-2\sin^2(2C)-\sin^2(3C)
 =-4(4c^4+1)\sigma^2.
\]
\end{proof}

The ordinary ancestor is equally explicit.  Let $\mathsf S_-^{(n-1)}$
be the lower shift of order $n-1$ and set
\begin{equation}\label{eq:cosine-anchored-parent}
 S_n:=\diag\!\left(1,I-2c\mathsf S_-^{(n-1)}
                         +(\mathsf S_-^{(n-1)})^2\right),
 \qquad M_n:=S_nA_nS_n^T.
\end{equation}
The zero-dimensional tail block is permitted when $n=1$.

\begin{theorem}\label{thm:cosine-ordinary-parent}
The matrices $M_n$ are symmetric pentadiagonal leading sections of one fixed array.  Their
Toeplitz bulk has coefficients
\begin{equation}\label{eq:cosine-ordinary-bulk}
 (m_0,m_1,m_2)=x(2+4c^2,-4c,1),\qquad m_{-j}=m_j.
\end{equation}
Every entry outside the northwest $3\times3$ block agrees with that
bulk.  For $n\ge3$ the complete exceptional block of the parent is
\begin{equation}\label{eq:cosine-ordinary-corner}
 M_n[1:3,1:3]=
 \begin{bmatrix}
 x+u_2&u_3&-u_2\\
 u_3&x+u_4&-2cx-u_3\\
 -u_2&-2cx-u_3&(1+4c^2)x+u_2
 \end{bmatrix}.
\end{equation}
For $n\ge2$, both layers of the order-two descendant are preserved:
\begin{equation}\label{eq:cosine-descendant-equality}
 \operatorname{Cond}^{n-2}(M_n,J)
 =\operatorname{Cond}^{n-2}(A_n,J).
\end{equation}
Consequently $M_n$ is an ordinary ancestor of that complete descendant,
of vaporisation depth $n-2$, and $\det M_n=D_n$.
\end{theorem}

\begin{proof}
The kernel satisfies
$u_t-2cu_{t-1}+u_{t-2}=0$.  The first three rows of $S_nK_n$ are
therefore the first, second, and negative first rows of $K_n$;
all later rows vanish.  Applying the same calculation to the columns
shows that $S_nK_nS_n^T$ is supported in the northwest $3\times3$
corner, where it equals
\[
 \begin{bmatrix}
 u_2&u_3&-u_2\\u_3&u_4&-u_3\\-u_2&-u_3&u_2
 \end{bmatrix}.
\]
The other summand is $xS_nS_n^T$.  Multiplying its fixed three-term
row expressions gives \eqref{eq:cosine-ordinary-bulk} and the remaining
entries of \eqref{eq:cosine-ordinary-corner}.  The formulas use only
preceding indices and are independent of the total size.

Finally, $S_n=\diag(1,S_*)$ with $S_*$ unit lower triangular.
Proposition~\ref{prop:anchored-minor-preservation} with $h=1$ preserves
all four boundary minors and the central denominator in
\eqref{eq:cosine-descendant-equality}.
\end{proof}

For $n\ge5$, three condensations give the displayed pentadiagonal
block; another $n-5$ give the order-two block realized by $M_n$ at
vaporisation depth $n-2$.

\begin{theorem}\label{thm:cosine-recurrence}
For generic $x$ and $c$, the determinant sequence in
\eqref{eq:cosine-input} has minimal eventual constant-coefficient
recurrence order five.  With $\mathcal E$ denoting forward shift, an annihilator is
\begin{equation}\label{eq:cosine-annihilator}
 P(\mathcal E),\qquad
 P(t):=(t-x)^3\bigl(t^2-2x\cos(2C)t+x^2\bigr).
\end{equation}
The ordinary-parent proof gives $P(\mathcal E)D_n=0$ for $n\ge5$.
Equivalently, for $N\ge10$,
\begin{align}\label{eq:cosine-backward-recurrence}
 D_N={}&x(1+4c^2)D_{N-1}
       -2x^2(6c^2-1)D_{N-2}
       +2x^3(6c^2-1)D_{N-3}\notag\\
      &-x^4(1+4c^2)D_{N-4}+x^5D_{N-5}.
\end{align}
These annihilator identities extend to every specialization of $x,c$.
The direct evaluation in Appendix~\ref{sec:cosine-comparison} gives
\eqref{eq:cosine-backward-recurrence} already for $N\ge5$.
\end{theorem}

\begin{proof}
Let $m(z):=\sum_{j=-2}^2m_jz^j$ be the bulk symbol.  Its polynomial is
\begin{equation}\label{eq:cosine-bulk-polynomial}
 z^2m(z)=x(z^2-2cz+1)^2.
\end{equation}
Write $\zeta:=e^{\mathrm i C}$, so $c=(\zeta+\zeta^{-1})/2$.  For generic
parameters the order-four bulk companion $T$ has one Jordan block
of order two at each of $\zeta$ and $\zeta^{-1}$.  By
Lemma~\ref{lem:exact-boundary-companion}, the ordinary determinant
transfer is $x\Lambda^2T$, after the fixed corner.  Its mixed
summand is the tensor product of the two Jordan blocks and has largest
Jordan block size three at eigenvalue $x$.  The other two summands are
one-dimensional, with eigenvalues $x\zeta^2$ and $x\zeta^{-2}$.
This proves \eqref{eq:cosine-annihilator}.  The corner has size three
and the semibandwidth is two, giving the stated threshold $n\ge5$.

For a small exact minimality certificate, specialize $c=1/2$ and $x=1$.
Compute the principal determinants from the ordinary matrices
\eqref{eq:cosine-ordinary-bulk}--\eqref{eq:cosine-ordinary-corner}.
For $H:=(D_{5+i+j})_{0\le i,j\le4}$, this gives
\begin{equation}\label{eq:cosine-Hankel-certificate}
 4H=
 \begin{bmatrix}
 -24&-32&-46&-63&-77\\
 -32&-46&-63&-77&-97\\
 -46&-63&-77&-97&-120\\
 -63&-77&-97&-120&-140\\
 -77&-97&-120&-140&-166
 \end{bmatrix},\qquad
 \det H=-\frac{1701}{128}\ne0.
\end{equation}
This Hankel submatrix uses only indices after the constant transfer begins, and $P(0)=-1$
at this specialization.  Any smaller eventual constant-coefficient
recurrence would propagate backwards along the invertible recurrence
$P(\mathcal E)$ to these indices, contradicting \eqref{eq:cosine-Hankel-certificate}.
The Hankel determinant is a polynomial in $x,c$, so the lower bound
is generic.  Finally, each fixed-index recurrence identity is
polynomial in $x,c$; the proved generic identity therefore specializes
to all parameter values.
\end{proof}

The separate calculation in Section~\ref{sec:cosine-comparison}
gives a closed determinant formula and identifies the following
specializations and variants.

\begin{remark}\label{rem:cosine-degenerate-parameters}
If $x=0$, the original rank-two kernel has
$D_1=\cos(2C)$, $D_2=-\sin^2 C$, and $D_n=0$ for $n\ge3$.
If $\sin C=0$, the kernel has rank one and
\[
 D_n=x^{n-1}(x+n)\qquad(n\ge1).
\]
For $x\ne0$ its minimal constant-coefficient recurrence has order two,
with polynomial $(t-x)^2$.  The first and full tridiagonal layers then
appear after one and two condensation steps, respectively; a
first-order anchored filter $\diag(1,I-c\mathsf S_-)$ gives the
corresponding ordinary tridiagonal ancestor with a northwest
$2\times2$ correction.

The remaining coalescence of eigenvalues for real $C$ occurs at $\cos C=0$.  For $x\ne0$,
\[
 D_n=x^{n-2}\left[
 -\frac{n^2}{4}+x^2-\frac{x}{2}+\frac18
 +\frac{4x-1}{8}(-1)^n\right].
\]
For $x\ne0,1/4$ the minimal order is four, with polynomial
$(t-x)^3(t+x)$; at $x=1/4$ it is three, with polynomial $(t-x)^3$.
The reduction comes from a vanishing boundary coefficient or coincident
eigenvalues, consistently with the generic order-five statement.
\end{remark}

\paragraph{Prescribing the diagonal instead.}
If the intended matrix has $\widehat A_{ii}=x$ and
$\widehat A_{ij}=u_{i+j}$ for $i\ne j$, put
\[
 d_i:=x-u_{2i},\qquad
 \widehat A=\diag(d_i)+K.
\]
Its upper and lower separated ranks are still two.  The same
condensation support statement and the same anchored filter apply,
but the ordinary pentadiagonal bulk is now position-dependent.  Put
$\widehat M_n:=S_n\widehat A_nS_n^T$.
Outside the northwest $3\times3$ corner its entries are
\begin{align}\label{eq:cosine-prescribed-diagonal-bands}
 \widehat m_0[i]&=d_i+4c^2d_{i-1}+d_{i-2}&& (i\ge4),\notag\\
 \widehat m_1[i]&=-2c(d_i+d_{i-1})&& (i\ge3),\\
 \widehat m_2[i]&=d_i&& (i\ge2),\notag
\end{align}
with $\widehat m_{-s}[i]=\widehat m_s[i]$.  Its complete exceptional northwest block is
\[
 \widehat M[1:3,1:3]=
 \begin{bmatrix}
 x&u_3&-u_2\\
 u_3&x&-2cd_2-u_3\\
 -u_2&-2cd_2-u_3&d_3+4c^2d_2+u_2
 \end{bmatrix}.
\]
The exact equality of both layers of the order-two descendant still holds.  The
appropriate recurrence construction is the position-dependent
six-state ordinary row--column transfer and its observation-row
dependences from Section~\ref{sec:background} and
\eqref{eq:local-vapor-observation-rows}; the constant order-five
polynomial \eqref{eq:cosine-annihilator} is not asserted for this
different family.

More generally, $\cos(\phi(i)+\psi(j))$ is a separated rank-two
kernel for arbitrary phase functions.  When the indicated
sines do not vanish, its local row annihilator for $i\ge4$ is
\[
 \alpha_{i,1}:=-\frac{\sin(\phi(i)-\phi(i-2))}
                       {\sin(\phi(i-1)-\phi(i-2))},\qquad
 \alpha_{i,2}:=\frac{\sin(\phi(i)-\phi(i-1))}
                      {\sin(\phi(i-1)-\phi(i-2))}.
\]
The column coefficients have the analogous formulas with $\psi,j$.
Thus the general local anchored theorem supplies a pentadiagonal
ordinary ancestor and a variable six-state transfer.  Constant
annihilator coefficients in the size variable require additional
structure; separated rank alone does not supply them.

\subsection{Polynomial kernels}

\subsubsection{Absolute powers of the distance}

Work in characteristic zero, let $p\ge1$ be an integer, and put
\[
 A_n(x,p):=T_n(b),\qquad
 b_0:=x,\qquad b_s:=|s|^p\quad(s\ne0),\qquad r:=p+1.
\]
For this Toeplitz family, write
\begin{equation}\label{eq:polyexplicit-minors}
 C_t(s):=\det\bigl(b_{s+j-i}\bigr)_{0\le i,j<t}.
\end{equation}
These are the connected minors of
\eqref{eq:application-connected-minors}, with the starting position
suppressed.  Symmetry gives $C_t(-s)=C_t(s)$.

\begin{theorem}\label{thm:polyexplicit-compression}
After $r$ ordinary condensation steps, the first layer of $A_n(x,p)$ is
symmetric Toeplitz and $(r,r)$-banded.  After $r+1$ steps, both layers are
symmetric Toeplitz and $(r,r)$-banded.  More precisely, put
\begin{equation}\label{eq:polyexplicit-constant}
 \kappa_p:=(p!)^{p+1}.
\end{equation}
Then
\begin{align}
 C_r(s)&=\kappa_p &&(|s|\ge r),
 \label{eq:polyexplicit-fartail}\\
 C_{r+1}(s)&=C_{r+2}(s)=0 &&(|s|>r),
 \label{eq:polyexplicit-support}\\
 C_{r+1}(\pm r)&=(-1)^r\kappa_p x,
 & C_{r+2}(\pm r)&=\kappa_p x^2.
 \label{eq:polyexplicit-edges}
\end{align}
Thus the denominator layer at depth $r$ still has nonzero constant far
tails.  The full banded block at depth $r+1=p+2$ has symbols
$C_{r+2}(s)$ and $C_{r+1}(s)$ for its first and second layers, respectively.
\end{theorem}

\begin{proof}
On each strict triangular part, the binomial expansion of $(j-i)^p$ or
$(i-j)^p$ gives a sum of $r$ separated kernels.  The upper and lower
quasiseparable ranks are therefore at most $r$.
Theorem~\ref{thm:structural-band-compression} gives
\eqref{eq:polyexplicit-support} and the stated depths.

For $s\ge r$, every entry in the size-$r$ minor lies strictly above the
original diagonal.  Factor its entries as
\[
 (s+j-i)^p=\sum_{a=0}^p(-i)^a\binom pa(s+j)^{p-a}.
\]
The two evaluation determinants are Vandermonde determinants with the same
sign $(-1)^{p(p+1)/2}$.  Their signs cancel, and the remaining factors give
\[
 C_r(s)=
 \left(\prod_{j=1}^{p}j!\right)^2
 \prod_{a=0}^{p}\binom pa
 =(p!)^{p+1}.
\]
Transpose symmetry proves \eqref{eq:polyexplicit-fartail} on the negative
side as well.

The matrix defining $C_{r+1}(r)$ differs from the rank-at-most-$r$
polynomial matrix $((r+j-i)^p)$ only at its bottom-left entry, where $0$ is
replaced by $x$.  Expanding in this correction gives
\[
 C_{r+1}(r)=(-1)^r x C_r(r+1)=(-1)^r\kappa_p x.
\]
Desnanot--Jacobi at shift $r$ now yields
\[
 C_{r+2}(r)C_r(r)
 =C_{r+1}(r)^2-C_{r+1}(r-1)C_{r+1}(r+1)
 =\kappa_p^2x^2.
\]
Division by the nonzero constant $\kappa_p$ proves the second outer-edge
formula; symmetry gives its negative counterpart.
\end{proof}

For $p=1$, the complete block at depth three has the following symmetric
symbols, with all coefficients outside $[-2,2]$ equal to zero:
\begin{align}
 (C_3(0),C_3(1),C_3(2))
 & =\bigl((x-2)(x^2+2x-2),\ 3x^2-4x+2,\ x\bigr),
 \label{eq:polyexplicit-poneB}\\
 (C_4(0),C_4(1),C_4(2))
 & =\bigl((x^2-4x+2)(x^2+4x-6),
 \nonumber\\[-2pt]
 &\hspace{20mm}-2(2x-1)(x^2-2x+2),\ x^2\bigr).
 \label{eq:polyexplicit-poneA}
\end{align}
These formulas follow by expanding the fixed minors
\eqref{eq:polyexplicit-minors}.

The symmetric ordinary parent below supplies the recurrence directly.
We first verify that the subsequent in-band condensations can be
performed over $\F(x)$.  For $0\le s\le r$ and $t\ge r+1$, the coefficient of
$x^{t-s}$ in $C_t(s)$ is
\begin{equation}\label{eq:polyexplicit-leading-coefficient}
 \coef{x^{t-s}}C_t(s)=(-1)^{s(t-s)}
 \det\bigl((t+v-u)^p\bigr)_{0\le u,v<s},
\end{equation}
with the empty determinant equal to $1$.  The displayed determinant is a
positive integer when $s\ge1$.  To see this, reverse its rows and write
its entries as $(\xi_i+\eta_j)^p$, where
$\xi_i:=t-s+1+i>0$ and $\eta_j:=j$.
Expand using the terms $\binom pa\xi_i^a\eta_j^{p-a}$.  Cauchy--Binet
expresses the determinant as a sum of products of generalized Vandermonde
determinants.  Reversing the selected powers
$p-a$ contributes the same sign as the row reversal, so every resulting
summand in the original determinant is nonnegative.  The choice of powers
$a=p-s+1,\ldots,p$ is strictly positive.  Thus every in-band connected minor
needed by subsequent condensation is nonzero in $\F(x)$.

The original determinants are monic polynomials in $x$.  The recurrence
bounds proved below therefore extend by polynomial Hankel identities
to every specialized value of $x$, including vanishing intermediate minors.

\subsubsection*{An explicit ordinary parent for polynomial-distance matrices}

Put $r:=p+1$ and $\alpha_a:=(-1)^a\binom ra$ for $0\le a\le r$.
The first coordinate is fixed in the following ordinary parent.  Define
\begin{equation}\label{eq:poly-anchored-transform}
 S_n:=\diag\bigl(1,(I-\mathsf S_-^{(n-1)})^r\bigr),\qquad
 M_n:=S_n A_n(x,p)S_n^T.
\end{equation}
The zero-dimensional block is allowed when $n=1$.

\begin{theorem}\label{thm:poly-anchored-vapor}
For every $n\ge3$,
\begin{equation}\label{eq:poly-anchored-descendant}
 \operatorname{Cond}^{n-2}(A_n(x,p),J)
 =\operatorname{Cond}^{n-2}(M_n,J).
\end{equation}
The ordinary matrices $M_n$ are symmetric $(r,r)$-banded leading
sections of one fixed array.  Their Toeplitz bulk is
\begin{equation}\label{eq:poly-anchored-bulk}
 c_s:=\sum_{a=-r}^r(-1)^a\binom{2r}{r+a}b_{s+a},\qquad
 b_0=x,\quad b_s=|s|^p\ (s\ne0),
\end{equation}
with $c_s=0$ for $|s|>r$ and $c_r=(-1)^r x$.  The correction
$E_n:=M_n-T_n(c)$ is supported in the fixed northwest
$(r+1)\times(r+1)$ block.  Every entry of the constructed parent is
given explicitly by
\begin{equation}\label{eq:poly-anchored-entries}
\begin{aligned}
 (M_n)_{11}&=x,\\
 (M_n)_{1j}&=\sum_{a=0}^{\min(r,j-2)}\alpha_a b_{j-1-a}
       &&(j\ge2),\\
 (M_n)_{ij}&=
 \sum_{a=0}^{\min(r,i-2)}\sum_{b=0}^{\min(r,j-2)}
        \alpha_a\alpha_b b_{j-i+a-b}
       &&(i,j\ge2),
\end{aligned}
\end{equation}
and symmetry.  Hence $\det M_n=\det A_n(x,p)$ by exact ordinary
vaporisation, with the order-two block retained.
\end{theorem}

\begin{proof}
The matrix $S_n$ in \eqref{eq:poly-anchored-transform} is
unit lower triangular with an isolated first coordinate.
Proposition~\ref{prop:anchored-minor-preservation}, with $h=1$,
therefore gives \eqref{eq:poly-anchored-descendant}, including its
central denominator entry.

The entry formulas follow by expanding the displayed row operations.
Whenever $i,j\ge r+2$, both sums have their full ranges; collecting terms with
the same difference gives \eqref{eq:poly-anchored-bulk}.  If $|s|>r$, all
arguments $s+a$ lie strictly on the same side of zero, so a difference of
order $2r$ annihilates the polynomial of degree $p=r-1$.  At $s=r$, only
the replacement of the polynomial value $0$ at the origin by $x$ survives,
giving $c_r=(-1)^r x$.

To locate the northwest correction exactly, suppose $j\ge r+2$.  If $i=1$, the complete
column difference in \eqref{eq:poly-anchored-entries} acts on positive
arguments and gives zero, as does $c_{j-1}$.  If $i\ge2$, compare the
truncated row sum with the complete bulk sum.  Every omitted row index has
$a\ge i-1$.  For each omitted term the full column difference involves
arguments at least
\[
 j-r-i+a\ge j-r-1\ge1,
\]
and therefore vanishes.  Thus $(M_n)_{ij}=c_{j-i}$ whenever $j\ge r+2$;
transpose symmetry proves the same when $i\ge r+2$.  This gives the
$(r+1)\times(r+1)$ corner bound.  Inside that corner the largest possible
distance from the diagonal is $r$, proving that the whole matrix is $(r,r)$-banded.
All formulas are polynomial in $x$, so no nonzero-minor qualification is
needed for the final identities.
\end{proof}

The construction preserves an actual two-layer descendant.  For $n\ge r+3$
one can first perform $r+1$ condensation steps to reach the explicit banded
block, continue for $n-r-3$ steps to order two, and vaporise that block for
$n-2$ steps into $M_n$.  If retaining a larger final block is useful, replace
$S_n$ by
\[
 S_{n,h}:=\diag\bigl(I_h,(I-\mathsf S_-^{(n-h)})^r\bigr).
\]
Proposition~\ref{prop:anchored-minor-preservation} then gives equality of
both layers of the order-$(h+1)$ descendant.  The same calculation bounds the
corner by $(h+r)\times(h+r)$; the Toeplitz bulk remains
\eqref{eq:poly-anchored-bulk}.  These larger-block statements do not assert
that the parent preserves the entire growing block at the first banded
depth.

For $p=1$, the ordinary parent has symmetric pentadiagonal bulk
\[
 (c_0,c_1,c_2)=(6x-4,2-4x,x),
\]
and, for $n\ge3$, its complete exceptional northwest block is
\begin{equation}\label{eq:poly-anchored-pone-corner}
 M_n[1:3,1:3]=
 \begin{bmatrix}
 x&1&0\\
 1&x&1-2x\\
 0&1-2x&5x-4
 \end{bmatrix}.
\end{equation}
Equivalently,
\[
 E_n[1:3,1:3]=
 \begin{bmatrix}
 4-5x&4x-1&-x\\
 4x-1&4-5x&2x-1\\
 -x&2x-1&-x
 \end{bmatrix}.
\]
All entries outside the displayed corner agree with the Toeplitz bulk.

For $p=3$, the symmetric bulk has lower and upper semibandwidth four and
\[
 (c_0,c_1,c_2,c_3,c_4)
 =(70x+32,-56x-18,28x,2-8x,x).
\]
For $n\ge5$, the exceptional block is
\begin{equation}\label{eq:poly-anchored-pthree-corner}
\begin{bmatrix}
 x&1&4&1&0\\
 1&x&1-4x&6x+4&1-4x\\
 4&1-4x&17x-8&-28x-9&22x-4\\
 1&6x+4&-28x-9&53x+40&-52x-19\\
 0&1-4x&22x-4&-52x-19&69x+32
\end{bmatrix}.
\end{equation}
The general formula \eqref{eq:poly-anchored-entries} supplies a compact
description at every higher odd power.

\begin{corollary}\label{cor:polydistance}
For odd $p\ge1$, the determinant sequence of $A_n(x,p)$ is eventually
annihilated by a homogeneous constant-coefficient recurrence of order at most
\begin{equation}\label{eq:polydistancebound}
 \frac{3^{p+1}+1}{2}.
\end{equation}
\end{corollary}

\begin{proof}
The tails satisfy the $r$th-difference recurrences with
$m_-=m_+=r$ and $\alpha_a=\beta_a=(-1)^a\binom ra$.
Theorem~\ref{thm:poly-anchored-vapor} gives their explicit symmetric
parent and fixed corner.  Applying the symmetric bound in
Corollary~\ref{cor:anchored-tail-recurrence} proves
\eqref{eq:polydistancebound}, including repeated-root specializations.
\end{proof}

\subsubsection*{The polynomial-tail specialization of the ordinary transfer}

The symmetric reflection bound is uniform over all symmetric band symbols.
The particular bulk in \eqref{eq:poly-anchored-bulk} has additional structure,
which permits a smaller bound beginning at $p=5$.

Lemma~\ref{lem:exact-boundary-companion} applies directly at repeated
roots.  For the anchored parent its corner has size $r+1$, and its leading
determinant is monic in $x$, so that leading block is invertible
over $\mathbb Q(x)$.  Write $c(z):=\sum_{s=-r}^r c_s z^s$.

For odd $p=2m-1$, put $r=2m$.  Define the Eulerian polynomial $A_p(z)$ by
\begin{equation}\label{eq:poly-eulerian-definition}
 \sum_{j\ge1}j^p z^j=\frac{zA_p(z)}{(1-z)^{p+1}}.
\end{equation}
It has degree $p-1$, is reciprocal, and $A_p(1)=p!$.  The ordinary bulk
symbol in \eqref{eq:poly-anchored-bulk} satisfies
\begin{equation}\label{eq:poly-bulk-unit-root-factor}
 z^r c(z)=(1-z)^r Q_x(z),\qquad
 Q_x(z):=x(1-z)^r+2zA_{r-1}(z).
\end{equation}
One way to verify this finite Laurent-polynomial identity is to multiply
the positive and negative tail generating functions by
$(1-z)^r(1-z^{-1})^r$ and use
\eqref{eq:poly-eulerian-definition}.  Since $r$ is even, the two rational
tail contributions agree.  The diagonal correction contributes
$xz^{-r}(1-z)^{2r}$.  This gives
\eqref{eq:poly-bulk-unit-root-factor} directly.

\begin{theorem}\label{thm:polydistance-improved-bound}
For $p=2m-1\ge1$, the sequence $D_n(x):=\det A_n(x,p)$ has an eventual
constant-coefficient recurrence over $\mathbb Q(x)$ of order at most
\begin{equation}\label{eq:polydistance-improved-degree}
 d_m:=\frac{(2m^2+1)3^m+1}{2}.
\end{equation}
The same order bound holds after specialization of $x$ to any
characteristic-zero value.  For $p=1,3,5$, the generic minimal orders are
$5,41,257$, respectively.
\end{theorem}

\begin{proof}
The polynomials $(1-z)^r$ and $2zA_{r-1}(z)$ are coprime, since
$A_{r-1}(1)=(r-1)!\ne0$.  Their nonconstant pencil $Q_x$ is squarefree for
generic $x$ in characteristic zero.  It is reciprocal, has nonzero
constant term for $x\ne0$, and $Q_x(1)=2(r-1)!$.  Its generic roots can
therefore be written
\[
 \rho_1,\rho_1^{-1},\ldots,\rho_m,\rho_m^{-1}.
\]
By \eqref{eq:poly-bulk-unit-root-factor}, the bulk companion is similar over
a splitting field to
\[
 J_{2m}(1)\oplus
 \diag(\rho_1,\rho_1^{-1},\ldots,\rho_m,\rho_m^{-1}).
\]
Let $T$ denote this bulk companion.  The ordinary determinant transfer in
\eqref{eq:exact-boundary-shooting} is $x\Lambda^{2m}T$.

Label a product of the nonunit roots by
$\epsilon\in\{-1,0,1\}^m$, write $|\epsilon|$ for its number of
nonzero coordinates, and put $j:=|\epsilon|$.    Such a product uses one root from each of those $j$ pairs.
From any of the remaining $m-j$ pairs it may use both roots or neither.
If it uses $f$ complete pairs, the number of selected unit-root vectors is
\[
 k=2m-j-2f,\qquad 0\le f\le m-j.
\]
Lemma~\ref{lem:exterior-jordan-length} bounds the largest Jordan block size by
$k(2m-k)+1$.  Maximizing over these choices gives
\begin{equation}\label{eq:poly-mode-exponent}
 \nu_j:=\begin{cases}
 m^2+1,&j\equiv m\pmod2,\\
 m^2,&j\not\equiv m\pmod2.
 \end{cases}
\end{equation}
Consequently the following monic polynomial annihilates the ordinary
determinant transfer:
\begin{equation}\label{eq:poly-mode-annihilator}
 \mathcal A_{p,x}(z):=
 \prod_{\epsilon\in\{-1,0,1\}^m}
 \left(z-x\prod_{a=1}^{m}\rho_a^{\epsilon_a}\right)^\nu_{|\epsilon|}.
\end{equation}
Its coefficients are fixed under permutations and inversions of the
reciprocal pairs and thus lie in $\mathbb Q(x)$.  More explicitly, the
invariance makes them symmetric polynomials in
$\rho_a+\rho_a^{-1}$, so their only possible parameter denominators are
powers of $x$.  Its constant term is $(-1)^{d_m}x^{d_m}=-x^{d_m}$,
because the nonunit root factors cancel in pairs.

There are $\binom mj2^j$ labels with $j$ nonzero coordinates.  The sum of
the exponents is
\begin{align*}
 \deg\mathcal A_{p,x}
 &=m^2\sum_{j=0}^m\binom mj2^j
   +\sum_{\substack{0\le j\le m\\j\equiv m\ ({\rm mod}\ 2)}}
                       \binom mj2^j\\
 &=m^2 3^m+\frac{3^m+1}{2}=d_m.
\end{align*}
Equation~\eqref{eq:exact-boundary-shooting} transfers this
annihilator to $D_n$ beyond a fixed initial boundary.  Because every
$D_n(x)$ is a polynomial, the finite Hankel-minor identities giving this
order bound are polynomial identities in $x$; they extend to every
specialization, including $x=0$.

For the three stated minimality claims, an exact certificate is recorded
below.  At $x=17$, form
\[
 H_{p}:=(D_{a+i+j}(17))_{0\le i,j<d_m},\qquad a:=4r+10.
\]
The resulting integer determinants have the following nonzero reductions
modulo $1000000007$:
\[
\begin{array}{c|r|r|r}
 p&d_m&a&\det H_p\pmod{1000000007}\\ \hline
 1&5&18&945293776\\
 3&41&26&503118073\\
 5&257&34&504544247
\end{array}
\]
The certificates are computed from integer matrices, with complete modular
sequences independently checked by an ordinary banded elimination and a
symmetric Toeplitz Schur-complement recursion on the original symbol.
Dense exact original determinants check the initial orders $1$ through
$18$ as well.  The scripts, modular sequences, and elimination pivots
are included with the manuscript.

These witnesses occur after the constant bulk transfer begins:
\eqref{eq:exact-boundary-shooting} applies for $n\ge2r+1$.
The annihilator has nonzero constant coefficient at $x=17$.  Therefore
any eventual lower-order recurrence would propagate backwards along this
invertible recurrence to the displayed Hankel submatrix, contradicting its
nonzero determinant.  Nonvanishing at $x=17$ also makes the corresponding
Hankel determinant a nonzero polynomial in $x$.  Thus each lower bound is
generic in characteristic zero.  Combining it with the proved upper bound
gives minimal orders $5,41,257$.
\end{proof}

The uniform symmetric-band bound $(3^{p+1}+1)/2$ is attained for
$p=1,3$.  For $p=5$ the generic order is $257$, below the uniform bound
$365$.  Formula~\eqref{eq:polydistance-improved-degree} remains an upper
bound for general odd $p$; equality at every odd power is open.

For $p=1$, the annihilator is particularly simple:
\begin{equation}\label{eq:distance-linear-quintic}
 \mathcal A_{1,x}(z)
 =(z-x)\bigl(z^2-2(x-1)z+x^2\bigr)^2.
\end{equation}
Its expanded recurrence is
\begin{align}\label{eq:distance-linear-orderfive}
 D_n={}&(5x-4)D_{n-1}-(10x^2-12x+4)D_{n-2}\notag\\
 &+x(10x^2-12x+4)D_{n-3}
  -x^3(5x-4)D_{n-4}+x^5D_{n-5},\qquad n\ge5.
\end{align}
The initial values are
\[
 D_0=1,\quad D_1=x,\quad D_2=x^2-1,\quad
 D_3=x^3-6x+4,\quad D_4=x^4-20x^2+32x-12.
\]
The boundary formula starts at $n=5$, so the annihilator proves the
recurrence for $n\ge10$.  Direct exact polynomial determinants verify
the remaining five cases $n=5,\ldots,9$; the accompanying symbolic
verification script records these computations.
The quintic is generically minimal by the certificate above; multiplying
it by $z-x$ gives a valid but nonminimal degree-six annihilator.

\begin{example}\label{ex:cubic-distance-vaporisation}
Consider the cubic-distance family
\[
 A_n:=\bigl(b_{j-i}\bigr)_{1\le i,j\le n},\qquad
 b_0:=x,\quad b_s:=|s|^3\ (s\ne0),\qquad D_n:=\det A_n.
\]
This example makes each stage of the construction explicit.  In the
notation \eqref{eq:polyexplicit-minors}, the separated rank of either
strict triangular tail is four, and
\[
 C_4(s)=1296\quad(|s|\ge4),\qquad
 C_5(s)=C_6(s)=0\quad(|s|>4).
\]
The surviving outer entries are
\[
 C_5(\pm4)=1296x,\qquad C_6(\pm4)=1296x^2.
\]
Thus four condensation steps produce a first layer of semibandwidth four,
while its denominator layer still has the nonzero constant tails
$C_4(s)$.  Five steps produce a full banded two-layer block, whose symbols
are $C_6$ and $C_5$.  These values follow from Theorem~\ref{thm:polyexplicit-compression}.

The ordinary parent is
\[
 S_n:=\diag\bigl(1,(I-\mathsf S_-^{(n-1)})^4\bigr),\qquad
 M_n:=S_nA_nS_n^T.
\]
Theorem~\ref{thm:poly-anchored-vapor} identifies its complete order-two
descendant with that of $A_n$.  For $n\ge7$, the route is therefore five
condensation steps to the banded block, another $n-7$ steps to order two,
and $n-2$ vaporisation steps to the ordinary matrix $M_n$.  The parent
has the nine-diagonal Toeplitz bulk
\[
 (c_0,c_1,c_2,c_3,c_4)
 =(70x+32,-56x-18,28x,2-8x,x)
\]
and the complete northwest $5\times5$ block in
\eqref{eq:poly-anchored-pthree-corner}.  These formulas specify every
entry of the ordinary parent, its retained layers, and the construction depth.

The ordinary row--column state has dimension $\binom84=70$.  Its smaller
annihilator follows from the factorization of its bulk polynomial:
\[
 z^4c(z)=(1-z)^4Q_x(z),\qquad
 Q_x(z):=x(1-z)^4+2z(1+4z+z^2).
\]
For generic $x$, let the roots of $Q_x$ be
$\rho_1,\rho_1^{-1},\rho_2,\rho_2^{-1}$.  By
Lemma~\ref{lem:exact-boundary-companion}, the determinant transfer is the
fourth compound of the bulk companion, multiplied by $x$.  The companion
has one block $J_4(1)$ and these four simple nonunit roots.
Lemma~\ref{lem:exterior-jordan-length} shows that the largest Jordan block at eigenvalue $x$ has
size at most five; those at $x\rho_i^{\pm1}$ have size at most four;
and those at $x\rho_1^{\pm1}\rho_2^{\pm1}$ have size at most five.

These three groups give an explicit degree-$41$ annihilator without
expanding forty-one coefficients.  Define
\begin{align}
 P_x(z)&:=z^4+(2-4x)z^3+(6x^2+8x)z^2
                  +(2x^2-4x^3)z+x^4,
 \label{eq:cubic-mode-quartic-single}\\
 R_x(z)&:=z^4-(4x+8)z^3+(6x^2-32x+4)z^2
                  -(4x^3+8x^2)z+x^4.
 \label{eq:cubic-mode-quartic-cross}
\end{align}
Then
\begin{equation}\label{eq:cubic-mode-annihilator}
 \mathcal A_{3,x}(z)=(z-x)^5P_x(z)^4R_x(z)^5.
\end{equation}
To verify the two quartics, put $y_i:=\rho_i+\rho_i^{-1}$.  Dividing
$Q_x(z)$ by $z^2$ gives
\[
 x(y-2)^2+2y+8=0,
 \qquad y_1+y_2=4-2/x,\quad y_1y_2=4+8/x.
\]
The single-root product is $P_x(z)=x^3Q_x(z/x)$.  For the cross products,
write $a:=\rho_1\rho_2$ and $b:=\rho_1/\rho_2$.  Then
\[
 (a+a^{-1})+(b+b^{-1})=y_1y_2,
 \qquad
 (a+a^{-1})(b+b^{-1})=y_1^2+y_2^2-4.
\]
Expanding $(z-xa)(z-x/a)(z-xb)(z-x/b)$ now gives
\eqref{eq:cubic-mode-quartic-cross}.  This proves
\eqref{eq:cubic-mode-annihilator} by the ordinary transfer.

In terms of the forward shift $\mathcal E$, the resulting
recurrence is
\begin{equation}\label{eq:cubic-compact-recurrence}
 (\mathcal E-x)^5P_x(\mathcal E)^4R_x(\mathcal E)^5D_n=0
 \qquad(n\ge9).
\end{equation}
The threshold $9=5+4$ follows from the fixed corner and the boundary
companion identity.  All coefficients in
\eqref{eq:cubic-mode-annihilator} are polynomials in $x$, so this identity
extends to every specialization of $x$.  The certificate in
Theorem~\ref{thm:polydistance-improved-bound} proves that its degree $41$
is generically minimal.  For example, the first backward recurrence
coefficient is $41x+32$ and the last is $x^{41}$.

For comparison, the conventional unanchored finite differences
\[
 \widehat S_n:=(I-\mathsf S_-^{(n)})^4,\qquad
 \widehat M_n:=\widehat S_nA_n\widehat S_n^T
\]
give the same Toeplitz bulk and a northwest $4\times4$ corner, as in
\eqref{eq:polyexplicit-corner}.  Unitriangularity proves the determinant
identity for this second representation directly.  The anchored parent
additionally satisfies the connected-minor equalities that make it an
ordinary ancestor of the specified order-two block.  The two calculations
therefore explain the same scalar recurrence while displaying the
different boundary information retained by the constructions.
\end{example}

\subsubsection{Signed powers}

For the signed polynomial kernel, the anchored ordinary parent has
a bulk companion with a single repeated eigenvalue.  It gives the recurrence by
vaporisation and the ordinary transfer.  A moment determinant then
proves that the scalar sequence requires the full predicted recurrence order.

\begin{theorem}
\label{thm:signedpolydistance}
Let $p\ge1$ be an integer, put $r:=p+1$, and over a field of characteristic
zero define
\[
 A_n^{\rm sgn}(x,p)=(a_{ij}),\qquad
 a_{ii}=x,\qquad a_{ij}=(i-j)^p\quad(i\ne j).
\]
Assume $x\ne0$ and put
\begin{equation}\label{eq:signedpolydegree}
 d:=r^2+1=(p+1)^2+1.
\end{equation}
Then
\begin{equation}\label{eq:signedfactor}
 \det A_n^{\rm sgn}(x,p)=x^{n-r}Q_p(n;x),
\end{equation}
where $Q_p(n;x)$ is a polynomial in $n$ of exact degree $r^2$.
Consequently the determinant sequence satisfies
\begin{equation}\label{eq:signedpoly}
 \sum_{j=0}^{d}(-x)^j\binom dj
 \det A_{n-j}^{\rm sgn}(x,p)=0
\end{equation}
for all sufficiently large $n$, and its minimal homogeneous
constant-coefficient recurrence order is exactly $d$.
\end{theorem}

\begin{proof}
Use the anchored $r$th differences of
Proposition~\ref{prop:anchored-minor-preservation} on the two sides of
the signed Toeplitz symbol $b_s=(-s)^p$ for $s\ne0$, $b_0=x$.
The global polynomial part disappears from the bulk.  The resulting
ordinary parent retains the complete order-two descendant,
has a correction in its northwest $(r+1)\times(r+1)$ block, and has bulk
\[
 c_s=(-1)^s\binom{2r}{r+s}x\quad(|s|\le r),\qquad
 z^r c(z)=(-1)^r x(1-z)^{2r}.
\]
Its exact determinant transfer from
Lemma~\ref{lem:exact-boundary-companion} is $x\Lambda^rJ_{2r}(1)$.
Lemma~\ref{lem:exterior-jordan-length} therefore gives the annihilator
$(\mathcal E-x)^{r^2+1}$.  This proves the
recurrence through the ordinary ancestor.  We next verify that the
polynomial factor has the full asserted degree and prove
\eqref{eq:signedfactor}.

Let $U_n$ be the $n\times r$ polynomial-evaluation matrix
\begin{equation}\label{eq:signedU}
 (U_n)_{i,a}=i^a,
 \qquad 1\le i\le n,\quad 0\le a\le p,
\end{equation}
and let $\mathsf C_p$ be the $r\times r$ anti-diagonal matrix
\begin{equation}\label{eq:signedCp}
 (\mathsf C_p)_{a,b}=
 \begin{cases}
 (-1)^b\binom pa,&a+b=p,\\
 0,&a+b\ne p.
 \end{cases}
\end{equation}
The binomial theorem gives the exact rank-$r$ factorization
\begin{equation}\label{eq:signedlowrank}
 \bigl((i-j)^p\bigr)_{1\le i,j\le n}
 =U_n\mathsf C_pU_n^T.
\end{equation}
Since the kernel in \eqref{eq:signedlowrank} has zero diagonal,
\[
 A_n^{\rm sgn}(x,p)=xI_n+U_n\mathsf C_pU_n^T.
\]
Put
\begin{equation}\label{eq:signedGram}
 G_n:=U_n^TU_n,
 \qquad
 (G_n)_{a,b}=\sum_{i=1}^n i^{a+b}.
\end{equation}
Sylvester's determinant identity gives, in the field with $x$ inverted,
\begin{equation}\label{eq:signedSylvester}
 \det A_n^{\rm sgn}(x,p)
 =x^{n-r}\det(xI_r+\mathsf C_pG_n).
\end{equation}
Thus $Q_p(n;x):=\det(xI_r+\mathsf C_pG_n)$ gives
\eqref{eq:signedfactor}.

By Faulhaber's formula,
\begin{equation}\label{eq:signedFaulhaber}
 (G_n)_{a,b}
 =\frac{n^{a+b+1}}{a+b+1}+\text{lower powers of }n.
\end{equation}
Hence
\begin{equation}\label{eq:signeddetGdegree}
 \deg_n\det G_n=r^2,
\end{equation}
with leading coefficient
\[
 \det\left(\frac1{a+b+1}\right)_{0\le a,b<r},
\]
the nonzero Hilbert determinant.  The anti-diagonal matrix $\mathsf C_p$ is
invertible, so the term
$\det(\mathsf C_pG_n)=\det(\mathsf C_p)\det G_n$ contributes a nonzero
coefficient of degree $r^2$ to $Q_p$.

No other term can cancel it.  Indeed, $\mathsf C_p$ only reverses and rescales
the rows of $G_n$.  A $k\times k$ minor of $G_n$ has degree at most
\[
 \max_{|I|=|J|=k}
 \left(\sum_{a\in I}a+\sum_{b\in J}b+k\right)
 =k(2r-k),
\]
which is strictly less than $r^2$ when $k<r$.  Expanding
$\det(xI_r+\mathsf C_pG_n)$ by the number of entries selected from the second
summand therefore proves
\[
 \deg_nQ_p(n;x)=r^2=d-1.
\]

Finally write $D_n:=\det A_n^{\rm sgn}(x,p)=x^{n-r}Q_p(n;x)$.  Then the
left side of \eqref{eq:signedpoly} equals
\[
 x^{n-r}\sum_{j=0}^{d}(-1)^j\binom djQ_p(n-j;x),
\]
the $d$th backward finite difference of a polynomial of degree $d-1$, hence
zero.  Since that degree is exact, the $(d-1)$st difference is nonzero.
Equivalently, $D_n$ is a nonzero geometric factor $x^n$ times a polynomial
of exact degree $d-1$, whose minimal constant-coefficient annihilator is
$(\mathcal E-x)^d$.  Thus the minimal recurrence order is exactly $d$.
\end{proof}

\begin{remark}
The same square $(p+1)^2$ has two explicit realizations: it is the
maximal polynomial degree in the middle compound of a single
unipotent bulk block, and it is the degree of the moment determinant
\eqref{eq:signedSylvester}.  The ordinary-parent proof supplies the
annihilator, while the moment calculation proves that no shorter
constant-coefficient recurrence is possible.
\end{remark}

\subsection{Finite exponential kernels and their smaller recurrence}

Finite exponential kernels give another family whose ordinary ancestor
has repeated companion roots and a recurrence smaller than the general
binomial state bound.
Let $r\ge1$, let $v_1,\ldots,v_r$ be nonzero, and put
\begin{equation}\label{eq:exponential-kernel}
 H_{ij}:=\begin{cases}
 x,&i=j,\\
 \displaystyle\sum_{s=1}^r u_s v_s^{i-j},&i\ne j,
 \end{cases}
 \qquad \delta:=x-\sum_{s=1}^r u_s,\qquad
 D_n:=\det H_n.
\end{equation}
Write $(\mathcal E D)_n:=D_{n+1}$ for the forward shift.

\begin{theorem}\label{thm:exponential-vapor-order}
In characteristic zero the sequence \eqref{eq:exponential-kernel} is
eventually annihilated by
\begin{equation}\label{eq:exponential-annihilator}
 \mathcal P_r(\mathcal E),\qquad
 \mathcal P_r(z):=
 \prod_{\substack{e\in\{-1,0,1\}^r\\\sum_s e_s=0}}
 \left(z-\delta\prod_{s=1}^r v_s^{e_s}\right)^{1+\#\{s:e_s=0\}}.
\end{equation}
Its degree is
\begin{equation}\label{eq:exponential-order}
 d_r:=\sum_{k=0}^{\lfloor r/2\rfloor}
 (r-2k+1)\frac{r!}{k!^2(r-2k)!}
 =\sum_{j=0}^r\binom rj\frac{j!}{\lfloor j/2\rfloor!^2}.
\end{equation}
The order $d_r$ is generically minimal.  It begins $2,5,16,47,146$.
\end{theorem}

\begin{proof}
The off-diagonal kernel has rank at most $r$.  Its first layer is
$(r,r)$-banded after $r$ condensation steps, and both layers are banded
after $r+1$ steps.  Choose the two local annihilator polynomials
\[
 A(z):=\prod_{s=1}^r(1-v_s z),\qquad
 B(z):=\prod_{s=1}^r(1-v_s^{-1}z).
\]
Anchor the corresponding lower and upper triangular operators at the
first coordinate, as in Proposition~\ref{prop:anchored-minor-preservation}.
They define an ordinary $(r,r)$-banded family $M_n$ preserving
both layers of $\operatorname{Cond}^{n-2}(H_n,J)$.  Since the generators satisfy these recurrences, its transformed part is confined to the first
$(r+1)\times(r+1)$ corner.  Outside that corner the bulk symbol is
\begin{equation}\label{eq:exponential-vapor-bulk}
 c(z)=\delta A(z^{-1})B(z),\qquad
 z^r c(z)=\frac{(-1)^r\delta}{v_1\cdots v_r}
                 \prod_{s=1}^r(z-v_s)^2.
\end{equation}
This specifies an ordinary parent, including its fixed northwest
correction, before applying the row--column transfer.

Assume first that the $v_s$ are distinct and $\delta\ne0$.
The ordinary companion has one size-two Jordan block at each $v_s$.
By Lemma~\ref{lem:exact-boundary-companion}, the determinant transfer is
the middle compound multiplied by $\delta/(v_1\cdots v_r)$.
Choose $m_s\in\{0,1,2\}$ vectors from the $s$th Jordan block, with
$\sum_s m_s=r$.  The resulting eigenvalue is
$\delta\prod_s v_s^{m_s-1}$.  If $e_s:=m_s-1$, the zero entries of
$e$ correspond precisely to the size-two factors selected once.
Their tensor product has largest Jordan length
$1+\#\{s:e_s=0\}$.  This proves
\eqref{eq:exponential-annihilator}.  Counting the choices with $k$
positive and $k$ negative entries gives the first sum in
\eqref{eq:exponential-order}; separating the even and odd terms in the
second sum gives the same expression.  The latter sequence is recorded
as OEIS~A163865 \cite{OEISA163865}.

To determine which geometric terms actually occur in the determinant sequence,
set $C_{ab}:=v_b/(v_b-v_a)$ for $a\ne b$.  A finite-rank determinant
expansion gives
\[
 D_n=\delta^{n-r}\det(\delta I_r+G_n),\qquad
 (G_n)_{ab}:=u_b\sum_{i=1}^n(v_b/v_a)^i.
\]
Its main diagonal entries are $u_a n$; its off-diagonal entries are
$u_b C_{ab}((v_b/v_a)^n-1)$.  For a fixed exponent vector $e$, let
$P:=\{s:e_s=1\}$, $N:=\{s:e_s=-1\}$, and
$z_e:=\#\{s:e_s=0\}$.  The coefficient of the geometric factor
$(\prod_s v_s^{e_s})^n$ inside the small determinant is a polynomial
in $n$ of degree at most $z_e$.
Its leading coefficient is
\begin{equation}\label{eq:exponential-mode-leading}
 \left(\prod_{s=1}^r u_s\right)
     \det C[N,P]\det C[P,N],
\end{equation}
with the two empty determinants equal to one.  Indeed, maximal degree
requires the diagonal term $u_s n$ at every zero coordinate.  The
remaining exponential entries match rows $N$ to columns $P$, and the
remaining constant entries match rows $P$ to columns $N$.  The
block-permutation sign cancels the signs of those constant entries,
giving \eqref{eq:exponential-mode-leading}.

Both displayed minors are nonzero Cauchy determinants when the $v_s$
are distinct.  Thus every geometric factor has a coefficient of the stated maximal polynomial
degree when all $u_s\ne0$.  On the further open set where the finitely
many mode products are distinct, all factors of
\eqref{eq:exponential-annihilator} are necessary.  This proves generic
minimality.  The annihilator identity itself extends to coincident
$v_s$ and other special parameter values after clearing powers of
the nonzero $v_s$; its degree may then exceed the minimal order.
\end{proof}

For $r=1$, the formula is
\[
 D_n=2\delta D_{n-1}-\delta^2D_{n-2}.
\]
For $r=2$, put $\eta:=v_1/v_2+v_2/v_1$.  The factored annihilator is
\begin{equation}\label{eq:exponential-ranktwo-factor}
 \mathcal P_2(z)=(z-\delta)^3(z^2-\eta\delta z+\delta^2),
\end{equation}
and hence
\begin{align}\label{eq:exponential-ranktwo-recurrence}
 D_n={}&(3+\eta)\delta D_{n-1}-(4+3\eta)\delta^2D_{n-2}
 +(4+3\eta)\delta^3D_{n-3}\notag\\
 &-(3+\eta)\delta^4D_{n-4}+\delta^5D_{n-5}.
\end{align}
These are constant-ratio formulas.  For independently varying separable
generators, the ordered products of successive ordinary transfers in
\eqref{eq:local-vapor-observation-rows} supply the appropriate local
coefficients.

\section{Functional kernels and the sampling domain}
\label{sec:functional-classification}

We now classify structural hypotheses stated directly in terms of the
kernel $H_{ij}=f(i,j)$.  Finite-dimensional spaces of shifted sequences
give an exact criterion for Toeplitz tails.  For functions of two real variables, the classification
depends on the sampling domain and steps.  Agreement of the upper and lower formulas on the main diagonal
can further narrow the band; the final exponential-distance example
carries this refinement through to an explicit ordinary recurrence.

\subsection{Separated kernels and finite shift spaces}

The separated-tail criterion applies to functions well beyond
monomials.  For example, for arbitrary scalar samples $u_i,v_j$,
\[
 \cos(u_i+v_j)=\cos u_i\cos v_j-\sin u_i\sin v_j,
\]
so either triangle has separated rank at most two.  The sampling
functions need not be linear, and the two triangles can use
different functions.  Adding an arbitrary diagonal gives a
pentadiagonal first layer after two steps and support within the same five diagonals
in both layers after three steps.  These are upper bounds on the required
depths and widths.  With an
arbitrary core of semibandwidths $(b_-,b_+)$, the corresponding
widths are $(b_-+2,b_++2)$.

The following separated-rank bounds hold over $\mathbb R$ or
$\mathbb C$ for arbitrary samples.  Here $P,P_\lambda$ are polynomials
and sums over $\lambda$ are finite:
\[
\begin{array}{c|c}
 \text{kernel on one triangle}&\text{rank bound}\\ \hline
 f(i)g(j)&1\\
 \exp(u_i+v_j)&1\\
 \cos(u_i\pm v_j),\ \sin(u_i\pm v_j)&2\\
 \cosh(u_i+v_j),\ \sinh(u_i+v_j)&2\\
 P(u_i+v_j)&\deg P+1\\
 P(u_i+v_j)\exp(\lambda(u_i+v_j))&\deg P+1\\
 \sum_{\lambda}P_\lambda(u_i+v_j)
                 \exp(\lambda(u_i+v_j))
   &\sum_\lambda(\deg P_\lambda+1)\\
 P(u_iv_j)&\#\{\text{nonzero monomials of }P\}.
\end{array}
\]
Each assertion follows from an addition formula or the binomial
theorem.  These are bounds: coincident sampled frequencies,
vanishing coefficients, or dependent sampled functions can reduce
the actual rank.  For instance, a single monomial $(u_iv_j)^d$
has rank one although $(u_i+v_j)^d$ generically has rank $d+1$.

For the special forms $a_{i+j}$ and $a_{j-i}$ there is also a
complete finite-rank classification.  The following elementary
discrete result is a counterpart of finite-dimensional
translation-invariance theory; a general reference is
\cite[Theorem~1]{Engert1970}.

\begin{theorem}\label{thm:bilateral-shift-classification}
Let $a=(a_n)_{n\in\mathbb Z}$ be a nonzero sequence over a
characteristic-zero field $\mathbb F$, and let
$(\mathcal E a)_n:=a_{n+1}$.  The following are equivalent:
\begin{enumerate}
 \item The shift space
 $V(a):=\operatorname{span}_{\mathbb F}\{\mathcal E^ka:k\in\mathbb Z\}$
 has finite dimension $r$.
 \item The bilateral Hankel matrix $(a_{i+j})_{i,j\in\mathbb Z}$
 has finite rank $r$, where rank is the supremum of the ranks
 of its finite submatrices.
 \item The bilateral Toeplitz matrix $(a_{j-i})_{i,j\in\mathbb Z}$
 has finite rank $r$.
 \item The minimal annihilating polynomial $\chi$ in
 $\chi(\mathcal E)a=0$ has degree $r$ and nonzero constant coefficient.
 \item Over an algebraic closure,
 \begin{equation}\label{eq:bilateral-exponential-polynomial}
  a_n=\sum_{\ell=1}^s P_\ell(n)\lambda_\ell^n
       \qquad(n\in\mathbb Z),
 \end{equation}
 where the $\lambda_\ell$ are distinct and nonzero, the
 $P_\ell$ are nonzero polynomials, and
 \[
  r=\sum_{\ell=1}^s(\deg P_\ell+1).
 \]
\end{enumerate}
In this case
\[
 \chi(z)=\prod_{\ell=1}^s
              (z-\lambda_\ell)^{\deg P_\ell+1}.
\]
\end{theorem}

\begin{proof}
The Hankel columns are the shifts of $a$, and row reversal gives the
Toeplitz matrix.  Independence of finitely many sequences is detected
on finitely many coordinates, so both ranks equal $\dim V(a)$.
The shift is invertible on this space.  A degree-$r$ annihilator with
nonzero constant coefficient expresses every forward and backward shift
through $a,\mathcal E a,\ldots,\mathcal E^{r-1}a$; minimality is equivalent to independence
of these $r$ vectors.

Over an algebraic closure, powers of each Jordan block of the cyclic
shift give a polynomial times its nonzero eigenvalue to the integer
power.  Evaluation at coordinate zero gives
\eqref{eq:bilateral-exponential-polynomial}.
Conversely, $\mathcal E-\lambda$ lowers the polynomial degree multiplying
$\lambda^n$ by one until it vanishes.  The factors for distinct
$\lambda$ are coprime, proving the asserted minimal annihilator and
its degree.
\end{proof}

This classifies bilateral sampled sequences.  Functions away from
their samples remain unrestricted, and one-sided sequences may have
finite initial transients.  For arbitrary arrays, a vanishing connected
minor layer alone need not imply finite rank.

For a Toeplitz symbol whose upper and lower tails are independently
exponential-polynomial, put
\[
 b_s=\sum_\ell P_\ell^+(s)(\lambda_\ell^+)^s
       \quad(s>b_+),\qquad
 b_{-s}=\sum_\ell P_\ell^-(s)(\lambda_\ell^-)^s
       \quad(s>b_-).
\]
Each expression in $j-i$ separates by the binomial theorem.
Consequently the upper and lower ranks are bounded by
$r_\pm:=\sum_\ell(\deg P_\ell^\pm+1)$, and
Theorem~\ref{thm:band-core-compression} gives first-layer depth
$\max(r_-,r_+)$ and width
$(b_-+r_-,b_++r_+)$.  This bound may improve when the formulas
extend through part of the prescribed core, when the effective
ranks are smaller, or when a sharper tail recurrence is known.

Finite sums of sampled trigonometric functions are included by
writing sine and cosine as sums of two complex exponentials.
The rank calculation can be made over $\mathbb C$ without changing
the rank of a real matrix.  A nonresonant term
$n^d\cos(Cn)$ has rank $2(d+1)$; when $e^{iC}=e^{-iC}$, the
two frequencies coincide and the rank can drop.  Polynomial
modulation and repeated frequencies are therefore accounted for
by the multiplicities in \eqref{eq:bilateral-exponential-polynomial}.

\subsection{Exact condensation depth for Toeplitz tails}
\label{sec:toeplitz-tail-classification}

For a one-sided sequence $a=(a_n)_{n\ge0}$, its eventual recurrence
order is the least degree $r$ of a monic polynomial $\chi$ such that
$\chi(0)\ne0$ and $(\chi(\mathcal E)a)_n=0$ for every sufficiently
large $n$.
An eventually zero sequence has order zero; if no such polynomial
exists, the order is infinite.  In characteristic zero a finite
eventual order is equivalent to an exponential-polynomial formula
after a finite initial segment.

\begin{lemma}\label{lem:eventual-connected-hankel-rank}
For $t\ge1$ set
\[
 h_t(n):=\det(a_{n+i+j})_{0\le i,j<t},\qquad h_0(n):=1.
\]
The eventual recurrence order of $a$ is at most $k$ if and only if
$h_{k+1}(n)=0$ for every sufficiently large $n$.
\end{lemma}

\begin{proof}
An order-$r$ recurrence bounds every sufficiently late Hankel rank by
$r$.  Conversely, choose the least $r\le k$ such that $h_{r+1}$ is
eventually zero.  The case $r=0$ is immediate.  For $r\ge1$,
Desnanot--Jacobi gives
\[
 h_r(n+1)^2=h_r(n)h_r(n+2)
\]
for all sufficiently large $n$.  Minimality supplies arbitrarily late
nonzero terms.  Each interior nonzero term forces both neighbours to
be nonzero, so $h_r(n)\ne0$ throughout a sufficiently late ray $n\ge N$.

Every finite rectangle of $(a_{N+i+j})_{i,j\ge0}$ has zero connected
$(r+1)$-minors and nonzero connected $r$-minors.  By
Lemma~\ref{lem:local-connected-rank}, its first $r$ columns form a
basis whenever the rectangle is large enough.  The unique coefficients
on their first $r$ rows therefore give one relation on the entire array:
\[
 a_{n+r}=\sum_{j=0}^{r-1}c_j a_{n+j}\qquad(n\ge N).
\]
Substituting into the adjacent $r\times r$ minor gives
\[
 c_0=(-1)^{r-1}\frac{h_r(N+1)}{h_r(N)}\ne0.
\]
A smaller recurrence would force $h_r$ eventually to vanish.  Hence
this recurrence has the required minimal order.
\end{proof}

\begin{theorem}\label{thm:toeplitz-bandability-classification}
Let $b=(b_s)_{s\in\mathbb Z}$ be an infinite Toeplitz symbol, and
consider all its growing finite principal sections.  Condensation
has a fixed-width first layer at a fixed depth, uniformly in the
section size, if and only if both one-sided tails $(b_s)_{s\ge0}$
and $(b_{-s})_{s\ge0}$ have finite eventual recurrence orders
$r_+,r_-$.  If these orders are finite, the exact minimum depth is
\[
 R:=\max(r_+,r_-).
\]
The exact minimum depth at which both layers have fixed-width
support is $R+1$.  Each minimum is taken over all finite bandwidths,
with depth and width independent of the growing section size.
\end{theorem}

\begin{proof}
Put
\[
 C_t(s;b):=\det(b_{s+j-i})_{0\le i,j<t}.
\]
Reversing its rows gives
\[
 C_t(s;b)=(-1)^{t(t-1)/2}
       \det(b_{s-t+1+i+j})_{0\le i,j<t}.
\]
Thus $C_t(s;b)$ vanishes for every sufficiently large positive
$s$ if and only if the positive tail has eventual recurrence
order at most $t-1$, by
Lemma~\ref{lem:eventual-connected-hankel-rank}.
Transposition gives the identical equivalence for the negative
tail and sufficiently negative $s$.

At depth $k$, the first layer has symbol $C_{k+1}(s;b)$.
It has finite support exactly when $k\ge r_+$ and $k\ge r_-$,
which proves the first assertion and its exact minimum.
For $k\ge1$ the second symbol is $C_k(s;b)$, so both have
finite support exactly when $k-1\ge R$.  At depth zero the
ordinary denominator has symbol $C_0(s;b)=1$ and is not
finitely supported.  Hence the minimum two-layer depth is $R+1$.
\end{proof}

The theorem uses the ordinary connected-minor interpretation for all
leading sections of one infinite symbol, including singular steps;
a single finite-size zero pattern gives no such converse.  The depth
is exact, while the width depends on the initial tail transients.
If separated formulas start beyond offsets $b_\pm$,
Theorem~\ref{thm:band-core-compression} bounds that width by
$(b_-+r_-,b_++r_+)$.

For example, in characteristic zero the tails of the symbol
$b_s=|s|^p$ for $s\ne0$, with any prescribed $b_0=x$ and integer
$p\ge0$, have minimal eventual recurrence order $p+1$.
Its exact first-layer and two-layer compression depths are therefore
$p+1$ and $p+2$, for every $x$.  A special value of $x$ can change
the surviving entries or width, but does not change these minimum
depths over all fixed widths.

\begin{example}\label{ex:noncompressing-quadratic-exponential}
The function being explicit or analytic does not by itself imply
finite-depth compression.  Let $q\ne0$, assume $q^2$ is not a root
of unity, and take $H_{ij}:=q^{(i+j)^2}$.  For consecutive row
indices $u_a=i+a$ and column indices $v_b=j+b$,
\[
 q^{(u_a+v_b)^2}=q^{u_a^2}q^{v_b^2}
                 (q^{2u_a})^{j+b}.
\]
Removing the nonzero row and column factors leaves a Vandermonde
matrix in the distinct nodes $q^{2u_a}$.  Every connected square
minor is therefore nonzero.  No fixed condensation depth produces
a fixed-width band in arbitrarily large sections.  The same
argument, with nodes $q^{-2u_a}$, applies to the Toeplitz kernel
$q^{(j-i)^2}$.
\end{example}

\subsection{Functional equations and the sampling domain}
\label{sec:functional-equations}

There are three sampling requirements: an identity on the integer
grid, an identity at every real base point with a fixed unit step, and
an identity with independently variable row and column steps.  Their
solution classes differ.

Continuity alone places no restriction on the values of a matrix sampled
on $\mathbb Z^2$.  Indeed, choose a smooth function $\phi$ supported in
$(-1/3,1/3)$ with $\phi(0)=1$.  For any prescribed array $(a_{ij})$,
the locally finite sum
\[
 f(x,y):=\sum_{i,j\in\mathbb Z}a_{ij}\phi(x-i)\phi(y-j)
\]
is smooth and satisfies $f(i,j)=a_{ij}$.  Adding
$\sin(\pi x)H(x,y)$ for continuous $H$ also leaves every integer sample
unchanged.  A functional classification therefore requires identities
beyond that grid or a specified function class.

\subsubsection*{Functional equations at a fixed lattice step}
\label{sec:feq-fixed-step}

Let $f\colon\mathbb R^2\to\mathbb R$ be continuous, and write
\[
 W_f(x,y):=\det\bigl(f(x+i,y+j)\bigr)_{0\le i,j\le1},\qquad
 Z_f(x,y):=\det\bigl(f(x+i,y+j)\bigr)_{0\le i,j\le2}.
\]
The equations below hold at every real base point, with step one.
A kernel has \emph{separated rank at most $r$} if
$f(x,y)=\sum_{k=1}^r u_k(x)v_k(y)$ on its whole domain.  Its
\emph{unit-lattice orbit} at $(x,y)$ is the array
$\bigl(f(x+m,y+n)\bigr)_{m,n\in\mathbb Z}$.

\begin{proposition}\label{prop:feq-unit-rank-one}
Suppose that $f$ is nowhere zero.  Then $W_f\equiv0$ if and only if
\begin{equation}\label{eq:feq-periodic-log-form}
 f(x,y)=\varepsilon\exp\bigl(P(x,y)+Q(x,y)\bigr),
 \qquad \varepsilon\in\{1,-1\},
\end{equation}
where $P,Q$ are continuous and
\[
 P(x+1,y)=P(x,y),\qquad Q(x,y+1)=Q(x,y).
\]
Every unit-lattice orbit then has rank one.  The functions $P$ and $Q$
need not be functions of only one variable.
\end{proposition}

\begin{proof}
The sign $\varepsilon$ is constant.  Put $g:=\log|f|$ and define
$\Delta_xg(x,y):=g(x+1,y)-g(x,y)$, with $\Delta_y$ defined analogously.
The equation $W_f=0$ is equivalent to $\Delta_x\Delta_yg=0$.
Thus $h:=\Delta_xg$ is one-periodic in $y$.  There is a continuous
one-periodic-in-$y$ solution of $\Delta_xQ=h$: first put
$Q(s,y):=s h(0,y)$ for $0\le s\le1$, and extend in both directions
by $Q(x+1,y):=Q(x,y)+h(x,y)$.  The endpoint values agree, so this
extension is continuous.  Now $P:=g-Q$ is one-periodic in $x$.
Conversely, either periodicity makes the corresponding mixed difference
zero.  Iterating the mixed-difference identity gives, for all integers
$m,n$,
\[
 f(x+m,y+n)f(x,y)=f(x+m,y)f(x,y+n),
\]
which proves the orbit-rank assertion.
\end{proof}

\begin{example}\label{ex:feq-analytic-rank-one-orbits}
The positive real-analytic function
\[
 f(x,y):=\exp\bigl(x\sin(2\pi y)\bigr)
\]
has $W_f\equiv0$, since it is one-periodic in $y$, but has infinite
separated rank.  Indeed, choosing arbitrarily many $y_j$ with distinct
$t_j:=\sin(2\pi y_j)$ gives the linearly independent functions
$x\mapsto e^{t_jx}$.  Their independence follows, for example, from
the Vandermonde matrix of their derivatives at zero.  On the integer
grid the same function is identically one.
\end{example}

\begin{proposition}\label{prop:feq-unit-rank-two}
Suppose that $W_f$ is nowhere zero.  Then $Z_f\equiv0$ if and only if
there are continuous functions $a,b$, one-periodic in $y$, with $b$
nowhere zero, such that
\begin{equation}\label{eq:feq-transverse-periodic-recurrence}
 f(x+2,y)=a(x,y)f(x+1,y)+b(x,y)f(x,y).
\end{equation}
They are uniquely determined by
\begin{align}\label{eq:feq-rank-two-coefficients}
 a(x,y)&:=
 \frac{f(x,y)f(x+2,y+1)-f(x,y+1)f(x+2,y)}{W_f(x,y)},\notag\\
 b(x,y)&:=-\frac{W_f(x+1,y)}{W_f(x,y)}.
\end{align}
Every unit-lattice orbit has rank exactly two.  The analogous column
recurrence has coefficients one-periodic in $x$.
\end{proposition}

\begin{proof}
Cramer's rule expresses the third row through the first two on columns
$y,y+1$, with the displayed coefficients.  If $Z_f(x,y)=0$, this
relation also holds at $y+2$.  Since $W_f(x,y+1)\ne0$, uniqueness on
the overlapping columns gives $a(x,y+1)=a(x,y)$ and
$b(x,y+1)=b(x,y)$.  Conversely, this periodicity makes the row
relation valid on all three columns and gives $Z_f=0$.

On each orbit the coefficients depend only on the row index.  Since
$b\ne0$, the recurrence propagates a two-row span in both directions;
$W_f\ne0$ makes its dimension exactly two.  Transposition gives the
column assertion.
\end{proof}

\begin{example}\label{ex:feq-analytic-rank-two-orbits}
Put $s(y):=\sin(2\pi y)$ and
\[
 f(x,y):=e^{x s(y)}\bigl(1+e^{x+y}\bigr).
\]
This positive real-analytic function satisfies
\[
 Z_f\equiv0,\qquad
 W_f(x,y)=(e-1)^2
 \exp\bigl((2s(y)+1)x+y+s(y)\bigr)>0.
\]
Indeed, on a unit-lattice orbit the periodic function $s(y)$ is fixed,
and the displayed kernel is a nonzero row factor times a rank-two sum.
Nevertheless its separated rank is infinite.  Choose arbitrarily many
$y_j$ for which the distinct values $s(y_j)$ lie in $(-1/4,1/4)$.
The slices
\[
 f(x,y_j)=e^{s(y_j)x}+e^{y_j}e^{(s(y_j)+1)x}
\]
have pairwise disjoint sets of exponential exponents and are linearly
independent.  Thus even positivity, real analyticity, and a nowhere-zero
preceding minor do not remove the fixed-step freedom.
\end{example}

\begin{proposition}\label{prop:feq-analytic-orbit-rank}
If $f$ is real analytic and $Z_f\equiv0$, every unit-lattice orbit has
rank at most two.  If $W_f\equiv0$, every such orbit has rank at most
one.
\end{proposition}

\begin{proof}
If $W_f\not\equiv0$, any finite rectangle of integer offsets has all
its translated $W_f$ simultaneously nonzero on a dense open set of
base points.  Row propagation there gives rank at most two.  Every
three-minor with fixed integer offsets then vanishes everywhere by
continuity.  Since the offsets are arbitrary, every orbit has rank at
most two.  If $W_f\equiv0$ but $f\not\equiv0$, use the dense open set
where the finitely many translated entries are nonzero and apply the
rank-one argument.  The zero kernel is immediate.
\end{proof}

\begin{remark}\label{rem:feq-zero-loci-and-core}
The analytic hypothesis in Proposition~\ref{prop:feq-analytic-orbit-rank}
cannot be replaced by smoothness.  Choose a smooth function $\phi$
supported in $(-1/4,1/4)$ with $\phi(0)=1$, and put
\[
 f(x,y):=\sum_{k\ge1}\phi(x-2k)\phi(y-2k).
\]
The sum is locally finite.  Every unit square has at most one nonzero
corner, so $W_f\equiv0$, whereas $f(2i,2j)=\delta_{ij}$ for positive
integers $i,j$.  Its integer-grid rank is unbounded.

There is also a distinction between a globally analytic kernel and a
matrix with a separately prescribed diagonal or band core.  For globally
analytic $f$, vanishing of $W_f$ or $Z_f$ on any nonempty open region of
real base points forces its vanishing everywhere by analytic continuation.
The corresponding sampled matrices then have rank at most one or two.
An off-band identity for all real base points therefore imposes more than
the discrete support condition used for a freely prescribed band core.
\end{remark}

\begin{proposition}\label{prop:feq-rational-classification}
Let $K$ be a field of characteristic zero, let $f\in K(x,y)$, and fix
$r\ge1$.  Then the rational identity
\begin{equation}\label{eq:feq-rational-unit-minor}
 \det\bigl(f(x+i,y+j)\bigr)_{0\le i,j\le r}=0
\end{equation}
holds if and only if
\begin{equation}\label{eq:feq-rational-separated-form}
 f(x,y)=\sum_{k=1}^{r}u_k(x)v_k(y),
 \qquad u_k\in K(x),\quad v_k\in K(y),
\end{equation}
where zero summands are allowed.  If $f$ is polynomial, the factors may
be chosen polynomial.  Equivalently, for a polynomial
$f(x,y)=\sum_{i,j}c_{ij}x^iy^j$, the coefficient matrix $(c_{ij})$ has
rank at most $r$.
\end{proposition}

\begin{proof}
The separated form makes every $(r+1)$-minor vanish.  Conversely,
assume $f\ne0$ and choose the least $s\le r$ for which the connected
minor of order $s+1$ is identically zero.  Its order-$s$ predecessor
and all its translates are nonzero over $K(x,y)$.
Lemma~\ref{lem:local-connected-rank} therefore bounds the rank of every
finite rectangle of integer offsets by $s$.  All its $(s+1)$-minors,
with arbitrary integer row and column offsets, vanish rationally.

Clear denominators in $\det(f(X_i,Y_j))_{0\le i,j\le s}$, obtaining
a polynomial $F$ in its $2s+2$ arguments.  For every integer pair of
vectors $m,n$, the preceding identities give
$F(x+m_0,\ldots,x+m_s,y+n_0,\ldots,y+n_s)=0$.
Setting $x=y=0$ shows that $F$ vanishes on the integer product grid,
which is Zariski dense in characteristic zero by successive use of the
univariate root bound.  Thus $F=0$.

Choose a nonsingular minor $B$ of maximum order $t\le s$ at scalar
points avoiding the poles; these exist because $K$ is infinite.
The vanishing bordered minors give
\[
 f(x,y)=
 \bigl(f(x,y_1),\ldots,f(x,y_t)\bigr)B^{-1}
 \bigl(f(x_1,y),\ldots,f(x_t,y)\bigr)^T.
\]
These factors are rational, and polynomial when $f$ is polynomial.
Factoring the polynomial coefficient matrix proves the final
characterization.  The zero kernel needs no argument.
\end{proof}

\begin{remark}\label{rem:feq-rational-denominators}
In particular, the rational functions satisfying a fixed-step equation
of some finite order are precisely those whose reduced denominator is
a product of a polynomial in $x$ and a polynomial in $y$.  A finite
separated sum has such a common denominator; conversely, division of a
polynomial numerator by such a product gives a finite separated sum.
Thus the fixed-step equations classify rational and polynomial kernels
more rigidly than continuous or analytic kernels.
\end{remark}

\subsubsection*{Sampling with variable steps}

Let $I,J\subset\mathbb R$ be nonempty open intervals,
$\mathbb K\in\{\mathbb R,\mathbb C\}$, and
$f:I\times J\to\mathbb K$.  For $m\ge1$ define
\begin{equation}\label{feq:variable-step-determinant}
 \mathcal D_m(f;x,y;h,k)
 :=\det\bigl(f(x+ih,y+jk)\bigr)_{0\le i,j<m},
\end{equation}
whenever all the indicated arguments belong to $I\times J$.
The steps $h,k$ may be chosen independently.  We write separated representations in the form
\begin{equation}\label{feq:separated-representation}
 f(x,y)=\sum_{a=1}^r u_a(x)v_a(y).
\end{equation}
Terms equal to zero are allowed.  When the required derivatives exist, put
\begin{equation}\label{feq:mixed-wronskian}
 \mathcal W_m(f;x,y)
 :=\det\bigl(\partial_x^i\partial_y^j f(x,y)\bigr)_{0\le i,j<m}.
\end{equation}
Here $C^{r,r}$ means that all mixed derivatives with at most $r$
derivatives in each variable exist and are continuous.

\begin{lemma}\label{feq:regular-differential-classification}
Let $r\ge1$ and $f\in C^{r,r}(I\times J)$.  If
$\mathcal W_{r+1}(f)=0$ and $\mathcal W_r(f)\ne0$ everywhere on
$I\times J$, then $f$ has separated rank exactly $r$ there.
More precisely, all functions $x\mapsto f(x,y)$ satisfy one common
ordinary differential equation
\begin{equation}\label{feq:common-ode}
 \partial_x^r f(x,y)
 =\sum_{i=0}^{r-1}a_i(x)\partial_x^i f(x,y),
\end{equation}
whose coefficients are continuous on $I$.
\end{lemma}

\begin{proof}
At every $(x,y)$, solve the nonsingular system
\begin{equation}\label{feq:ode-coefficient-system}
 \partial_x^r\partial_y^j f
 =\sum_{i=0}^{r-1}a_i(x,y)\partial_x^i\partial_y^j f,
 \qquad 0\le j<r.
\end{equation}
The equation $\mathcal W_{r+1}(f)=0$ supplies the same identity for
$j=r$.  Differentiate the identity with index $j<r$ in $y$ and
subtract the identity with index $j+1$.  This gives
\[
 \sum_{i=0}^{r-1}(\partial_y a_i)
       \partial_x^i\partial_y^j f=0,
 \qquad 0\le j<r.
\]
Since the coefficient determinant is $\mathcal W_r(f)\ne0$, every
$\partial_y a_i$ is zero.  Thus $a_i=a_i(x)$ and
\eqref{feq:common-ode} follows.  Choose a fundamental system
$u_1,\ldots,u_r$ of this ordinary differential equation on $I$.
Uniqueness for initial data at a fixed $x_0\in I$ expresses each
$f(\,\cdot\,,y)$ in that system, with coefficients $v_a(y)$.
This proves \eqref{feq:separated-representation}.  Its rank cannot be
smaller than $r$, because $\mathcal W_r(f)$ is nonzero.
\end{proof}

\begin{theorem}\label{feq:analytic-variable-step-classification}
Let $r\ge1$ and let $f:I\times J\to\mathbb K$ be jointly real analytic.
The following conditions are equivalent:
\begin{enumerate}[label=\textup{(\roman*)}]
 \item $\mathcal D_{r+1}(f;x,y;h,k)=0$ for every admissible
       $x,y$ and every admissible pair of nonzero steps $h,k$;
 \item $\mathcal W_{r+1}(f;x,y)=0$ for every $(x,y)\in I\times J$;
 \item $f$ has separated rank at most $r$ on $I\times J$.
\end{enumerate}
In \textup{(iii)}, the functions $u_a,v_a$ may be chosen real analytic.
\end{theorem}

\begin{proof}
For \textup{(i)}$\Rightarrow$\textup{(ii)}, successive row and column
forward differences preserve the determinant.  Divide row $i$ by
$h^i$ and column $j$ by $k^j$, then take limits:
\begin{equation}\label{feq:difference-wronskian-limit}
 \lim_{h,k\to0}
 \frac{\mathcal D_{r+1}(f;x,y;h,k)}
      {(hk)^{r(r+1)/2}}
 =\mathcal W_{r+1}(f;x,y).
\end{equation}

For \textup{(ii)}$\Rightarrow$\textup{(iii)}, induct on $r$, with
$\mathcal W_1(f)=f=0$ as base case.  If $\mathcal W_r$ is identically
zero, apply induction.  Otherwise
Lemma~\ref{feq:regular-differential-classification} gives rank $r$ on
an open rectangle where $\mathcal W_r\ne0$.  Thus
$\det(f(x_i,y_j))_{0\le i,j\le r}$ vanishes there for independently
chosen arguments.  As an analytic function on the connected domain
$I^{r+1}\times J^{r+1}$, it vanishes everywhere.

Choose a nonzero minor $M:=(f(\xi_i,\eta_j))_{1\le i,j\le s}$ of
maximum order $s\le r$.  Vanishing of its bordered minors gives
\begin{equation}\label{feq:pivot-representation}
 f(x,y)
 =\begin{bmatrix}f(x,\eta_1)&\cdots&f(x,\eta_s)\end{bmatrix}
 M^{-1}
 \begin{bmatrix}f(\xi_1,y)\\ \vdots\\ f(\xi_s,y)\end{bmatrix}.
\end{equation}
This gives analytic separated factors; take $s=0$ when $f=0$.
The bordered-minor argument itself requires no regularity.
Finally, \textup{(iii)}$\Rightarrow$\textup{(i)} follows by factoring
each sampled matrix through an $r$-dimensional space.
\end{proof}

For the first condensation step, arbitrary $h,k$ simply allow every
pair of row points and every pair of column points.  Thus separated
rank at most one follows without continuity or analyticity.  On an
open rectangle where a real-valued $f$ is nonzero, the corresponding
differential equation is
\[
 f f_{xy}-f_x f_y=0,
 \qquad\text{equivalently}\qquad
 (\log|f|)_{xy}=0.
\]
For the second step, the analytic classification is exactly
$f(x,y)=u_1(x)v_1(y)+u_2(x)v_2(y)$, and its differential equation is
\begin{equation}\label{feq:second-step-pde}
 \det\begin{bmatrix}
 f&f_y&f_{yy}\\
 f_x&f_{xy}&f_{xyy}\\
 f_{xx}&f_{xxy}&f_{xxyy}
 \end{bmatrix}=0.
\end{equation}
For a $C^{2,2}$ function this classification remains valid on any
rectangle on which $f f_{xy}-f_xf_y$ is everywhere nonzero.
In the analytic theorem, in contrast, zeros of this lower determinant
need not be excluded.

\begin{example}\label{feq:smooth-arithmetic-counterexample}
Smoothness alone does not give the global conclusion of
Theorem~\ref{feq:analytic-variable-step-classification}, even at the
second step.  Fix $0<\delta<1/4$ and a smooth real function $\chi$
supported in $(-\delta,\delta)$, with $\chi(0)=1$.  Put
\begin{equation}\label{feq:infinite-rank-smooth-kernel}
 f(x,y):=\sum_{m\ge1}\chi(x-3^m)e^{my}.
\end{equation}
The sum is locally finite, so $f\in C^\infty(\mathbb R^2)$.
No three distinct support intervals can be met by a three-term
arithmetic progression.  Indeed, that would give distinct exponents
$a<b<c$ with
\[
 |3^a+3^c-2\cdot3^b|<4\delta<1.
\]
The integer inside the absolute value is nonzero: after division by
$3^a$, its residue modulo $3$ is $1$.  Hence every arithmetic triple
of row points supplies vectors in at most two of the one-dimensional
spaces spanned by $e^{my}$, with zero rows allowed.  Therefore
$\mathcal D_3(f;x,y;h,k)=0$ for all $x,y,h,k$.
Nevertheless, $f(3^m,y)=e^{my}$, and these functions are linearly
independent for arbitrarily many values of $m$.  The separated rank
of $f$ is infinite.

The failure is not confined to kernels whose lower mixed determinant
vanishes everywhere.  Replacing the summand in
\eqref{feq:infinite-rank-smooth-kernel} by
\[
 \chi(x-3^m)\bigl(e^{2my}+(x-3^m)e^{(2m+1)y}\bigr)
\]
gives the same counterexample, with
$\mathcal W_2=\chi(x-3^m)^2e^{(4m+1)y}$ on the $m$th strip.
Three nonzero sampled rows must now lie in one strip: three distinct
strips are excluded as above, and a progression meeting one strip
twice has too small a step to meet a second strip.  Its rank is
therefore at most two.  This also explains why regularity at a single
point is insufficient for a global smooth classification.
\end{example}

For any function, vanishing of all $(r+1)$-minors at independently
chosen row and column points is equivalent to separated rank at most
$r$, by \eqref{feq:pivot-representation}.  Arithmetic sampling is
weaker; analyticity makes the two conditions equivalent.

\subsubsection*{Sum and difference kernels}

For kernels $f(x,y)=g(x+y)$, the continuous all-step problem has a
classical complete solution.  Its one-variable equation is the
Popoviciu--Ionescu equation
\begin{equation}\label{feq:popoviciu-ionescu}
 \det\bigl(g(t+(i+j)h)\bigr)_{0\le i,j\le r}=0
 \qquad(t,h\in\mathbb R).
\end{equation}
The solution theory for continuous real-valued $g$ gives exponential
polynomials; see \cite{Almira2016}.  More precisely, its solutions have
the form
\begin{equation}\label{feq:sum-exponential-classification}
 g(t)=\sum_{\nu=1}^{s}P_\nu(t)e^{\lambda_\nu t},
 \qquad
 \sum_{\nu=1}^{s}(\deg P_\nu+1)\le r,
\end{equation}
with distinct complex $\lambda_\nu$, nonzero polynomials $P_\nu$, and
conjugate terms paired when needed to give real values; the zero function
is allowed.  To see the order bound from the exponential-polynomial
conclusion, let $d$ be the sum of the multiplicities.  The addition
formula and linear independence of exponential monomials give separated
rank exactly $d$ for $g(x+y)$.  Taking $h\to0$ in
\eqref{feq:popoviciu-ionescu} makes its mixed determinant
$\mathcal W_{r+1}$ vanish, so
Theorem~\ref{feq:analytic-variable-step-classification} gives $d\le r$.
Conversely, the separated form makes every displayed determinant zero.

For $r=1$, the nonzero real solutions are $g(t)=Ae^{\lambda t}$.
For $r=2$, the real solutions, including lower-rank degenerations, are
\begin{align*}
 g(t)&=Ae^{\lambda t}+Be^{\mu t},\\
 g(t)&=(A+Bt)e^{\lambda t},\\
 g(t)&=e^{\alpha t}\bigl(A\cos(\beta t)+B\sin(\beta t)\bigr),
 \qquad \beta\ne0,
\end{align*}
with real parameters.  This explains the exponential and trigonometric
families as solutions of a functional equation, in addition to their
direct separated formulas.  For $g(y-x)$ the same conclusion follows
by reversing the order of the sampled rows.

The quantifier over $h$ is necessary.  For example,
$g(t):=\exp(t\sin(2\pi t))$ satisfies the first equation with the
single step $h=1$, since $g(t+n)=\exp((t+n)\sin(2\pi t))$ for
integer $n$, but it is not a single exponential.  Finally, a globally
finite-rank kernel gives zero determinants in large orders.  Nontrivial
determinant sequences arise by prescribing a separate diagonal or band
core, or by matching different finite-rank formulas on the two sides
of the diagonal, as in the following application.

\subsection{Matching triangular kernels and sharper bands}
\label{sec:continuous-matching-kernels}

Continuity forces the two triangular formulas to agree on the sampled
diagonal, narrowing the band obtained with an independently prescribed
diagonal.

\begin{theorem}\label{thm:continuous-triangular-compression}
Let $I\subseteq\mathbb R$ be an interval, let $f:I^2\to\mathbb C$ be
continuous, and suppose that continuous one-variable functions give
representations
\[
 f(x,y)=\sum_{\nu=1}^{r_+}u_\nu^+(x)v_\nu^+(y)\quad(x<y),\qquad
 f(x,y)=\sum_{\nu=1}^{r_-}u_\nu^-(x)v_\nu^-(y)\quad(x>y),
\]
where $r_+,r_-\ge1$.  Sample at any strictly increasing points
$t_1<t_2<\cdots$ in $I$, and put $H_{ij}:=f(t_i,t_j)$ and
\[
 C_k(i,j):=\det H[i:i+k-1,j:j+k-1].
\]
Then
\begin{align}\label{eq:continuous-triangular-support}
 C_k(i,j)&=0&&\text{if }k>r_+\text{ and }j-i\ge r_+,\notag\\
 C_k(i,j)&=0&&\text{if }k>r_-\text{ and }i-j\ge r_-.
\end{align}
Consequently, with $R:=\max(r_-,r_+)$, the first layer after $R$
condensations has semibandwidths at most $(r_--1,r_+-1)$; after
$R+1$ condensations both layers have these semibandwidths, whenever
the indicated layers exist.
\end{theorem}

\begin{proof}
Continuity extends the upper formula to $x=y$, and likewise for the
lower formula.  Set $s:=j-i\ge r_+$.  If $k\le s+1$, the entire
$k\times k$ submatrix is in the closed upper triangle and has rank
at most $r_+<k$.  If $k>s+1$, its first $s+1$ rows lie in that
triangle and have rank at most $r_+$.  The remaining $k-s-1$ rows
increase the rank by at most that number, so its rank is at most
$r_++k-s-1<k$.  This proves the upper vanishing statement; transpose
for the lower statement.  At depth $R$ the first layer consists of
$C_{R+1}$; at depth $R+1$ the two layers consist of $C_{R+2}$ and
centrally shifted $C_{R+1}$.
\end{proof}

The reduction by one uses the matching diagonal; an independently
assigned diagonal is governed by Theorem~\ref{thm:band-core-compression}.
Both statements require finite triangular rank.

\paragraph{A rank-one continuous kernel.}
For continuous $u,v:I\to\mathbb C$, consider
\[
 f(x,y):=u(\min(x,y))v(\max(x,y)).
\]
Write $u_i:=u(t_i)$, $v_i:=v(t_i)$, and
$\delta_i:=u_{i+1}v_i-u_iv_{i+1}$.  The first condensed layer is
diagonal, with
\[
 C_2(i,i)=u_iv_{i+1}\delta_i,
 \qquad C_2(i,j)=0\quad(i\ne j).
\]
More generally,
\begin{equation}\label{eq:continuous-green-minors}
 C_k(i,i)=u_iv_{i+k-1}\prod_{\ell=i}^{i+k-2}\delta_\ell,
 \qquad C_k(i,j)=0\quad(i\ne j,\ k\ge2).
\end{equation}
The diagonal formula follows inductively from Desnanot--Jacobi, first
where its denominators are nonzero and then as a polynomial identity.
In particular,
\begin{equation}\label{eq:continuous-green-determinant}
 D_n:=\det H_n=u_1v_n\prod_{\ell=1}^{n-1}\delta_\ell.
\end{equation}
Thus continuity and triangular separation do not force the original
matrix to have finite rank: for instance $u(t)=e^t$, $v(t)=e^{-t}$
gives $f(x,y)=e^{-|x-y|}$ and nonzero determinants at all orders.

\begin{proposition}\label{prop:matched-mixture-exact-depth}
Let $c_1,\ldots,c_r>0$, and let $q_1,\ldots,q_r\in(0,1)$ be
distinct.  For the integer samples of the continuous kernel
\[
 f(x,y):=\sum_{\nu=1}^r c_\nu q_\nu^{|x-y|},
\]
put $H_n:=(f(i,j))_{i,j=1}^n$ and $D_n:=\det H_n$.  The
minimum uniform depths for a first banded layer and for two
banded layers are $r$ and $r+1$, respectively.  At these depths
the surviving semibandwidth is exactly $r-1$.
\end{proposition}

\begin{proof}
The upper bounds follow from
Theorem~\ref{thm:continuous-triangular-compression}.  Write
\[
 h_s:=\sum_\nu c_\nu q_\nu^{|s|},\qquad
 C_k(s):=\det(h_{s+j-i})_{i,j=0}^{k-1}.
\]
For $k\le r$ and
$s\ge k-1$, Cauchy--Binet gives
\begin{equation}\label{eq:matched-mixture-far-minors}
 C_k(s)=(-1)^{k(k-1)/2}
 \sum_{\substack{I\subseteq\{1,\ldots,r\}\\|I|=k}}
 \left(\prod_{\nu\in I}c_\nu q_\nu^{s-k+1}\right)
 \prod_{\substack{\mu<\nu\\\mu,\nu\in I}}(q_\nu-q_\mu)^2.
\end{equation}
Every summand has the same sign and is nonzero.  Thus no layer
before depth $r$ has a fixed finite width; at depth $r$ its
denominator still has the nonzero far minors $C_r(s)$.

Put $\delta:=\sum_\nu c_\nu(q_\nu-q_\nu^{-1})<0$.  In the
matrix for $C_{r+1}(r-1)$ only its southwest entry differs from
the rank-$r$ extension of the upper-triangle formula.  Cofactor
expansion therefore gives
\begin{equation}\label{eq:matched-mixture-edge}
 C_{r+1}(r-1)=(-1)^r\delta C_r(r)\ne0.
\end{equation}
This proves the exact first-layer width.  The same minor occurs
in the denominator layer after $r+1$ steps.  Moreover,
Desnanot--Jacobi and $C_{r+1}(r)=0$ give
\[
 C_{r+2}(r-1)
 =\frac{C_{r+1}(r-1)^2}{C_r(r-1)}\ne0.
\]
Hence both layers have that width.  Symmetry handles the opposite
edge.  All depths here are uniform over growing principal sections.
\end{proof}

\begin{corollary}\label{cor:matched-mixture-ordinary-parent}
For the family in Proposition~\ref{prop:matched-mixture-exact-depth},
put
\[
 P(z):=\prod_{\nu=1}^r(1-q_\nu z),\qquad
 Q_\nu(z):=\prod_{\mu\ne\nu}(1-q_\mu z),
\]
let $\mathsf S_-^{(n-1)}$ be the lower one-step nilpotent shift,
and define
\[
 S_n:=\diag\bigl(1,P(\mathsf S_-^{(n-1)})\bigr),
 \qquad M_n:=S_nH_nS_n^T.
\]
Then the $M_n$ are leading sections of one fixed ordinary array and
are ancestors of the complete late order-two descendants of $H_n$,
at depth $n-2$ for $n\ge2$.  Their Toeplitz bulk has Laurent
symbol
\begin{equation}\label{eq:matched-mixture-ordinary-symbol}
 c(z):=\sum_{\nu=1}^r c_\nu(1-q_\nu^2)
                       Q_\nu(z)Q_\nu(z^{-1}),
\end{equation}
and every deviation from this bulk is in the northwest
$(r+1)\times(r+1)$ corner.  The bulk semibandwidth is exactly
$r-1$.  Consequently the principal determinant sequence has an
eventual constant-coefficient recurrence of order at most
$\binom{2r-2}{r-1}$.
\end{corollary}

\begin{proof}
For one exponential, direct multiplication of Laurent series gives
\[
 \sum_{s\in\mathbb Z}q^{|s|}z^s
 =\frac{1-q^2}{(1-qz)(1-qz^{-1})}
 \qquad(q<|z|<q^{-1}).
\]
Multiplying by $P(z)P(z^{-1})$ proves the bulk formula, initially
on their common annulus and hence coefficientwise.  The coefficient
of $z^{r-1}$ in \eqref{eq:matched-mixture-ordinary-symbol} is
\[
 (-1)^{r-1}\sum_{\nu=1}^r c_\nu(1-q_\nu^2)
                            \prod_{\mu\ne\nu}q_\mu,
\]
which is nonzero under the positive hypotheses.

Compared with multiplication by the full Laurent polynomials, the
finite anchored row and column combinations omit only terms whose row
or column index is at most one.  If the other index exceeds $r+1$, the
complete combination in that direction annihilates each omitted term,
since it belongs to a one-sided exponential tail.  Thus these
truncations alter only the stated corner.  The fixed first
coordinate and unit lower triangular trailing block of $S_n$
allow Proposition~\ref{prop:anchored-minor-preservation} to identify
both late order-two layers.  This proves ordinary ancestry and
$\det M_n=\det H_n$.  Once the fixed corner is passed, the ordinary
row--column transfer for semibandwidth $r-1$ has dimension
$\binom{2r-2}{r-1}$ and constant coefficients, proving the bound.
\end{proof}

\paragraph{A matched rank-two kernel.}
Fix positive $a,b$ and distinct $\rho,\tau\in(0,1)$.  The continuous
kernel
\[
 f(x,y):=a\rho^{|x-y|}+b\tau^{|x-y|}
\]
has separated rank two on each closed triangle.  For integer sampling
define
\begin{align}\label{eq:matched-exponential-parameters}
 h_s&:=a\rho^{|s|}+b\tau^{|s|},&
 H_n&:=(h_{j-i})_{i,j=1}^n,&D_n&:=\det H_n,\notag\\
 u&:=a(1-\rho^2),&v&:=b(1-\tau^2),&
 \kappa&:=ab(\rho-\tau)^2,\notag\\
 \eta&:=u\tau+v\rho,&
 \theta&:=u(1+\tau^2)+v(1+\rho^2).&&
\end{align}
Put $D_0:=1$.  For this Toeplitz family, write
$C_k(s):=\det(h_{s+j-i})_{i,j=0}^{k-1}$.

\begin{proposition}\label{prop:matched-exponential-layers}
The first three nontrivial connected-minor layers are
\begin{align}\label{eq:matched-exponential-layers}
 C_2(s)&=
 \begin{cases}
  D_2=(a+b)^2-(a\rho+b\tau)^2,&s=0,\\
  -\kappa(\rho\tau)^{|s|-1},&s\ne0,
 \end{cases}\notag\\
 C_3(s)&=
 \begin{cases}
  D_3=(u+v)((a+b)(u+v)+2\kappa),&s=0,\\
  \kappa\eta,&|s|=1,\\
  0,&|s|\ge2,
 \end{cases}\notag\\
 C_4(s)&=
 \begin{cases}
  D_4=\theta D_3-\eta^2D_2,&s=0,\\
  -\kappa\eta^2,&|s|=1,\\
  0,&|s|\ge2.
 \end{cases}
\end{align}
Uniformly over the growing principal sections, the first tridiagonal
layer occurs after exactly two condensations, and both layers are
tridiagonal after exactly three.  Neither
tridiagonal layer can be replaced by a diagonal one at these steps.
\end{proposition}

\begin{proof}
Expand $C_2(s)=h_s^2-h_{s-1}h_{s+1}$.  The two pure exponential
terms cancel, leaving the stated mixed term for $s\ne0$.
Theorem~\ref{thm:continuous-triangular-compression} gives the support
of $C_3$ and all subsequent layers.  Direct expansion gives
$C_3(0)$ and $C_3(1)$; symmetry gives $C_3(-1)$.  Equivalently,
the entry in the southwest corner of the matrix for $C_3(1)$ differs
from the rank-two upper-triangle extension by
\[
 h_1-a\rho^{-1}-b\tau^{-1}=-\frac{\eta}{\rho\tau}.
\]
Its cofactor is $C_2(2)=-\kappa\rho\tau$, so
$C_3(1)=\kappa\eta$.  Desnanot--Jacobi gives
$C_4(1)=C_3(1)^2/C_2(1)=-\kappa\eta^2$ and
$C_4(0)=(D_3^2-\kappa^2\eta^2)/D_2$; substitution gives the
displayed $D_4$.  All these equalities are polynomial identities
and therefore extend through zero denominators.  Under the stated
positive hypotheses, $\kappa,\eta>0$ and $C_2(s)\ne0$ at every
displacement.  This proves the asserted depths and surviving width.
\end{proof}

\begin{theorem}\label{thm:matched-exponential-recurrence}
With the notation above,
\begin{equation}\label{eq:matched-exponential-recurrence}
 D_n=\theta D_{n-1}-\eta^2D_{n-2}\qquad(n\ge3),
\end{equation}
where $D_1=a+b$ and $D_2=(a+b)^2-(a\rho+b\tau)^2$.
The minimal eventual constant-coefficient recurrence order is two.
\end{theorem}

\begin{proof}
Set $s:=\rho+\tau$, $p:=\rho\tau$, let $\mathsf S_-^{(n-1)}$
be the lower shift of order $n-1$, and put
\begin{equation}\label{eq:matched-exponential-parent}
 S_n:=\diag\bigl(1,I-s\mathsf S_-^{(n-1)}
                   +p(\mathsf S_-^{(n-1)})^2\bigr),
 \qquad M_n:=S_nH_nS_n^T.
\end{equation}
For $n=1$ the trailing block is empty.  The sequence $h_j$ satisfies
$h_{j+2}=sh_{j+1}-ph_j$ for $j\ge0$.  Applying these row and column combinations
shows that $M_n$ has a tridiagonal Toeplitz bulk with main-diagonal
entry $\theta$ and entries $-\eta$ on both adjacent diagonals.  Its complete exceptional northwest
block is
\begin{equation}\label{eq:matched-exponential-corner}
 M_n[1:3,1:3]=
 \begin{bmatrix}
 h_0&h_1&-ph_0\\
 h_1&h_0&h_1-sh_0\\
 -ph_0&h_1-sh_0&\theta+p^2h_0
 \end{bmatrix}\qquad(n\ge3).
\end{equation}
In particular, outside this corner every entry at distance at least
two from the main diagonal vanishes; inside it the only possible entries
at that distance are $M_{13}=M_{31}=-ph_0$.
For example, $h_2-sh_1=-ph_0$ gives $M_{13}$, while
\[
 M_{34}=(1+p)h_1-sh_0=-\eta,\qquad
 \theta=(1+s^2-p^2)h_0-2sh_1
\]
gives the two bulk coefficients and the last corner entry.

The first coordinate of $S_n$ is fixed, and its remaining block is
unit lower triangular.  Proposition~\ref{prop:anchored-minor-preservation}
therefore gives
\begin{equation}\label{eq:matched-exponential-late-block}
 \operatorname{Cond}^{n-2}(M_n,J)
 =\operatorname{Cond}^{n-2}(H_n,J)\qquad(n\ge2).
\end{equation}
Thus $M_n$ is an ordinary ancestor of the complete late order-two
block, at inverse depth $n-2$.  Unitriangularity also
gives $\det M_n=D_n$, and every entry is independent of $n$.
For $n\ge4$ an ordinary last-row continuant expansion outside the
fixed corner gives
$\det M_n=\theta\det M_{n-1}-\eta^2\det M_{n-2}$.
The formulas for $D_1,D_2,D_3$ directly verify the same identity at
$n=3$.  This proves \eqref{eq:matched-exponential-recurrence} by
the ordinary row--column method.

Finally, the initial determinants satisfy
\begin{equation}\label{eq:matched-exponential-hankel}
 D_1D_3-D_2^2=-\kappa^2\ne0.
\end{equation}
Because $\eta\ne0$, the order-two recurrence is invertible, and
its consecutive Hankel determinants are multiplied by $\eta^2$
under each forward shift.  None can vanish, so an eventual
recurrence of order one is impossible.
\end{proof}

The same anchoring gives a particularly simple ordinary parent for
the preceding rank-one continuous kernel.  When every $v_i\ne0$,
leave rows $1,2$ unchanged and simultaneously replace row $i$, for
$i\ge3$, by $H_{i,*}-(v_i/v_{i-1})H_{i-1,*}$.  Then apply the
transposed column combinations simultaneously to the resulting matrix,
leaving columns $1,2$ unchanged.  The parent consists of the northwest $2\times2$
block $H_2$ and the main-diagonal entries
\[
 m_0[i]=\frac{v_i}{v_{i-1}}\delta_{i-1}\qquad(i\ge3).
\]
It preserves both late order-two layers by the same anchoring
proposition, and its determinant is
\eqref{eq:continuous-green-determinant}.  These rational constructions
are defined when their displayed denominators are nonzero.  The
connected-minor and determinant identities extend through the excluded
specializations after denominators are cleared.

\section{Further questions}\label{sec:conclusion}

The structural route has three distinct inputs.  Ranks or functional
relations explain the vanishing connected minors; anchored
reconstruction gives the entries of a suitable ordinary parent; and
its bulk and northwest correction determine the available determinant
recurrence.  The rank-one example and the later families provide all
three pieces explicitly.  Preserving a complete late order-two block suffices for
this purpose.  Reconstructing an arbitrary earlier block is a stronger
inverse problem and is not required by these recurrence proofs.

Several questions remain within the structural route.  The general
rank bounds are sharp for independent generators, but matching
conditions or special sampled functions can lower either the depth
or the surviving width.  The continuous triangular-matching theorem
gives one systematic improvement.  A broader classification of such
improvements would distinguish genuine rank reduction from cancellations
specific to the diagonal core.  The functional counterexamples also
show why any extension beyond rational or analytic kernels must specify
its sampling domain and zero sets.

The ordinary parent often allows a determinant computation with fewer
auxiliary quantities than a direct treatment of the dense matrix.
Choosing those quantities efficiently is a separate question from
proving a finite recurrence.  The explicit
local construction, the minors used to recover the scalar recurrence,
and the pentadiagonal
transfer in Appendix~\ref{app:variable-pentadiagonal-transfer} provide
a constructive starting point.  If a minor used in that recovery vanishes, one must retain the full
transfer recurrence or prove that a smaller set of auxiliary quantities
suffices.  Direct substitution into the generic scalar coefficients
may be undefined.
Similarly, replacing variable coefficients by constant ones requires
a transfer matrix independent of position, periodic data for the entire
transfer recurrence, or another
proved finite-dimensional shift space.  Smooth, polynomial, or
exponential-polynomial entries by themselves do not imply it.

For odd polynomial distance $p=2m-1$, the repeated-root bound
$((2m^2+1)3^m+1)/2$ is attained generically at $p=1,3,5$.
Whether it is attained at every higher odd power remains open.
The boundary-companion formulation separates two possible sources
of a smaller order: coincidences among the compound eigenvalues and
vanishing boundary components.  The same distinction is relevant
to other recurrent-tail symbols.

Finally, the rank identities for products and structured changes of
basis provide sufficient condensation bounds in the given coordinate
order.  They do not classify every dense representation of a tractable
determinant sequence.  When the factor determinants already have
transfer matrices independent of position, their tensor product supplies a recurrence for
the finite matrix product.  Understanding when the product itself
admits a narrower banded descendant or an ordinary parent with a
smaller northwest correction is an additional structural problem.

\section*{Acknowledgements}

The authors acknowledge the use of ChatGPT (OpenAI) for assistance during
the preparation of this manuscript. The authors reviewed and verified the
material and take full responsibility for its content.

\appendix
\clearpage
\section{Interleaved layers for hand calculation}
\label{app:interleaved-layers}

The two layers of a Dodgson block can be displayed in one interleaved
array.  Place each entry of the denominator layer at the center of the
four entries of the first layer to which it applies, and print that
denominator entry in bold.  Double square brackets distinguish this
display from an ordinary matrix.  Thus, when $\beta\ne0$,
\[
 \Ddet\left[\!\left[
 \begin{array}{ccc}
 a&&b\\
 &\boldsymbol{\beta}&\\
 c&&d
 \end{array}
 \right]\!\right]
 =\frac{ad-bc}{\beta}.
\]
The brackets encode the same pair $(A,B)$ used in the main text;
interleaving changes only its presentation.  In a larger array every
bold entry is the divisor for the four first-layer entries surrounding
it.  After one step, the internal entries of the former first layer
become the new bold denominator layer.

Here is a calculation with two zero entries in the original matrix.
Replace both zeros by the same indeterminate $\eps$ and use the
all-ones denominator layer.  The starting block is
\begin{equation}\label{eq:interleaved-initial}
 \mathcal U^{(0)}_\eps:=
 \left[\!\left[
 \begin{array}{rrrrrrrrr}
 4&&3&&9&&1&&2\\
 &\boldsymbol{1}&&\boldsymbol{1}&&\boldsymbol{1}&&\boldsymbol{1}&\\
 1&&5&&7&&8&&3\\
 &\boldsymbol{1}&&\boldsymbol{1}&&\boldsymbol{1}&&\boldsymbol{1}&\\
 3&&\eps&&\eps&&4&&6\\
 &\boldsymbol{1}&&\boldsymbol{1}&&\boldsymbol{1}&&\boldsymbol{1}&\\
 6&&1&&2&&5&&7\\
 &\boldsymbol{1}&&\boldsymbol{1}&&\boldsymbol{1}&&\boldsymbol{1}&\\
 8&&1&&1&&7&&2
 \end{array}
 \right]\!\right].
\end{equation}
Write $H_\eps$ for its first layer.  Ordinary condensation gives
$\Ddet\mathcal U^{(0)}_\eps=\det H_\eps$.  Write
$\mathcal U^{(k)}_\eps:=\operatorname{Cond}^k(\mathcal U^{(0)}_\eps)$.
The first descendant is
\begin{equation}\label{eq:interleaved-first}
 \mathcal U^{(1)}_\eps=
 \left[\!\left[
 \begin{array}{rrrrrrr}
 17&&-24&&65&&-13\\
 &\boldsymbol{5}&&\boldsymbol{7}&&\boldsymbol{8}&\\
 -15+\eps&&-2\eps&&28-8\eps&&36\\
 &\boldsymbol{\eps}&&\boldsymbol{\eps}&&\boldsymbol{4}&\\
 3-6\eps&&\eps&&-8+5\eps&&-2\\
 &\boldsymbol{1}&&\boldsymbol{2}&&\boldsymbol{5}&\\
 -2&&-1&&9&&-39
 \end{array}
 \right]\!\right].
\end{equation}
For example, its first-layer entry in row two and column one is
$1\cdot\eps-5\cdot3=-15+\eps$.  At the next step some divisors
are $\eps$, so specialization must wait until the quotients have been
simplified in the rational function field.

\clearpage
The next two descendants, with every quotient simplified, are
\begin{equation}\label{eq:interleaved-second}
 \mathcal U^{(2)}_\eps=
 \left[\!\left[
 \begin{array}{rrrrr}
 -72-2\eps&&-96+46\eps&&338-13\eps\\
 &\boldsymbol{-2\eps}&&\boldsymbol{28-8\eps}&\\
 -9-11\eps&&-12-2\eps&&58-41\eps\\
 &\boldsymbol{\eps}&&\boldsymbol{-8+5\eps}&\\
 -3+8\eps&&-4+7\eps&&66-39\eps
 \end{array}
 \right]\!\right],
\end{equation}
\begin{equation}\label{eq:interleaved-third}
 \mathcal U^{(3)}_\eps=
 \left[\!\left[
 \begin{array}{rrr}
 237-255\eps&&-54+239\eps\\
 &\boldsymbol{-12-2\eps}&\\
 71-61\eps&&70+73\eps
 \end{array}
 \right]\!\right].
\end{equation}
The cancellation in the first-layer entry in row two and column one of
\eqref{eq:interleaved-second}, for instance, is
\[
 \frac{(-15+\eps)\eps-(-2\eps)(3-6\eps)}{\eps}
 =-9-11\eps.
\]
This gives a finite polynomial value at $\eps=0$, although the
uncancelled quotient is then $0/0$.

The last step gives
\begin{align}\label{eq:interleaved-terminal}
 \Ddet\mathcal U^{(3)}_\eps
 &=\frac{(237-255\eps)(70+73\eps)
       -(-54+239\eps)(71-61\eps)}{-12-2\eps}\notag\\
 &=2018\eps-1702.
\end{align}
Consequently,
\begin{align}\label{eq:interleaved-specialization}
 \det H_0
 &=\lim_{\eps\to0}\Ddet\mathcal U^{(0)}_\eps
  =\Ddet\left[\!\left[
    \begin{array}{rrr}
     237&&-54\\
     &\boldsymbol{-12}&\\
     71&&70
    \end{array}
    \right]\!\right]\notag\\
 &=\frac{237\cdot70-(-54)\cdot71}{-12}=-1702.
\end{align}
The same perturbation was propagated from the original matrix through
all four steps.  In particular, the bold entries in later displays
were inherited from that calculation, not independently perturbed.
The polynomial connected-minor identities justify the specialization.
The interleaved display makes the required divisors visible throughout;
the separate-layer notation in the main text remains equivalent.

\clearpage
\section{Sharpness and structured transformations}
\label{app:some-aspects}

We prove the sharpness and regular converse of the compression bounds,
then extend them through products and structured changes of basis.
The final result gives compact coefficients for geometric rank-two
generators.

\subsection{Sharpness and a regular converse}

\begin{proposition}\label{prop:band-core-generic-sharpness}
For the separated model \eqref{eq:band-core-separated-data} with
independent generator and band entries, the depth and width bounds in
Theorem~\ref{thm:band-core-compression} are generically sharp, provided
the relevant minors fit in the finite matrix.  More precisely, for
$1\le t\le r_+$ every connected size-$t$ minor lying wholly beyond
the upper band is a nonzero polynomial in the generator data.  For
$t>r_+$ the extreme supported minor
$C_t(i,i+w_+)$ is also a nonzero polynomial.  The corresponding
lower statements hold by transposition.
\end{proposition}

\begin{proof}
For the first assertion, specialize all but $t$ generator components
to zero and choose the remaining row and column evaluations to be
Vandermonde matrices.  The resulting size-$t$ determinant is nonzero,
so its polynomial is not identically zero.

For the second assertion write $r:=r_+$ and $b:=b_+$.  Regard the
outer band entries $d_z:=H_{z,z+b}$ as independent variables.
In the block defining $C_t(i,i+b+r)$, the variable $d_z$ occurs
for exactly the rows $z=i+r,\ldots,i+t-1$.  Selecting all these
entries leaves the first $r$ rows and last $r$ columns.  Thus
\begin{equation}\label{eq:band-core-outer-coefficient}
\begin{split}
 \coef{\prod_{z=i+r}^{i+t-1}d_z}C_t(i,i+b+r)
 ={}&(-1)^{r(t-r)}
     \det F^+[i:i+r-1,:]\\
 &\quad{}\times
     \det G^+[i+b+t:i+b+t+r-1,:].
\end{split}
\end{equation}
The empty determinants for $r=0$ equal one.  The right-hand side
is nonzero for generic generators.  This proves exact upper
semibandwidth $b+r$; the lower edge is analogous.

On a side whose rank is $R$, all layers of sizes at most $R$
have generically nonzero connected minors at arbitrarily large
offsets in growing sections.  No smaller depth therefore guarantees
a fixed first-layer band.  At depth $R$ the denominator still
contains size-$R$ minors, so it is generically not banded.
For a fixed finite size the finitely many claimed nonvanishings
hold on a nonempty Zariski-open set.
\end{proof}

Generic sharpness is a statement about the full class.  Individual
matrices can have smaller ranks, narrower bands, or an earlier
collapse.  In particular, if the entire matrix has rank $r$,
then every size-$(r+1)$ minor is zero: the depth-$r$ first layer
collapses completely.  A diagonal or band perturbation is what
usually turns this global low-rank collapse into a useful
nontrivial banded descendant.

\subsubsection*{A converse with nonzero preceding minors}

Vanishing connected minors are the exact support test for a specified
condensation layer.  They do not, by themselves, characterize the
original rectangle ranks.  The following local rank statement
specifies a useful regularity condition under which a converse is valid.

\begin{lemma}\label{lem:local-connected-rank}
Let $M$ be a rectangular matrix over a field, and let $r\ge1$.
If every contiguous $(r+1)\times(r+1)$ minor is zero and every
contiguous $r\times r$ minor is nonzero, then
$\operatorname{rank}M=r$, provided both dimensions are at least $r$.
\end{lemma}

\begin{proof}
The first $r$ columns are independent because their first $r$ rows
form a nonsingular minor.  Induct on the subsequent columns.
For a new column $j>r$, solve for its linear combination of columns
$j-r,\ldots,j-1$ on the first $r$ rows.  Suppose this same relation
has been verified through row $a+r-1$.  In the contiguous block
on rows $a,\ldots,a+r$ and columns $j-r,\ldots,j$, subtract that
combination from column $j$.  Its first $r$ entries are zero.
The zero $(r+1)$-minor and the nonzero preceding $r$-minor force
the last entry to vanish as well.  Induction extends the relation
through every row.  Every column is consequently in the span of
the first $r$, which proves the result.  If a dimension is exactly
$r$, the conclusion is immediate from the nonzero $r$-minor.
\end{proof}

\begin{corollary}\label{cor:regular-band-core-converse}
Fix $b\ge0$ and $r\ge1$.  Suppose every valid upper connected
size-$r$ minor lying wholly beyond the diagonal band is nonzero:
\[
 C_r(i,j)\ne0\qquad(j-i\ge b+r).
\]
Then the following conditions are equivalent:
\[
 \mathfrak r_+(b;H)\le r,
 \qquad
 C_{r+1}(i,j)=0\quad(j-i>b+r).
\]
The reflected statement holds for the lower part.
\end{corollary}

\begin{proof}
The forward implication is Theorem~\ref{thm:band-core-compression}.
For the converse, apply Lemma~\ref{lem:local-connected-rank}
to each rectangle $H[1:c,c+b+1:N]$ having both dimensions at least
$r$.  Every contiguous $r$-minor of this rectangle satisfies the
stated nonvanishing condition, and every contiguous $(r+1)$-minor
satisfies the stated vanishing condition.  Rectangles with a
smaller dimension have rank at most $r$ automatically.
\end{proof}

The nonvanishing hypothesis cannot be omitted.  For any $q\ge1$,
take an order-$4q$ matrix whose only nonzero entries are
\[
 H_{2a,\,2q+2a}:=1\qquad(1\le a\le q).
\]
Every pair of consecutive rows contains a zero row, so every
connected $2\times2$ minor vanishes.  Nevertheless
$H[1:2q,2q+1:4q]$ has rank $q$.  Thus even a completely zero
first condensed layer does not imply a uniform off-diagonal rank
bound without a regularity assumption.

\subsection{Products, inverses, and structured similarities}
\label{sec:structured-closure}

Classical identities for quasiseparable ranks
\cite{EidelmanGohbergHaimovici2014} extend the support theorem to products
and suitably structured changes of basis.  We record the elementary
rank arguments and the resulting condensation depths.  For a square
matrix $H$, put
\[
 q_\pm(H):=\mathfrak r_\pm(0;H).
\]
Splitting off the first $b_+$ columns of an upper cut rectangle, or the
first $b_-$ rows of a lower cut rectangle, gives
\begin{equation}\label{eq:closure-cut-bound}
 q_\pm(H)\le b_\pm+\mathfrak r_\pm(b_\pm;H).
\end{equation}

\begin{lemma}\label{lem:closure-cut-ranks}
For square matrices $X,Y$ of the same order and an invertible matrix $S$,
\begin{equation}\label{eq:closure-cut-algebra}
 q_\pm(XY)\le q_\pm(X)+q_\pm(Y),
 \qquad q_\pm(S^{-1})=q_\pm(S).
\end{equation}
In the inverse identity the two ranks agree at every individual cut.
\end{lemma}

\begin{proof}
Partition at a cut $c$.  The upper-right block of $XY$ is
$X_{11}Y_{12}+X_{12}Y_{22}$, whose rank is at most
$\operatorname{rank}Y_{12}+\operatorname{rank}X_{12}$; the lower-left
block is analogous.

Write $T:=S^{-1}$.  If $v\in\ker S_{12}$, then
$S(0,v)=(0,S_{22}v)$, and applying $T$ shows that
$S_{22}v\in\ker T_{12}$ and $T_{22}S_{22}v=v$.
Interchanging $S,T$ proves that the restrictions of $S_{22}$ and
$T_{22}$ are mutually inverse maps between these two kernels.
Their dimensions agree, so $S_{12}$ and $T_{12}$ have the same rank.
The same argument with vectors $(v,0)$ proves equality for the
lower-left blocks.  Taking maxima over $c$ proves the lemma.
\end{proof}

Retaining the band core gives a sharper product bound than first
replacing every matrix by its cut-rank bound in
\eqref{eq:closure-cut-bound}.

\begin{proposition}\label{prop:closure-band-product}
Let $X,Y$ be square matrices of the same order, and suppose
\[
 \mathfrak r_\pm(b_\pm^X;X)\le r_\pm^X,
 \qquad
 \mathfrak r_\pm(b_\pm^Y;Y)\le r_\pm^Y,
\]
where the band widths and rank bounds are nonnegative integers.  Then
\begin{equation}\label{eq:closure-band-product}
 \mathfrak r_\pm(b_\pm^X+b_\pm^Y;XY)
 \le r_\pm^X+r_\pm^Y.
\end{equation}
\end{proposition}

\begin{proof}
Write $n$ for the common order.  For an upper rectangle put
$I:=\{1,\ldots,c\}$,
$J:=\{c+b_+^X+b_+^Y+1,\ldots,n\}$, and
$K:=\{1,\ldots,c+b_+^X\}$; an empty rectangle needs no argument.
With $K^c$ the complement in $\{1,\ldots,n\}$,
\[
 (XY)[I,J]=X[I,K]Y[K,J]+X[I,K^c]Y[K^c,J].
\]
The first summand has rank at most $r_+^Y$, using the cut
$c+b_+^X$ in $Y$; the second has rank at most $r_+^X$, using
the cut $c$ in $X$.  Transposition, with the factor order reversed,
proves the lower assertion.
\end{proof}

\begin{corollary}\label{cor:closure-product-depth}
Let
$H_n:=X_n^{(1)}\cdots X_n^{(k)}$, where $k$ is fixed and
\[
 \mathfrak r_\pm(b_\pm^{(j)};X_n^{(j)})\le r_\pm^{(j)}
 \qquad(1\le j\le k)
\]
with bounds independent of $n$.  Set
\[
\begin{gathered}
 b_\pm:=\sum_{j=1}^k b_\pm^{(j)},\qquad
 r_\pm:=\sum_{j=1}^k r_\pm^{(j)},\\
 w_\pm:=b_\pm+r_\pm,\qquad R:=\max(r_-,r_+).
\end{gathered}
\]
The first layer of $H_n$ is $(w_-,w_+)$-banded after $R$
ordinary condensation steps, and both layers are banded after $R+1$
steps, whenever the indicated layers are present.  If $H_n$ is
invertible, $H_n^{-1}$ has the same width bound after
$\max(w_-,w_+)$ steps for its first layer and one further step for
both layers, whenever these layers are present.
\end{corollary}

\begin{proof}
Iterate Proposition~\ref{prop:closure-band-product} and apply
Theorem~\ref{thm:band-core-compression}.  For the inverse,
\eqref{eq:closure-cut-bound} and Lemma~\ref{lem:closure-cut-ranks}
give $q_\pm(H_n^{-1})\le w_\pm$; apply the same theorem with
core widths zero.  These are connected-minor statements, including
regularized specializations at zero intermediate denominators.
\end{proof}

\begin{proposition}\label{prop:structured-similarity-depth}
Suppose
$\mathfrak r_\pm(b_\pm;H_n)\le r_\pm$ and
$q_\pm(S_n)\le s_\pm$, with all bounds independent of $n$ and
each $S_n$ invertible.  For
$\widetilde H_n:=S_n^{-1}H_nS_n$,
\begin{equation}\label{eq:structured-similarity-ranks}
 \mathfrak r_\pm(b_\pm;\widetilde H_n)\le r_\pm+2s_\pm.
\end{equation}
Thus its first layer is banded after
$\widetilde R:=\max(r_-+2s_-,r_++2s_+)$ steps and both layers
after $\widetilde R+1$ steps, with semibandwidths
$\widetilde w_\pm:=b_\pm+r_\pm+2s_\pm$, whenever these layers are
present.
\end{proposition}

\begin{proof}
By Lemma~\ref{lem:closure-cut-ranks}, both $S_n$ and $S_n^{-1}$
have cut ranks at most $s_\pm$.  Apply
Proposition~\ref{prop:closure-band-product} twice, assigning core
widths zero to these two factors and retaining $b_\pm$ for $H_n$.
Theorem~\ref{thm:band-core-compression} gives the stated depths.
\end{proof}

For instance, an invertible tridiagonal $S_n$ has $s_\pm\le1$,
even though its inverse can be dense.  Conjugating an ordinary
$(b_-,b_+)$-banded family therefore gives a first layer with widths
at most $(b_-+2,b_++2)$ after two steps, and both layers after three.
The rank bounds certify compression directly in the transformed
ordering.  Constructing an ordinary parent and obtaining a constant-coefficient
recurrence require the additional hypotheses stated in the main text.

An arbitrary reordering can instead destroy every fixed-depth bound.

\begin{example}\label{ex:order-permutation}
Let $A_{2m}:=\operatorname{tridiag}(1,3,1)$, and let $P_{2m}$
place the odd indices first and the even indices second.  Put
$B_{2m}:=P_{2m}A_{2m}P_{2m}^T$.  Its upper-right $m\times m$
block is lower bidiagonal with main diagonal and subdiagonal entries $1$,
because its $(i,j)$ entry is $(A_{2m})_{2i-1,2j}$.  Denote its
connected minors by $C_t^{B_{2m}}(i,j)$.  They satisfy
\begin{equation}\label{eq:order-permutation-minors}
 C_t^{B_{2m}}(1,m+1)=1\qquad(1\le t\le m).
\end{equation}
At any fixed depth $k$, choose $m\ge k+1$.  Taking $t=k+1$
leaves a nonzero first-layer entry at displacement $m$.
No fixed depth therefore produces a width bounded independently
of $m$, although $A_{2m}$ has two banded layers after one step and
$\det B_{2m}=\det A_{2m}$.  The cut ranks of these permutations
grow with $m$, so Proposition~\ref{prop:structured-similarity-depth}
does not give a uniform bound for this ordering.
\end{example}

Some structural subclasses survive arbitrary simultaneous permutations.
For $H=\diag(\delta)+UV^T$, where $U,V$ have $r$ columns, such a
permutation only reorders the entries of $\delta$ and the rows of $U,V$.
The cut ranks remain at most $r$, so the corresponding compression
bound still applies.  If the diagonal part is $\gamma I$ for a scalar
$\gamma$, every invertible $S$ preserves this class under similarity:
\[
 S^{-1}(\gamma I+UV^T)S
 =\gamma I+(S^{-1}U)(V^TS).
\]
The counterexample therefore concerns general ordered structure, not
every low-rank model.

\paragraph{Determinant transport across sizes.}
For size-by-size similarities,
$\det(S_n^{-1}H_nS_n)=\det H_n$, so every already established
scalar recurrence transports unchanged.  For row and column
permutations, put $K_n:=P_nH_nQ_n$ and
$\sigma_n:=\det P_n\det Q_n$.  If
$D_n:=\det H_n$ satisfies
$D_n=\sum_{j=1}^d a_j(n)D_{n-j}$, then
$E_n:=\det K_n=\sigma_nD_n$ satisfies
\begin{equation}\label{eq:closure-sign-transport}
 E_n=\sum_{j=1}^d a_j(n)\frac{\sigma_n}{\sigma_{n-j}}E_{n-j}.
\end{equation}
Arbitrary size-dependent signs can make these coefficients variable;
removing the known sign recovers the original sequence.  Likewise,
the support bounds above alone leave the transfer position dependent.
A constant-coefficient conclusion requires a finite transfer matrix
independent of position, or a separately established scalar recurrence.

These identities concern transformations of each finite $H_n$.
If an ambient or infinite array is transformed first and its leading
sections are then taken, those sections need not be the corresponding
transforms of the original $H_n$.  Their determinant sequence requires
a separate identification; invariance of the full determinant does
not identify all principal sections.

\subsection{A compact six-term formula for geometric generators}

\begin{proposition}\label{prop:geometric-ranktwo-six}
Suppose
\[
 \begin{aligned}
 f_s(1)&:=F_s,& f_s(i+1)&:=\varphi_s f_s(i),\\
 g_s(1)&:=G_s,& g_s(i+1)&:=\chi_s g_s(i)
 \qquad(s=1,2).
 \end{aligned}
\]
and define
\[
 D_n:=\det\left(\sum_{s=1}^2f_s(i)g_s(j)
                         -x\mathbf1_{i=j}\right)_{i,j=1}^n,
 \qquad D_0:=1.
\]
Put
\[
 \rho:=\varphi_1\varphi_2\chi_1\chi_2,\qquad
 h_0:=1+\rho,
\]
\[
 h_1:=\varphi_1\chi_1+\varphi_2\chi_2,\qquad
 h_2:=\varphi_1\chi_2+\varphi_2\chi_1,
\]
and
\begin{align*}
 r_1&:=h_0+h_1+h_2,\\
 r_2&:=h_0h_1+h_0h_2+h_1h_2+3\rho,\\
 r_3&:=h_0h_1h_2+2\rho r_1,\\
 r_4&:=\rho r_2,\qquad
 r_5:=\rho^2r_1,\qquad r_6:=\rho^3.
\end{align*}
Then
\begin{equation}\label{eq:geometric-ranktwo-six}
 D_n=-\sum_{j=1}^6x^jr_jD_{n-j},\qquad n\ge6.
\end{equation}
Repeated or zero geometric ratios are allowed; the order need not be minimal.
\end{proposition}
\begin{proof}
The coefficient identity is
\[
 1+\sum_{j=1}^6r_jt^j
   =\prod_{a=0}^2(1+h_at+\rho t^2)
   =(1+t)(1+\rho t)
       \prod_{s=1}^2\prod_{u=1}^2(1+\varphi_s\chi_u t).
\]
For generic nonzero parameters, expand the determinant through the
two generators.  Each entry of the $2\times2$ matrix
$\bigl(\sum_{i=1}^n g_s(i)f_t(i)\bigr)_{s,t=1}^2$ is a geometric sum.  Its trace and determinant show that $D_n$ is a linear
combination of the six exponential sequences with bases
\[
 -x,\quad -x\varphi_1\chi_1,\quad -x\varphi_1\chi_2,\quad
 -x\varphi_2\chi_1,\quad -x\varphi_2\chi_2,\quad -x\rho.
\]
Their characteristic polynomial is
\[
 z^6+xr_1z^5+x^2r_2z^4+x^3r_3z^3
                 +x^4r_4z^2+x^5r_5z+x^6r_6,
\]
which gives \eqref{eq:geometric-ranktwo-six}.
For each fixed $n$, the identity is polynomial in the original
parameters.  It therefore extends to the exceptional parameter values.
\end{proof}

\needspace{11\baselineskip}
The constancy of the geometric ratios is essential in this proposition.
For example, take $f_1(i)=g_1(i)=1$, $f_2(i)=g_2(i)=i$, and $x=1$.
The first seven determinants are
\[
 (D_0,\ldots,D_6)=(1,1,-5,10,-13,9,9).
\]
Substituting only the current ratios
$\varphi_1=\chi_1=1$ and $\varphi_2=\chi_2=6/5$ into the coefficients
above gives
\[
 D_6+\sum_{j=1}^6 r_jD_{6-j}=\frac{16298}{15625}\ne0.
\]
Thus the same formulas cannot be used pointwise for arbitrary functions.
For position-dependent data the valid observation formula
\eqref{eq:local-vapor-scalar-coefficients} uses products of successive
transfer matrices and their appropriate shifts.

\clearpage
\section{Scalar determinant realizations for comparison}
\label{app:scalar-comparisons}

The calculations below give alternative determinant evaluations and, for the
unanchored filters, smaller boundary corrections.  The main proofs
additionally identify the common connected-minor descendant.

\subsection{The rank-one scalar comparison}\label{app:rankone-comparison}

Use the notation of Section~\ref{sec:rankone-model}.
For comparison, the determinant lemma gives
\begin{equation}\label{eq:rankone-comparison-formula}
 D_n=\left(\prod_{i=1}^n\delta_i\right)
      \left(1+\sum_{i=1}^n\frac{q_i}{\delta_i}\right),
\end{equation}
initially when the $\delta_i$ are nonzero, and then after clearing
denominators.  Eliminating the partial sum gives the same recurrence.
This familiar calculation verifies the scalar answer.  The explicit
parent construction also identifies an ordinary ancestor of the
order-two descendant.
Arbitrary position-dependent $f_i,g_i,h_i$ do not imply constant
coefficients.  For instance, the partial sums in the rank-one
characteristic polynomial can encode an arbitrary sequence.  Constant coefficients in Section~\ref{sec:geometrictails} come from
its Toeplitz structure.

\subsection{A cosine determinant comparison}\label{sec:cosine-comparison}
Use the notation of Section~\ref{sec:cosine-vaporisation}.  The recurrence
and generic minimality in Theorem~\ref{thm:cosine-recurrence} were obtained
from the explicit ordinary ancestor.  A separate low-rank calculation makes
the boundary coefficients and exceptional parameter values explicit.
For $x\sigma\ne0$ it gives, for every $n\ge0$,
\begin{equation}\label{eq:cosine-determinant-formula}
 D_n=x^{n-2}\left[
 x^2+x\frac{\sin(nC)\cos((n+1)C)}{\sin C}
 -\frac14\left(n^2-\frac{\sin^2(nC)}{\sin^2 C}\right)
 \right].
\end{equation}
Indeed, write $K_n=U_n\diag(1,-1)U_n^T$, where row $i$ of $U_n$ is
$(\cos(Ci),\sin(Ci))$.  Sylvester's determinant identity reduces
the calculation to a determinant of order two.  Its linear coefficient
is $\sum_{i=1}^n\cos(2Ci)$, while its constant coefficient is
\[
 -\det(U_n^TU_n)
 =-\frac14\left(n^2-
       \left|\sum_{j=1}^ne^{2\mathrm i Cj}\right|^2\right)
 =-\frac14\left(n^2-\frac{\sin^2(nC)}{\sin^2 C}\right)
\]
for real $C$; the corresponding algebraic identity uses the product
of the two reciprocal geometric sums.  This proves
\eqref{eq:cosine-determinant-formula} independently of the ancestor.

For a finite exponential expansion, put $\rho:=e^{2\mathrm i C}\ne1$ and define
\begin{align}\label{eq:cosine-mode-coefficients}
 a_0&:=x^2-\frac{x}{2}-\frac{\rho}{2(\rho-1)^2},\notag\\
 a_+&:=\frac{\rho\bigl(2x(\rho-1)+1\bigr)}{4(\rho-1)^2},
 &a_-&:=\frac{\rho-2x(\rho-1)}{4(\rho-1)^2}.
\end{align}
Then
\begin{equation}\label{eq:cosine-mode-expansion}
 D_n=x^{n-2}\left(-\frac{n^2}{4}+a_0+a_+\rho^n+a_-\rho^{-n}\right).
\end{equation}
It verifies the same order-five annihilator and also shows that the
backward recurrence \eqref{eq:cosine-backward-recurrence} already
holds for $N\ge5$.  By polynomial specialization this earlier
threshold is valid for every $x,c$.  For $\rho\ne\pm1$ and
$a_+a_-\ne0$, the three distinct modes have polynomial degrees
$2,0,0$, so the minimal order is exactly five.  In particular this
holds for all real $x\ne0$ and real $C$ with
$\sin C\cos C\ne0$.

\subsection{Unanchored recurrent-tail filters}

Use the Toeplitz symbol and tail hypotheses
\eqref{eq:tailannpos}--\eqref{eq:tailannneg}.
Define
\begin{equation}\label{eq:tailfilters}
 L_n:=\sum_{a=0}^{r_-}\alpha_a(\mathsf S_-^{(n)})^a,
 \qquad
 R_n:=\sum_{q=0}^{r_+}\beta_q(\mathsf S_+^{(n)})^q
\end{equation}
and
\begin{equation}\label{eq:tailfilteredmatrix}
 C_n:=L_nT_n(b)R_n.
\end{equation}
Finally put
\begin{equation}\label{eq:tailfilteredsymbol}
 c_s:=\sum_{a=0}^{r_-}\sum_{q=0}^{r_+}
 \alpha_a\beta_q b_{s+a-q}.
\end{equation}

\begin{theorem}
\label{thm:tailannihilator}
For every $n$ the matrices in \eqref{eq:tailfilters} are unitriangular and
\begin{equation}\label{eq:taildetpreserve}
 \det C_n=\det T_n(b).
\end{equation}
Moreover, in one-based indexing,
\begin{equation}\label{eq:tailbulkentry}
 (C_n)_{ij}=c_{j-i}
 \qquad(i>r_-,\ j>r_+),
\end{equation}
and
\begin{equation}\label{eq:tailbulkband}
 c_s=0\qquad(s>m_+\ \text{or}\ s<-m_-).
\end{equation}
Consequently, if $T_n(c)$ denotes the $(m_-,m_+)$-banded Toeplitz matrix
with symbol $c$, then
\begin{equation}\label{eq:taildefect}
 E_n:=C_n-T_n(c)
\end{equation}
is supported in the union of the first $r_-$ rows and the first $r_+$
columns.  In particular,
\begin{equation}\label{eq:taildefectrank}
 \operatorname{rank}E_n\le r_-+r_+,
\end{equation}
uniformly in $n$.  In fact the correction is contained in a fixed rectangle:
\begin{equation}\label{eq:tail-corner-size}
 E_{ij}=0\quad\text{if }i>H_-\text{ or }j>H_+,
 \qquad
 H_-:=\max(r_-,r_++m_-),\quad H_+:=\max(r_+,r_-+m_+).
\end{equation}
Its entries are explicitly
\begin{equation}\label{eq:tail-boundary-entry}
 E_{ij}=
 \sum_{a=0}^{\min(r_-,i-1)}\sum_{q=0}^{\min(r_+,j-1)}
 \alpha_a\beta_q b_{j-i+a-q}-c_{j-i}.
\end{equation}
If $r_+\ge m_+$ and $r_-\ge m_-$, the support is contained in the
smaller northwest $r_-\times r_+$ corner.
\end{theorem}

\begin{proof}
The filters are unitriangular, so they preserve the determinant.
Multiplication gives the truncated double sum in
\eqref{eq:tail-boundary-entry}; when both sums have their full ranges it is
$c_{j-i}$.  For $s>m_+$, apply the positive-tail recurrence at
$s+a$; for $s<-m_-$, apply the negative-tail recurrence at
$s-q$.  Equations~\eqref{eq:tailannpos}--\eqref{eq:tailannneg}
give \eqref{eq:tailbulkband}.

Only the first $r_-$ rows or first $r_+$ columns can have truncated
sums, proving the stated support union and rank bound.  In a top
row $i\le r_-$ with $j>H_+$, the full positive-tail recurrence is applied
and $j-i>m_+$, so both the entry and its bulk value vanish.
Transposition treats a left column with $i>H_-$, proving
\eqref{eq:tail-corner-size}.

If $r_\pm\ge m_\pm$, a row-sum term missing at $i\le r_-$
has $a\ge i$.  For $j>r_+$, its complete positive-tail sum is
taken at $j-i+a\ge j>m_+$ and vanishes.  The reflected argument
treats missing column terms.  Thus only the northwest
$r_-\times r_+$ rectangle can differ from the bulk.
\end{proof}

\begin{corollary}
\label{cor:rationaltaildet}
Under the hypotheses of Theorem~\ref{thm:tailannihilator}, the determinant
sequence
\[
 D_n:=\det T_n(b)
\]
is eventually annihilated by a homogeneous constant-coefficient recurrence
of order at most
\begin{equation}\label{eq:tailrowcolbound}
 \binom{m_-+m_+}{m_-}.
\end{equation}
In particular this applies whenever the positive and negative symbol tails
have rational generating functions after deletion of finite prefixes.

If, in addition, the original Toeplitz family is symmetric, $m_-=m_+=m$,
and the two filters can be chosen as transposes of one another, then
$C_n=L_nT_n(b)L_n^T$ is symmetric and the bound may be sharpened to
\begin{equation}\label{eq:tailsymbound}
 \frac{3^m+1}{2}.
\end{equation}
\end{corollary}

\begin{proof}
The displayed entries of $C_n$ are independent of $n$, and its correction
has the fixed rectangular support \eqref{eq:tail-corner-size}.
The ordinary boundary transfer of
Lemma~\ref{lem:exact-boundary-companion} therefore proves the binomial
bound.  In the symmetric case, the reciprocal-root argument in
Corollary~\ref{cor:anchored-tail-recurrence} applies to the same
fixed symmetric bulk and every boundary state, giving
\eqref{eq:tailsymbound}, including singular specializations.
Rational tail generating functions supply precisely the eventual
recurrences \eqref{eq:tailannpos}--\eqref{eq:tailannneg}.
\end{proof}

\begin{remark}\label{rem:detbanding-vs-vaporisation}
The unanchored filters in this appendix preserve the ordinary determinant.
The anchored construction of Theorem~\ref{thm:anchored-recurrent-tails}
also preserves both layers of a specified late descendant.  The two constructions are used with their stated boundary normalizations.
\end{remark}

\paragraph{Comparison for geometric tails.}
The determinant-preserving calculation for the same family produces a
second ordinary tridiagonal representative, with a correction confined to its
first main diagonal entry.  We retain the calculation to compare its scalar
equivalence with the exact descendant equality
\eqref{eq:geometric-corner-cond}.

\begin{proof}[Independent verification of recurrence~\eqref{eq:recB}]
Put
\[
 \rho_+:=\frac{b_2}{b_1},\qquad
 \rho_-:=\frac{b_{-2}}{b_{-1}},
\]
and use the first-order tail filters
\begin{equation}\label{eq:geotailfilters}
 L_n=I-\rho_-\mathsf S_-^{(n)},\qquad
 R_n=I-\rho_+\mathsf S_+^{(n)}.
\end{equation}
Theorem~\ref{thm:tailannihilator} shows that
$C_n:=L_nB_nR_n$ has the same determinant as $B_n$.  Here the boundary
cancellations are stronger than the general theorem requires: $C_n$ is
tridiagonal.  In signed-diagonal notation its main diagonal is
\begin{align}\label{eq:geotailalpha}
 c_0[1]&:=b_0,\notag\\
 c_0[i]&:=
 \frac{b_{-1}b_0b_1-b_{-1}^2b_2-b_{-2}b_1^2+b_{-2}b_0b_2}
      {b_{-1}b_1}\qquad(i\ge2),
\end{align}
and its superdiagonal and subdiagonal are constant:
\begin{equation}\label{eq:geotailul}
 c_1[i]:=\frac{b_1^2-b_0b_2}{b_1},\qquad
 c_{-1}[i]:=\frac{b_{-1}^2-b_0b_{-2}}{b_{-1}}.
\end{equation}
Expansion of the final row of this ordinary tridiagonal determinant gives
\[
 D_n=c_0[n]D_{n-1}-c_{-1}[n-1]c_1[n-1]D_{n-2}\qquad(n\ge2),
\]
which is exactly \eqref{eq:recB}.
\end{proof}

\subsection{Unanchored polynomial differences}

Let $A_n(x,p)$, $r=p+1$, and $b$ be as in the polynomial-distance
application.  The unanchored representative has the following explicit
smaller northwest correction.

\begin{proof}[Scalar-filter proof of Corollary~\ref{cor:polydistance}]
Put $\alpha_a:=(-1)^a\binom ra$ for $0\le a\le r$, and use the transpose
pair of filters
\begin{equation}\label{eq:polydistancefilter}
 L_n:=(I-\mathsf S_-^{(n)})^r,\qquad R_n:=L_n^T.
\end{equation}
They have determinant one.  Thus the symmetric family of leading sections
\begin{equation}\label{eq:polyexplicit-representative}
 C_n:=L_nA_n(x,p)L_n^T=T_n(c)+E_n
\end{equation}
has determinant $\det A_n(x,p)$.  Its Toeplitz bulk is explicit:
\begin{equation}\label{eq:polyexplicit-bulk}
 c_s:=\sum_{h=-r}^{r}(-1)^h\binom{2r}{r+h}b_{s+h}.
\end{equation}
Indeed, collecting the terms with $a-q=h$ in
$\sum_{a,q=0}^r\alpha_a\alpha_qb_{s+a-q}$ gives
\eqref{eq:polyexplicit-bulk} by the binomial convolution identity.  For
$|s|>r$, every argument lies on the same side of zero, and the $2r$th finite
difference annihilates that polynomial of degree $p$.  Hence $c_s=0$ there.

The correction is confined to the first $r\times r$ block, and its entries
are specified by the truncated sum:
\begin{equation}\label{eq:polyexplicit-corner}
 (E_n)_{ij}=
 \begin{cases}
 \displaystyle
 \sum_{a=0}^{i-1}\sum_{q=0}^{j-1}
 \alpha_a\alpha_q b_{j-i+a-q}-c_{j-i},
 &1\le i,j\le\min(n,r),\\[6pt]
 0,&i>r\ \text{or}\ j>r.
 \end{cases}
\end{equation}
To verify this support assertion, suppose $i>r$ and $j\le r$.  Every omitted
column-filter term has $q\ge j$, so $j-i-q\le-i<-r$.  Its complete
row-filter sum is an $r$th finite difference taken wholly on the negative
polynomial tail, and hence is zero.  The case $i\le r$, $j>r$ is the
transpose argument; when both indices exceed $r$, neither sum is
truncated.

The correction in \eqref{eq:polyexplicit-corner} is fixed as $n$ increases, so
it affects only the initial boundary state.  The eventual row--column
transfer is the ordinary symmetric $(r,r)$-banded Toeplitz transfer.
The symmetric refinement \eqref{eq:tailsymbound} of
Corollary~\ref{cor:rationaltaildet} therefore gives
$(3^r+1)/2$, proving \eqref{eq:polydistancebound}.
\end{proof}

For $p=1$, the ordinary representative
\eqref{eq:polyexplicit-representative} has
\begin{equation}\label{eq:polyexplicit-pone-parent}
 (c_0,c_1,c_2)=(6x-4,\ 2-4x,\ x),\qquad
 E_n[1{:}2,1{:}2]=
 \begin{bmatrix}
 4-5x&2x-1\\
 2x-1&-x
 \end{bmatrix}
 \quad(n\ge2).
\end{equation}
For $p=3$, its symmetric bulk has lower and upper semibandwidth four and coefficients
\begin{equation}\label{eq:polyexplicit-pthree-parent}
 (c_0,c_1,c_2,c_3,c_4)
 =(70x+32,\ -56x-18,\ 28x,\ 2-8x,\ x).
\end{equation}
The $4\times4$ correction is given by \eqref{eq:polyexplicit-corner}, and the
resulting order bound is $41$.  The two-layer descendant itself has outer
coefficients $1296x$ and $1296x^2$ by
\eqref{eq:polyexplicit-edges}.

The same bound can also be obtained from an ordinary parent that preserves
the order-two descendant.
\begin{proof}
Theorem~\ref{thm:poly-anchored-vapor} constructs a symmetric ordinary parent
with the same order-two descendant as $A_n(x,p)$ and hence the same
determinant.  Its Toeplitz bulk \eqref{eq:poly-anchored-bulk} is
$(r,r)$-banded, and its fixed northwest correction is confined to an
$(r+1)\times(r+1)$ block.  The two polynomial tails satisfy the same
$r$th-difference recurrence, with $m_-=m_+=r$ and
$\alpha_a=\beta_a=(-1)^a\binom ra$.  The symmetric bound in
Corollary~\ref{cor:anchored-tail-recurrence} therefore gives
\eqref{eq:polydistancebound}, including singular specializations.
\end{proof}

\clearpage
\section{Ordinary pentadiagonal transfer}
\label{app:variable-pentadiagonal-transfer}

An ordinary ancestor obtained by vaporisation may have position-dependent
entries.  Let $T=(t_{ij})_{i,j\ge1}$ be such a pentadiagonal array and write its entries by diagonal offset and northwest position:
\[
 t_s[i]:=t_{i+\max(-s,0),\,i+\max(s,0)},
 \qquad -2\le s\le2,\quad i\ge1.
\]
In the following formula put $t_{-1}[0]=t_{-2}[-1]=t_{-2}[0]=0$.
The six-dimensional transfer is explicitly
\begin{equation}\label{eq:penta-variable-ordinary-transfer}
 M_i:=\begin{pmatrix}
 t_0[i]&t_{-1}[i-1]&t_{-2}[i-2]&0&0&0\\
 -t_1[i]&0&0&t_{-1}[i-1]&t_{-2}[i-2]&0\\
 0&-t_1[i]&0&-t_0[i]&0&t_{-2}[i-2]\\
 t_2[i]&0&0&0&0&0\\
 0&t_2[i]&0&0&0&0\\
 0&0&0&t_2[i]&0&0
 \end{pmatrix}.
\end{equation}
If $e_1:=(1,0,0,0,0,0)^{\mathsf T}$, then
\begin{equation}\label{eq:penta-variable-ordinary-determinant}
 \det T[1:n,1:n]=e_1^{\mathsf T}M_n\cdots M_1e_1.
\end{equation}
When $t_2[i]\ne0$, the equation $Tv=0$ advances the
four-component state
\[
 (v_{i+1},v_i,v_{i-1},v_{i-2})^{\mathsf T}.
\]  The two left boundary conditions are
$v_{-1}=v_0=0$, and the two right conditions are
$v_{n+1}=v_{n+2}=0$.  Taking the second compound of the companion
update and multiplying it by $t_2[i]$ gives
\eqref{eq:penta-variable-ordinary-transfer}, in the order
$12,13,14,23,24,34$ of the two-element subsets.  These exterior-square coordinates differ from
the row--column signature coordinates used in
\cite{AlekseyevKhomovskyRecurrences2026}; both constructions give the same
principal determinants.  The product of the
factors $t_2[i]$ cancels the denominators in the companion construction.
Both sides of \eqref{eq:penta-variable-ordinary-determinant} are polynomial
in the entries, so the identity extends when some $t_2[i]$ vanish.

A scalar recurrence is obtained from this transfer by finite linear
algebra, with an explicit nonvanishing condition.  For $n\ge6$ put
\[
 \omega_0:=e_1^{\mathsf T},\qquad
 \omega_j:=e_1^{\mathsf T}M_{n-6+j}\cdots M_{n-5}
 \quad(1\le j\le6),
 \qquad
 \mathcal O_n:=\begin{pmatrix}\omega_0\\\omega_1\\\vdots\\\omega_5\end{pmatrix}.
\]
When $\det \mathcal O_n\ne0$, define
\[
 (\gamma_0(n),\ldots,\gamma_5(n)):=\omega_6\mathcal O_n^{-1}.
\]
Writing $d_n:=\det T[1:n,1:n]$ gives
\begin{equation}\label{eq:penta-variable-observation-recurrence}
 d_n=\sum_{j=0}^5\gamma_j(n)d_{n-6+j}.
\end{equation}
Indeed, all seven observations act on the same state
$M_{n-6}\cdots M_1e_1$.  This also describes a direct symbolic or exact
arithmetic implementation.  At a vanishing observation determinant,
one retains the state update.  A smaller independent set can be used
when it spans the required observation.  The displayed inverse cannot be
used when the observation determinant vanishes.

The coefficient formulas use several consecutive positions.  Substituting
only the current entries into a stationary Toeplitz formula is not valid.
For example, let all four off-diagonals be $1$ and let $t_0[i]=i+3$.
Then
\[
 d_0,\ldots,d_6=1,4,19,107,713,5500,48132.
\]
The stationary coefficients with current main-diagonal value $9$ would
instead give
\[
 9d_5-16d_3+9d_1-d_0=47823.
\]
Thus neither constant scalar coefficients nor their reciprocal relations
follow merely from a position-dependent pentadiagonal support.

The polynomial transfer applies to ordinary pentadiagonal matrices.
For two-layer input, one must first construct an ordinary representative
and identify the determinants of its leading sections.  The explicit
ordinary families in the main text provide those data.

\end{document}